\documentclass[12pt,reqno]{amsart}
\usepackage{amsthm,amsfonts,amssymb,euscript}

\newcommand{\R}{\ensuremath{\mathbb{R}}}
\newcommand{\Z}{\ensuremath{\mathbb{Z}}}

\newcommand{\om}{\ensuremath{\omega}}

\newcommand{\half}{\frac{1}{2}}

\newcommand{\FF}{\ensuremath{\mathcal{F}}}

\newcommand{\m}{\mu}
\newcommand{\x}{\xi}

\theoremstyle{definition}

\theoremstyle{remark}

\theoremstyle{proposition}

\theoremstyle{lemma}

\theoremstyle{corollary}

\numberwithin{equation}{section}
\numberwithin{lemma}{section}
\numberwithin{remark}{section}

\numberwithin{equation}{section}
\begin{document}

\title[Low-regularity solutions of the KP-I equation]
{Weighted low-regularity solutions of the KP-I initial-value problem}
\author{J. Colliander}

\address{University of Toronto}
\thanks{J.C. is supported in part by an NSERC Grant.}
\email{colliand@math.toronto.edu}
\author{A. D. Ionescu}
\address{University of Wisconsin--Madison}
\thanks{A.D.I. is supported in part by an NSF grant and a Packard Fellowship.}
\email{ionescu@math.wisc.edu}
\author{C. E. Kenig}
\thanks{C.E.K. is supported in part by an NSF grant.}
\address{University of Chicago}
\email{cek@math.uchicago.edu}
\author{G. Staffilani}
\thanks{G.S. is supported in part by an NSF grant.}
\address{Massachusetts Institute of Technology}
\email{gigliola@math.mit.edu}

\maketitle
\tableofcontents

\section{Introduction}\label{section1}

In this paper we consider the KP-I initial-value problem
\begin{equation}\label{eq-1}
\begin{cases}
\partial_tu+\partial_x^3u-\partial_x^{-1}\partial_y^2u+\partial_x(u^2/2)=0;\\
u(0)=\phi,
\end{cases}
\end{equation}
on $\mathbb{R}^2_{x,y}\times\mathbb{R}_t$. The dispersion  function for this dispersive equation is 
for  $(\xi,\mu)\in\mathbb{R}\setminus\{0\}\times\mathbb{R}$
\begin{equation*}
\omega(\xi,\mu)=\xi^3+\mu^2/ \xi.
\end{equation*}
In \cite{CoKeSt}  three of the authors studied \eqref{eq-1} with initial data $\phi$ in the space $E\cap P$ defined below.
(See the introduction and the references of  \cite{CoKeSt}  for a discussion of \eqref{eq-1}, its relationship to the corresponding IVP for the KPII equation, and a discussion of the spaces $E$ and $P$ in connection with  \eqref{eq-1}). The main result in \cite{CoKeSt}  is a weak form of local in time well-posedness (Theorem 1 in \cite{CoKeSt} ) for data which is small in  $E\cap P$. Unfortunately, A. Ionescu discovered a counterexample  to the main estimate used in  \cite{CoKeSt}
(Theorem 3 in  \cite{CoKeSt}) to establish Theorem 1. The example exhibits a logarithmic divergence in the estimate, which shows that the proof of Theorem 1 in  \cite{CoKeSt} is incorrect. The same applies to Theorem 2 in  \cite{CoKeSt}. The counterexample is explained in subsection \ref{count} below. Colliander, Kenig and Staffilani are very grateful to Ionescu for pointing  out this mistake and for joining them in this work. Here we obtain a strengthening  of Theorem 1 in  \cite{CoKeSt} which yields the strong form of local in time well-posedness for small data in  $E\cap P$. This is Theorem
\ref{Main1} below. The logarithmic  divergence  is avoided by
introducing new resolution spaces, inspired by those used by
Ionescu-Kenig (\cite{IK1, IK3, IK4}) in works on Benjamin-Ono equation and on the Schr\"odinger map problems. It seems very likely that using the tools developed here, a correct (and similarly strengthened) version of Theorem 2 in  \cite{CoKeSt} could also be obtained. We have felt, however, that this would increase substantially the technicalities in an already very technical paper and we have therefore not pursued this issue.

We conclude by mentioning that our main theorem does not give local well-posedness in  $E\cap P$  for large data; such a result would immediately yield global in time well-posedness.

\subsection{The counterexample}\label{count}
We start this section with some notation and by recalling some spaces of functions  introduced in \cite{CoKeSt}.  We denote the Fourier transform of a function $f(x,y)$ as 
\begin{equation}\label{FT}
\hat f(\xi,\mu)=\FF f(\xi,\mu)=\int_{\R^2}f(x,y)e^{i\langle(x,y)\cdot(\xi,\mu)\rangle}\, dx\, dy.
\end{equation}
Now let $\chi_A$ denote a smooth  characteristic function of the set $A$. 

{\bf{Definition.}} Let $\theta_{0}(s)=\chi_{[-1,1]}(s), \, \, \theta_{m}(s)=
\chi_{[2^{m-1},2^{m}]}(|s|), \, \, m\in \mathbb{N}$. For $(\x,\m)\in
 \mathbb{R}^{2}$
let $\chi_{1}(\x,\m)=\chi_{\{|\x|\geq \half \frac{|\mu|}{|\x|}\}}$, and
$\chi_{2}(\x,\m)=\chi_{\{|\x|< \half{\frac{|\mu|}{|\xi|}}\}}$.
Let $\chi_0(s)=\chi_{\{|s|< 1\}}, \, \chi_j(s)=\chi_{\{2^{j-1}\leq |s|< 2^j\}}$ and $w(\xi,\mu)=
(1+|\xi|+|\mu|/|\xi|)$. We define the space $X_{s,b}$ through the norm
\begin{eqnarray*}
    \label{xs}\|f\|_{X_{s,b}}&=&\sum_{j, m\geq 0}
    2^{jb}\left(\int_{\R^{3}}
    \chi_{j}(\tau-\om(\xi,\m))\chi_{1}(\xi,\m)\theta_{m}(\xi)
    w^{2s}|\hat{f}|^{2}(\x,\m,\tau)d\x d\m d\tau\right)^{\half}\\
    \nonumber&+&\sum_{j, m \geq 0}
    2^{jb}\left(\int_{\R^{3}}
    \chi_{j}(\tau-\om(\xi,\m))\chi_{2}(\xi,\m)\theta_{n}(\m)
    w^{2s}|\hat{f}|^{2}(\x,\m,\tau)d\x d\m d\tau\right)^{\half}.
    \end{eqnarray*}
We also define the space
\begin{equation*}
    \label{yrs}Y_{s,r,b}=\{f : tf \in X_{s,b}, \mbox{ and }  yf \in
    X_{r,b}\},
    \end{equation*}
    and the spaces
\begin{equation*}
    \label{zsb}Z_{s,b}=X_{s,b}\cap Y_{s,1-s,b}, \, \, \mbox{ and} \, \, 
    Z_{1-\epsilon}=Z_{1-\epsilon,\half}.
    \end{equation*}
    \label{spaces}
We  recall here the statement of Theorem 3 in \cite{CoKeSt}:

\newtheorem*{Old1}{Theorem}
\begin{Old1}\label{Old1}
Assume $0 < \epsilon_{0}<\frac{1}{8}$. Then for any 
$\frac{1}{4}<\epsilon<1$, we have
\begin{eqnarray}
\label{Xbilinear}
\|\partial_{x}(uv)\|_{X_{1-\epsilon_{0},-\half}}&\leq& C
\|u\|_{X_{1-\epsilon_{0},\half}}(\|v\|_{X_{1-\epsilon_{0},\half}}+
\|v\|_{X_{1-\epsilon_{0},\half}}^{1-\epsilon }
\|v\|_{Y_{1-\epsilon_{0},-\epsilon_{0}, \half}}^{\epsilon })\\\nonumber
&+&C
\|v\|_{X_{1-\epsilon_{0},\half}}(\|u\|_{X_{1-\epsilon_{0},\half}}+
\|u\|_{X_{1-\epsilon_{0},\half}}^{1-\epsilon }
\|u\|_{Y_{1-\epsilon_{0},-\epsilon_{0}, \half}}^{\epsilon })
\end{eqnarray}
\end{Old1}

This theorem unfortunately cannot hold since the following counterexample 
shows a logarithmic divergence.
Let $\psi:\mathbb{R}\to[0,1]$ denote a smooth function
supported in the interval $[-2,2]$ and equal to $1$ in the
interval $[-1,1]$. Assume $N\gg 1$ is very large,
$\omega(\xi,\mu)=\xi^3+\mu^2/\xi$, and define
\begin{equation}\label{w1}
\widehat{u}(\xi,\mu,\tau)=\psi(\xi-N)\psi((\mu-\sqrt{3}\xi^2)/\xi)\psi(\tau-\omega(\xi,\mu)),
\end{equation}
and
\begin{equation}\label{w2}
\widehat{v}(\xi,\mu,\tau)=\psi(\xi-4)\psi((\mu+\sqrt{3}\xi^2)/\xi)\psi(\tau-\omega(\xi,\mu)).
\end{equation}
Notice that in the definition \eqref{w1} $|\mu-\sqrt{3}N^2|\leq
CN$ and in the definition \eqref{w2} the variable $\xi$ is about
$1$ (bounded away from $0$). The functions
$\mu\to\psi((\mu-\sqrt{3}\xi^2)/\xi)$ in \eqref{w1} and
$\mu\to\psi((\mu+\sqrt{3}\xi^2)/\xi)$ in \eqref{w2} are
essentially the characteristic functions of the intervals
$[\sqrt{3}N^2-N,\sqrt{3}N^2+N]$ and
$[-16\sqrt{3}-1,-16\sqrt{3}+1]$ respectively. The precise formulas
$\psi((\mu-\sqrt{3}\xi^2)/\xi)$ and
$\psi((\mu+\sqrt{3}\xi^2)/\xi)$ are convenient for the nonlinear
change of variables \eqref{w11}. Then
\begin{equation}\label{w3}
\begin{split}
&||u||_{X_{1-\epsilon_0,1/2}}\approx N^{1-\epsilon_0}N^{1/2},\,||u||_{Y_{1-\epsilon_0,-\epsilon_0,1/2}}
\approx N^{1-\epsilon_0}N^{1/2},\\
&||v||_{X_{1-\epsilon_0,1/2}}\approx 1,\,||v||_{Y_{1-\epsilon_0,-\epsilon_0,1/2}}\approx 1.
\end{split}
\end{equation}
So the right-hand side in \eqref{Xbilinear} is
\begin{equation}\label{w4}
RHS\approx N^{1-\epsilon_0}N^{1/2}.
\end{equation}
We look now at the left-hand side of \eqref{Xbilinear}: the function $\widehat{u}\ast\widehat{v}$ is supported in the 
set $\{(\xi,\mu,\tau):|\xi-N|\leq C\text{ and }|\mu-\sqrt{3}N^2|\leq CN\}$. So,
\begin{equation*}
||\partial_x(uv)||_{X_{1-\epsilon_0,-1/2}}\approx N\cdot N^{1-\epsilon_0}
\sum_{j\geq 0}2^{-j/2}||(\widehat{u}\ast\widehat{v})(\xi,\mu,\tau)
\chi_j(\tau-\omega(\xi,\mu))||_{L^2_{\xi,\mu,\tau}},
\end{equation*}
where $\chi_j$ is the characteristic function of the set $\{s:|s|\in[2^{j-1},2^{j+1}]\}$. Using \eqref{w4}, it would follow from \eqref {Xbilinear} that
\begin{equation}\label{w5}
\sum_{j\geq 0}2^{-j/2}||(\widehat{u}\ast\widehat{v})(\xi,\mu,\tau)\chi_j(\tau-\omega(\xi,\mu))||_{L^2_{\xi,\mu,\tau}}\leq CN^{-1/2}.
\end{equation}
We show now that if $100\leq 2^j\leq N^{1/10}$ then
\begin{equation}\label{w6}
||(\widehat{u}\ast\widehat{v})(\xi,\mu,\tau)\chi_j(\tau-\omega(\xi,\mu))||_{L^2_{\xi,\mu,\tau}}
\geq c2^{j/2}N^{-1/2}.
\end{equation}
So the bound \eqref{w5} would fail by $\ln N$ since the sum in $j$ has $\approx \ln N$ terms.
To prove \eqref{w6}, by duality, it suffices to prove that if $100\leq 2^j\leq N^{1/10}$
\begin{equation}\label{w7}
\begin{split}
\int_{\mathbb{R}}\int_{\mathbb{R}}\int_{\mathbb{R}}&(\widehat{u}\ast\widehat{v})(\xi,\mu,\tau)\mathbf{1}_{[N-10,N+10]}(\xi)\\
&\mathbf{1}_{[\sqrt{3}N^2-100N,\sqrt{3}N^2+100N]}(\mu)\chi_j(\tau-\omega(\xi,\mu))\,d\xi d\mu d\tau\geq c2^{j},
\end{split}
\end{equation}
where $\mathbf{1}_A$ denotes the characteristic function of the set $A$. We substitute the formulas \eqref{w1} and \eqref{w2}; the left-hand side of \eqref{w7} becomes
\begin{equation}\label{w8}
\begin{split}
\int_{\mathbb{R}^6}&\psi(\xi_1-N)\psi((\mu_1-\sqrt{3}\xi_1^2)/\xi_1)\psi(\tau_1-\omega(\xi_1,\mu_1))\psi(\xi-\xi_1-4)\\
&\psi((\mu-\mu_1+\sqrt{3}(\xi-\xi_1)^2)/(\xi-\xi_1))\psi(\tau-\tau_1-\omega(\xi-\xi_1,\mu-\mu_1))\\
&\mathbf{1}_{[N-10,N+10]}(\xi)\mathbf{1}_{[\sqrt{3}N^2-100N,\sqrt{3}N^2+100N]}(\mu)\chi_j(\tau-\omega(\xi,\mu))\,d\xi_1 d\mu_1 d\tau_1d\xi d\mu d\tau.
\end{split}
\end{equation}
In this expression we make the change of variables
\begin{equation*}
\begin{split}
&\xi_1=\xi_1,\,\mu_1=\mu_1,\,\xi=\xi_1+\xi_2,\,\mu=\mu_1+\mu_2\\
&\tau_1=\mu_1+\omega(\xi_1,\mu_1),\,\tau=\mu_2+\mu_1+\omega(\xi_1,\mu_1)+\omega(\xi_2,\mu_2).
\end{split}
\end{equation*}
Then we notice that

\begin{equation*}
\begin{split}
&\psi(\xi_1-N)\psi(\xi_2-4)\mathbf{1}_{[N-10,N+10]}(\xi_1+\xi_2)=\psi(\xi_1-N)\psi(\xi_2-4);\\
&\psi((\mu_1-\sqrt{3}\xi_1^2)/\xi_1)\psi((\mu_2+\sqrt{3}\xi_2^2)/\xi_2)\mathbf{1}_{[\sqrt{3}N^2-100N,\sqrt{3}N^2+100N]}(\mu_1+\mu_2)\\
&=\psi((\mu_1-\sqrt{3}\xi_1^2)/\xi_1)\psi((\mu_2+\sqrt{3}\xi_2^2)/\xi_2)\text{ if }\xi_1\in[N-2,N+2]\text{ and }\xi_2\in[2,6].
\end{split}
\end{equation*}
Thus the expression in \eqref{w8} becomes
\begin{equation}\label{w10}
\begin{split}
\int_{\mathbb{R}^6}&\psi(\xi_1-N)\psi((\mu_1-\sqrt{3}\xi_1^2)/\xi_1)\psi(\mu_1)\psi(\xi_2-4)\psi((\mu_2+\sqrt{3}\xi_2^2)/\xi_2)\psi(\mu_2)\\
&\chi_j(\mu_1+\mu_2+\Omega(\xi_1,\mu_1,\xi_2,\mu_2))\,d\xi_1 d\mu_1 d\mu_1d\xi_2 d\mu_2 d\mu_2,
\end{split}
\end{equation}
where
\begin{equation}\label{w20}
\begin{split}
\Omega(\xi_1,\mu_1,\xi_2,\mu_2)&=\omega(\xi_1,\mu_1)+\omega(\xi_2,\mu_2)-\omega(\xi_1+\xi_2,\mu_1+\mu_2)\\
&=-\frac{\xi_1\xi_2}{\xi_1+\xi_2}\Big[(\sqrt{3}\xi_1+\sqrt{3}\xi_2)^2-(\mu_1/\xi_1-\mu_2/\xi_2)^2\Big].
\end{split}
\end{equation}
We make now the nonlinear change of variables
\begin{equation}\label{w11}
\mu_1=\sqrt{3}\xi_1^2+\beta_1\xi_1,\,\mu_2=-\sqrt{3}\xi_2^2+\beta_2\xi_2,
\end{equation}
with $d\mu_1d\mu_2=\xi_1\xi_2d\beta_1d\beta_2\approx Nd\beta_1d\beta_2$. 
The expression in \eqref{w10} is bounded from below by
\begin{equation}\label{w15}
\begin{split}
(N/2)\int_{\mathbb{R}^6}&\psi(\xi_1-N)\psi(\beta_1)\psi(\mu_1)\psi(\xi_2-4)\psi(\beta_2)\psi(\mu_2)\\
&\chi_j(\mu_1+\mu_2+\widetilde{\Omega}(\xi_1,\beta_1,\xi_2,\beta_2))\,d\xi_1 d\beta_1 d\mu_1d\xi_2 d\beta_2 d\mu_2,
\end{split}
\end{equation}
where, by \eqref{w20},
\begin{equation}\label{w16}
\widetilde{\Omega}(\xi_1,\beta_1,\xi_2,\beta_2)=(\beta_1-\beta_2)\xi_1\xi_2\Big(2\sqrt{3}+\frac{\beta_1-\beta_2}{\xi_1+\xi_2}\Big).
\end{equation}
It follows from \eqref{w16} that if $\xi_1\in[N-1/100,N+1/100]$,

$\xi_2\in[4-1/100,4+1/100]$,

$|\beta_1-\beta_2|\in[[1/(8\sqrt{3})-1/100]2^j/N,[1/(8\sqrt{3})+1/100]2^j/N]$,

$\mu_1,\mu_2\in[-2,2]$, and $2^j\in[100,N^{1/10}]$ then

\begin{equation*}
\chi_j(\mu_1+\mu_2+\widetilde{\Omega}(\xi_1,\beta_1,\xi_2,\beta_2))=1.
\end{equation*}
Thus the only nontrivial restriction in the integral in \eqref{w15} is
\begin{equation*}
|\beta_1-\beta_2|\in[[1/(8\sqrt{3})-1/100]2^j/N,[1/(8\sqrt{3})+1/100]2^j/N],
\end{equation*}
which shows that this integral is bounded from below by $cN\cdot
2^j/N=c2^j$. This is the bound \eqref{w7}, which implies \eqref{w6}.

\subsection{ The main theorem}
In this section we introduce again the spaces of functions $E$ and $P$ already defined in 
\cite{CoKeSt} and state the main result that replaces Theorem 1 in \cite{CoKeSt}.
We define the energy
space $E$,
\begin{equation}\label{eq-2}
E=\{\phi:\mathbb{R}\times\mathbb{R}\to\mathbb{C}:\,\|\phi\|_{E}
:=\|\widehat{\phi}(\xi,\mu)\cdot (1+|\xi|+|\mu/ \xi| )\|_{L^2_{\xi,\mu}}<\infty\},
\end{equation}
and the weighted space $P$,
\begin{equation}\label{eq-3}
P=\{\phi:\mathbb{R}\times\mathbb{R}\to\mathbb{C}:\,\|\phi\|_{P}
:=\| (y+i)\cdot \phi\|_{L^2}<\infty\}.
\end{equation}
In Section \ref{section2}, see \eqref{sp5}, we will define a Banach space $F\hookrightarrow C(\R:E\cap P)$; let
\begin{equation*}
F_1=\{u\in C([-1,1]:E\cap P):\, \|u\|_{F_1}=\inf_{\widetilde{u}=u\text{ on }\R^2\times[-1,1]}\|\widetilde{u}\|_F<\infty \}. 
\end{equation*}
For any Banach space $V$ and $r>0$ let $B(r,V)$ denote the open
ball $\{v\in V:||v||_V<r\}$. Our main theorem concerns local
well-posedness of
the KP-I initial value problem \eqref{eq-1} for small data in $E\cap P$.
\newtheorem{Main1}{Theorem}[section]
\begin{Main1}\label{Main1}
There are $\overline{r},\overline{R}\in(0,1]$, $\overline{r}\leq \overline{R}$, with the property that for any
$\phi\in B({\overline{r}},E\cap P)$ there is a unique $u\in
B({\overline{R}},F_1)$ such that
\begin{equation*}
\begin{cases}
(\partial_t+\partial_x^3-\partial_x^{-1}\partial_y^2)u+\partial_x(u^2/2)=0
\text{ in }C((-1,1):H^{-2});\\
u(0)=\phi.
\end{cases}
\end{equation*}
In addition, the mapping $\phi\to u$ is Lipschitz continuous from
$B(\overline{r},E\cap  P)$ to $B({\overline{R}},F_1)$.
\end{Main1}

The rest of the paper is organized as follows: in section \ref{section2} we define the main normed spaces $X_k$, $Y_k$, $V_k$, $W_k$, $F$, and $N$, and prove some of their basic properties. As explained in subsection \ref{count}, the use of standard $X^{s,b}$-type spaces seems to lead inevitably to logarithmic divergences in the modulation variable. To avoid these logarithmic divergences we work with high-frequency spaces that have two components: an $X^{s,b}$-type component measured  in the frequency space (see the space $X_k$) and a normalized $L^1_yL^2_{x,t}$ component measured in the physical space (see the space $Y_k$). As in \cite{IK1}, \cite{IK3}, and \cite{IK4}, for the physical space component we use a suitable normalization of the local smoothing space $L^1_yL^2_{x,t}$.  

In section \ref{linear} we prove two linear estimates. In section \ref{proofthm} we prove Theorem \ref{Main1}, using a direct perturbative argument in the Banach space $F_1$, and assuming the dyadic bilinear estimates \eqref{BIG1} and \eqref{BIG2}. The remaining sections are concerned with the proofs of \eqref{BIG1} and \eqref{BIG2}: in sections \ref{prop} and \ref{L2bi} we prove preliminary linear estimates and an $L^2$ bilinear estimate. In sections \ref{bilinear2}, \ref{bilinear1}, and \ref{bilinear3} we prove the dyadic bilinear estimate \eqref{BIG1}. In section \ref{bilinear4} we prove the dyadic bilinear estimate \eqref{BIG2}.

\section{The resolution spaces}\label{section2}

In this section we define the main normed spaces we will use in
the rest of the paper, and prove some of their basic properties. Let
$\mathbb{Z}_+=\mathbb{Z}\cap[0,\infty)$. For $k\in \mathbb{Z}$ let
$I_k=\{\xi:|\xi|\in[2^{k-1},2^{k+1}]\}$, $\widetilde{I}_k=I_k$ if
$k\geq 1$, $\widetilde{I}_k=[-2,2]$ if $k=0$, and
$\widetilde{I}_k=\emptyset$ if $k\leq -1$. Let
$\mu_0:\mathbb{R}\to[0,1]$ denote an even smooth function
supported in $[-8/5,8/5]$ and equal to $1$ in $[-5/4,5/4]$. For
$k\in\mathbb{Z}$ let
$\chi_k(\xi)=\mu_0(\xi/2^k)-\mu_0(\xi/2^{k-1})$. Let $\mu_k=\chi_k$ for $k\in\Z\cap[1,\infty)$ and $\mu_k=0$ for $k\in \Z\cap(-\infty,-1]$. Let
\begin{equation*}
\chi_{[k_1,k_2]}=\sum_{k=k_1}^{k_2}\chi_k\text{ for any }k_1\leq
k_2\in\mathbb{Z}.
\end{equation*}
and, for $j\in\Z$,
\begin{equation*}
\mu_{\leq j}=\sum_{j'=-\infty}^j\mu_{j'}\text{ and }\mu_{\geq j}=\sum_{j'=j}^\infty \mu_{j'}.
\end{equation*}

For $(\xi,\mu)\in\mathbb{R}\setminus\{0\}\times\mathbb{R}$ let
\begin{equation}\label{omega}
\omega(\xi,\mu)=\xi^3+\mu^2/ \xi.
\end{equation}
We define the relevant KP-I region\footnote{The main difficulties of the KP-I problem, including the counterexample of subsection \ref{count}, are caused by functions with Fourier support  in this region.}  
\begin{equation}\label{sp8}
R_{KP-I}=\{(\xi,\mu,\tau)\in\mathbb{R}^3:|\xi|\geq 1,\,|\mu|\in[
|\xi|^2/2^{20},2^{20}\cdot |\xi|^2],\,|\tau-\omega(\xi,\mu)|\leq
|\xi|\}.
\end{equation}

For $k\in\Z$ let $k_+=\max(k,0)$. We define the normed spaces $X_k$,
\begin{equation}\label{sp7}
\begin{split}
X_{k}=\{f\in L^2(I_k&\times\mathbb{R}\times\R):\|f\|_{X_k}=\sum_{j=
0}^{2k_+-1}2^{j/2}\|\mu_j(\tau-\omega(\xi,\mu)) \cdot f\|_{L^2}\\
&+\Big[\sum_{j\geq 2k_+}2^{2j-2k_+}\|\mu_j(\tau-\omega(\xi,\mu)) \cdot f\|_{L^2}^2\Big]^{1/2}<\infty\}.
\end{split}
\end{equation}
Notice that
\begin{equation}\label{sp7.7}
\|(\tau-\omega(\xi,\mu)+i)^{-1}\cdot \mu_{\geq J}(\tau-\omega(\xi,\mu))\cdot f\|_{X_k}\leq C(2^{-J/2}+2^{-2k_+/2})\cdot \|f\|_{L^2},
\end{equation}
for any $f\in L^2(\R^3)$ supported in $I_k\times\R\times\R$ and $J\in\Z_+$. 

The spaces $X_{k}$ are not sufficient for a fixed-point argument, due to
various logarithmic divergences.
For $k\geq 100$ we also define the normed spaces $Y_k$,
\begin{equation}\label{sp2}
\begin{split}
Y_{k}=&\{f\in L^2(\mathbb{R}^3): f\text{ supported in }
R_{KP-I}\cap I_k\times \mathbb{R}\times\mathbb{R}\text{ and }\\
&\| f \|_{Y_{k}}=2^{-k/2}\|\mathcal{F}^{-1}
[(\tau-\omega(\xi,\mu)+i)\cdot f(\xi,\mu,\tau)]\|_{L^1_yL^2_{x,t}}<\infty \}.
\end{split}
\end{equation}
For simplicity of notation, we define $Y_k=\{0\}$ for $k\leq 99$. Then we define the normed spaces $X_k+Y_k$, $k\in\Z$,
\begin{equation*}
\begin{split}
X_k+Y_{k}=&\{f\in L^2(\mathbb{R}^3): f\text{ supported in }I_k\times \mathbb{R}\times\mathbb{R}\text{ and }\\
&\| f \|_{X_k+Y_{k}}=\inf_{f=f_1+f_2}\|f_1\|_{X_k}+\|f_2\|_{Y_k}<\infty \}.
\end{split}
\end{equation*}

For
$k\in\mathbb{Z}$ we define the normed spaces $V_k$
\begin{equation}\label{sp1}
\begin{split}
V_{k}=\{f\in L^2&(\mathbb{R}^3): f\text{ supported in }
I_k\times \mathbb{R}\times\mathbb{R}\text{ and }\\
&\|f\|_{V_k}=\|f\cdot (1+2^k+i\mu/ 2^k)\|_{X_k+Y_k}<\infty\},
\end{split}
\end{equation}
and the normed spaces $W_k$,
\begin{equation}\label{sp4}
\begin{split}
W_{k}=\{f\in L^2(\mathbb{R}^3)&: f\text{ supported in }
I_k\times \mathbb{R}\times\mathbb{R}\text{ and }\\
&\|f\|_{W_k}=\|(\partial_\mu+I)f\|_{X_k+Y_k}<\infty\}.
\end{split}
\end{equation}
We define the (global) normed space $F=F(\mathbb{R}^3)$,
\begin{equation}\label{sp5}
\begin{split}
F= \{u\in L^2&(\mathbb{R}^3):\,u \text{ supported in }\R^2\times[-2,2]\text{ and }\\
&\|u\|_{F}^2=\sum_{k\in\mathbb{Z}}\|\chi_{k}(\xi)\cdot \mathcal{F}(u)\|_{V_k\cap W_k}^2<\infty\},
\end{split}
\end{equation}
and the normed space $N=N(\R^3)$, 
\begin{equation}\label{sp6}
\begin{split}
N=&\{ u\in C(\R:H^{-2}(\R^3)):\\
&\|u\|_{N}^2
=\sum_{k\in\mathbb{Z}}\|\chi_{k}(\xi)(\tau-\omega(\xi,\mu)+i)^{-1}
\cdot \mathcal{F}(u)\|_{V_k\cap W_k}^2<\infty\}.
\end{split}
\end{equation}

We start with a simple lemma concerning basic properties of our normed spaces.

\newtheorem{Lemmaa1}{Lemma}[section]
\begin{Lemmaa1}\label{Lemmaa1}
(a) If $m:\mathbb{R}\to\mathbb{C}$,
$m':\mathbb{R}^2\to\mathbb{C}$, $k\in\mathbb{Z}$, and $f$ is supported in $I_k\times\mathbb{R}\times\mathbb{R}$
then
\begin{equation}\label{lb1}
\begin{cases}
&||m(\mu)\cdot f||_{X_k+Y_k}\leq C||\mathcal{F}^{-1}(m)
||_{L^1(\mathbb{R})}\cdot ||f||_{X_k+Y_k};\\
&||m(\mu)\cdot f||_{V_k\cap W_k}\leq C(||\mathcal{F}^{-1}(m)
||_{L^1(\mathbb{R})}+\|\partial_\mu m\|_{L^\infty(\mathbb{R})})\cdot ||f||_{V_k\cap W_k};\\
&||m'(\xi,\tau)\cdot f||_{X_k+Y_k}\leq
C||m'||_{L^\infty(\mathbb{R}^2)}||f||_{X_k+Y_k};\\
&||m'(\xi,\tau)\cdot f||_{V_k\cap W_k}\leq
C||m'||_{L^\infty(\mathbb{R}^2)}||f||_{V_k\cap W_k}.
\end{cases}
\end{equation}

(b) If $k\in\mathbb{Z}$, $j\geq 0$, and $f_k\in X_k+Y_k$ then
\begin{equation}\label{lb2}
||\mu_j(\tau-\omega(\xi,\mu))\cdot f_k||_{X_k}\leq C||f_k||_{X_k+Y_k}.
\end{equation}
In particular, for any $J\in\Z_+$,
\begin{equation}\label{sp7.8}
\|\eta_{\geq J}(\tau-\omega(\xi,\mu))\cdot f_k\|_{L^2}\leq C2^{-J/2}(2^{(J-2k_+)/2}+1)^{-1}\cdot \|f_k\|_{X_k+Y_k},
\end{equation}
and
\begin{equation}\label{lb9}
||f_k||_{X_k}\leq C(1+k_+)||f_k||_{X_k+Y_k}.
\end{equation}

(c) If $k\geq 0$, $j\in[0,k]\cap\mathbb{Z}$, and $f$ is supported in the set
$$\{(\xi,\mu,\tau)\in\mathbb{R}^3:\xi\in I_k,\,|\mu|\in[2^{2k-100},2^{2k+100}]\},$$ then
\begin{equation}\label{lb7}
||\mathcal{F}^{-1}[\eta_{\leq j}(\tau-\omega(\xi,\mu))\cdot f]
||_{L^1_yL^2_{x,t}}\leq
C||\mathcal{F}^{-1}(f)||_{L^1_yL^2_{x,t}}.
\end{equation}
\end{Lemmaa1}

\begin{proof}[Proof of Lemma \ref{Lemmaa1}] Part (a) follows directly from the definitions. 

For part (b), we may assume $k\geq 100$, $f_k\in Y_k$, so $f_k$ can be written as
\begin{equation}\label{me1}
\begin{split}
f_k(\xi,\mu,\tau)=2^{k/2}&\mathbf{1}_{I_k}(\xi)\chi_{[2k-30,2k+30]}(\mu)\eta_{\leq k+1}(\tau-\omega(\xi,\mu))\\
&\times (\tau-\omega(\xi,\mu)+i)^{-1}\cdot \int_{\mathbb{R}} e^{-iy\cdot\mu}g_k(y,\xi,\tau)\,dy,
\end{split}
\end{equation}
with
\begin{equation}\label{me2}
\| f_k\|_{Y_k}=C\|g_k\|_{L^1_yL^2_{\xi,\tau}}.
\end{equation}
The bound \eqref{lb2} follows easily since $|\{\mu:|\tau-\omega(\xi,\mu)|\leq 2^{j+1}\}|\leq C2^{j-k}$ whenever $|\xi|\approx 2^k$, $|\mu|\approx 2^{2k}$, and $j\leq k+C$. 

For part (c), using Plancherel theorem, it suffices to prove that 
\begin{equation}\label{pr42}
\Big|\Big|\int_{\mathbb{R}}e^{iy\cdot \mu}\chi_{[k-1,k+1]}(\xi)\chi_{[2k-110,2k+110]}(\mu)\eta_{\leq j}(\tau-\omega(\xi,\mu))\,d\mu \Big|\Big|_{L^1_yL^\infty_{\xi,\tau}}\leq C.
\end{equation}
In proving \eqref{pr42} we may assume $k\geq 100$. Then the function in the left-hand side of \eqref{pr42} is not zero only if $|\tau-\xi^3|\approx 2^{3k}$. Simple estimates using integration by parts show that
\begin{equation*}
\Big|\int_{\mathbb{R}}e^{iy\cdot \mu}\chi_{[k-1,k+1]}(\xi)\chi_{[2k-110,2k+110]}(\mu)\eta_{\leq j}(\tau-\omega(\xi,\mu))\,d\mu\Big|\leq C\frac{2^{j-k}}{1+(2^{j-k}y)^2}
\end{equation*}
if $|\tau-\xi^3|\approx 2^{3k}$, which suffices to prove \eqref{pr42}.
\end{proof}

We show now that $F\hookrightarrow C(\R:E\cap P)$.

\newtheorem{Lemmaa2}[Lemmaa1]{Lemma}
\begin{Lemmaa2}\label{Lemmaa2}
If $u\in F$ then
\begin{equation}\label{sd1}
\sup_{t\in\R} \|u(.,.,t)\|_{E\cap P}\leq C\|u\|_{F}.
\end{equation}
Thus $F\hookrightarrow C(\R:E\cap P)$.
\end{Lemmaa2}

\begin{proof}[Proof of  Lemma \ref{Lemmaa2}] Let $f_k=\chi_k(\xi)\cdot \mathcal{F}(u)$, $k\in\Z$. In view of the definition \eqref{sp5}, it suffices to prove that for any $t\in \R$ and $k\in\Z$
\begin{equation*}
\|\mathcal{F}^{-1}(f_k)(.,.,t)\|_{E\cap P}\leq C\|f_k\|_{V_k\cap W_k}.
\end{equation*}
In view of the last bound in \eqref{lb1}, we may assume $t=0$. Thus it suffices to prove that if $k\in\Z$ and $f_k\in Z_k$ then
\begin{equation}\label{sd2}
\Big|\Big|\int_{\R^3}f_k(\xi,\mu,\tau)e^{ix\cdot \xi}e^{iy\cdot \mu}\,d\xi d\mu d\tau\Big|\Big|_{E\cap P}\leq C\|f_k\|_{V_k\cap W_k}.
\end{equation}

We show first that
\begin{equation}\label{sd3}
\Big|\Big|\int_{\R^3}f_k(\xi,\mu,\tau)e^{ix\cdot \xi}e^{iy\cdot \mu}\,d\xi d\mu d\tau\Big|\Big|_{E}\leq C\|f_k\|_{V_k}.
\end{equation}
Using the definition \eqref{eq-2}, it suffices to prove that
\begin{equation}\label{sd4}
\Big|\Big|(1+2^k+i\mu/ 2^k )\cdot \int_{\R}f_k(\xi,\mu,\tau)\,d\tau\Big|\Big|_{L^2_{\xi,\mu}}\leq C\|f_k\|_{V_k}.
\end{equation}
Using the definition \eqref{sp1}, it suffices to prove that 
\begin{equation}\label{sd4.1}
\Big|\Big|\int_{\R}f_k(\xi,\mu,\tau)\,d\tau\Big|\Big|_{L^2_{\xi,\mu}}\leq C\|f_k\|_{X_k+Y_k}.
\end{equation}
Assume first that $f_k\in X_k$ and write $f_k=\sum_{j\geq 0}f_k\cdot \eta_j(\tau-\omega(\xi,\mu))=\sum_{j\geq 0}f_{k,j}$. The left-hand side of \eqref{sd4.1} is dominated by
\begin{equation*}
\begin{split}
&C\sum_{j\geq 0}\Big|\Big| \int_{\R}|f_{k,j}(\xi,\mu,\tau)|\,d\tau\Big|\Big|_{L^2_{\xi,\mu}}\leq C\sum_{j\geq 0}2^{j/2}\|f_{k,j}\|_{L^2_{\xi,\mu,\tau}}\leq C\|f_k\|_{X_k},
\end{split}
\end{equation*}
as desired. Asssume now that $f_k\in Y_k$ (so $k\geq 100$) and write $f_k$ as in \eqref{me1}. With $g_k$ as in \eqref{me1} and \eqref{me2}, the left-hand side of \eqref{sd4.1} is dominated by
\begin{equation}\label{sd7}
\begin{split}
&C2^{k/2}\Big|\Big|\mathbf{1}_{I_k}(\xi)\chi_{[2k-30,2k+30]}(\mu)\int_{\mathbb{R}\times\R} e^{-i y\cdot \mu}\frac{\eta_{\leq k+1}(\tau-\omega(\xi,\mu))}{\tau-\omega(\xi,\mu)+i}\cdot g_k(y,\xi,\tau)\,dyd\tau\Big|\Big|_{L^2_{\xi,\mu}}.
\end{split} 
\end{equation}
We define the partial Hilbert transform operator 
\begin{equation*}
\mathcal{L}_k(g)(y,\xi,\nu)=\int_{\R}g(y,\xi,\tau)\cdot \eta_{\leq k+1}(\tau-\nu)\cdot (\tau-\nu+i)^{-1}\,d\tau.
\end{equation*}
Using the Minkowski inequality, the expression in \eqref{sd7} is dominated by
\begin{equation*}
C2^{k/2}\int_{\R}\Big|\Big|\mathbf{1}_{I_k}(\xi)\chi_{[2k-30,2k+30]}(\mu)\cdot \mathcal{L}_k(g_k)(y,\xi,\omega(\xi,\mu))\Big|\Big|_{L^2_{\xi,\mu}}\,dy.
\end{equation*}
A simple change of variables shows that this is dominated by
\begin{equation*}
C\int_{\R}\Big|\Big|\mathbf{1}_{I_k}(\xi)\mathcal{L}_k(g_k)(y,\xi,\nu)\Big|\Big|_{L^2_{\xi,\nu}}\,dy,
\end{equation*}
and the bound \eqref{sd4.1} follows from \eqref{me2} and the estimate $\|\mathcal{L}_k(g)(y,\xi,\nu)\|_{L^2_{\xi,\nu}}\leq C\|g(y,\xi,\tau)\|_{L^2_{\xi,\tau}}$.

We show now  that
\begin{equation}\label{sd10}
\Big|\Big|\int_{\R^3}f_k(\xi,\mu,\tau)e^{ix\cdot \xi}e^{iy\cdot \mu}\,d\xi d\mu d\tau\Big|\Big|_{P}\leq C\|f_k\|_{W_k}.
\end{equation}
Using the definition \eqref{eq-3} and Plancherel theorem, it suffices to prove that
\begin{equation*}
\Big|\Big|\int_{\R}(\partial_\mu+I)f_k(\xi,\mu,\tau)\,d\tau\Big|\Big|_{L^2_{\xi,\mu}}\leq C\|f_k\|_{W_k},
\end{equation*}
which follows from \eqref{sd4.1}. The bound \eqref{sd2} follows from \eqref{sd3} and \eqref{sd10}.
\end{proof}

\section{Linear estimates}\label{linear}

In this section we prove two linear estimates. For $\phi \in L^2(\mathbb{R}^2)$ let $W\phi \in C(\mathbb{R}:L^2_{x,y})$ denote
the solution of the free KP-I evolution given by
\begin{equation}\label{ni1}
W\phi(x,y,t)=C\int_{\R^2}e^{ix\cdot \xi}e^{iy\cdot \mu}e^{it\omega(\xi,\mu)}\widehat{\phi}(\xi,\mu)\,d\xi d\mu,
\end{equation}
where $\omega(\xi,\mu)$ is defined in \eqref{omega}. Let
$\psi=\widehat{\eta_0}\in\mathcal{S}(\mathbb{R})$.

\newtheorem{Lemmab1}{Proposition}[section]
\begin{Lemmab1}\label{Lemmab1}
If $\phi \in E\cap P$ then
\begin{equation*}
||\eta_0(t)\cdot W\phi||_{F}\leq C||\phi||_{E\cap P}.
\end{equation*}
\end{Lemmab1}

\begin{proof}[Proof of Proposition \ref{Lemmab1}] A straightforward
computation shows that
\begin{equation}\label{ni3}
\mathcal{F}[\eta_0(t)\cdot W\phi](\xi,\mu,\tau)=
\widehat{\phi}(\xi,\mu)\widehat{\eta_0}(\tau-\omega(\xi,\mu)).
\end{equation}
Then, directly from the definitions,
\begin{equation*}
\begin{split}
|&|\eta_0(t)\cdot W\phi ||_{F}^2
\leq C\sum_{k\in\mathbb{Z}}||\chi_k(\xi)\cdot \widehat{\phi}(\xi,\mu)
\cdot \psi(\tau-\omega(\xi,\mu))||_{V_k\cap W_k}^2\\
&\leq
C\sum_{k\in\mathbb{Z}}||\chi_k(\xi)\widehat{\phi}(\xi,\mu)\cdot
(1+2^k+|\mu|/2^k)||_{L^2_{\xi,\mu}}^2+C\sum_{k\in\mathbb{Z}}||\chi_k(\xi)
(\partial_\mu\widehat{\phi})(\xi,\mu)||_{L^2_{\xi,\mu}}^2\\
&\leq C(||\phi||_{E}^2+||\phi ||_{P}^2),
\end{split}
\end{equation*}
as desired.
\end{proof}

\newtheorem{Lemmab3}[Lemmab1]{Proposition}
\begin{Lemmab3}\label{Lemmab3}
If $u\in N$
then
\begin{equation*}
\Big|\Big|\eta_0(t)\cdot \int_0^t[Wu(s)](t-s)\,ds\Big|\Big|_{F}\leq C||u||_{N}.
\end{equation*}
\end{Lemmab3}

\begin{proof}[Proof of Proposition \ref{Lemmab3}] A direct computation shows that
\begin{equation}\label{ni2}
\begin{split}
\mathcal{F}&\Big[\eta_0(t)\cdot \int_0^t[Wu(s)](t-s)ds\Big](\xi,\mu,\tau)\\
&=C\int_\mathbb{R}\mathcal{F}(u)(\xi,\mu,\tau')\cdot \frac{\widehat{\eta_0}(\tau-\tau')-\widehat{\eta_0}(\tau-\omega(\xi,\mu))}{\tau'-\omega(\xi,\mu )}d\tau'.
\end{split}
\end{equation}
For $k\in\mathbb{Z}$ let
$$f_k(\xi,\mu,\tau')=\chi_k(\xi)(\tau'-\omega(\xi,\mu)+i)^{-1}\cdot \mathcal{F}(u)(\xi,\mu,\tau').$$
For $f_k\in V_k\cap W_k$ let
\begin{equation}\label{ar202}
T(f_k)(\xi,\mu,\tau)=\int_\mathbb{R}f_k(\xi,\mu,\tau')\frac{\psi(\tau-\tau')-
\psi(\tau-\omega(\xi,\mu))}{\tau'-\omega(\xi,\mu)}\cdot (\tau'-\omega(\xi,\mu)+i)\,d\tau'.
\end{equation}
In view of the definitions, it suffices to prove that
\begin{equation}\label{ni5}
||T||_{V_k\cap W_k\to V_k\cap W_k}\leq C\text{ uniformly in }k\in\mathbb{Z}.
\end{equation}

We prove first that
\begin{equation}\label{ni7}
||T(f_k)||_{X_k}\leq C\|f_k\|_{X_k}\text{ uniformly in }k\in \Z.
\end{equation}
We observe the elementary bound
\begin{equation*}
\Big|\frac{\psi(\theta-\theta')-
\psi(\theta)}{\theta'}(\theta'+i)\Big|\leq C[(1+|\theta|)^{-4}+(1+|\theta-\theta'|)^{-4}],
\end{equation*}
for any $\theta,\theta'\in\R$. Thus, for \eqref{ni7}, it suffices to prove that
\begin{equation}\label{ni7.3}
\Big|\Big|\int_{\R}|f_k(\xi,\mu,\tau')|\cdot [(1+|\tau-\tau'| )^{-4}+(1+|\tau-\omega(\xi,\mu)| )^{-4}]\,d\tau'\Big|\Big|_{X_k}\leq C\|f_k\|_{X_k},
\end{equation}
for any $f\in X_k$. For this, we notice first that
\begin{equation*}
\Big|\Big|(1+|\tau-\omega(\xi,\mu)| )^{-4}\int_{\R}|f_k(\xi,\mu,\tau')|\,d\tau'\Big|\Big|_{X_k}\leq C\|f_k\|_{X_k},
\end{equation*}
using \eqref{sd4.1}. In addition, for any $j\geq 0$,
\begin{equation*}
\begin{split}
\Big|\Big|\eta_j(\tau-\omega(\xi,\mu))&\cdot \int_{\R}|f_k(\xi,\mu,\tau')|\cdot (1+|\tau-\tau'| )^{-4}\,d\tau'\Big|\Big|_{L^2}\\
&\leq C\sum_{j'\in\Z}2^{-3|j-j'|}||\eta_{j'}(\tau'-\omega(\xi,\mu))\cdot f_k(\xi,\mu,\tau')||_{L^2}.
\end{split}
\end{equation*}
The bound  \eqref{ni7.3} follows from the definition \eqref{sp7}.

We show now that
\begin{equation}\label{ni30}
||T(f_k)||_{X_k+Y_k}\leq C\|f_k\|_{Y_k}\text{ uniformly in }k\in \Z.
\end{equation}
We may assume $k\geq 100$. Using \eqref{ni7} and Lemma \ref{Lemmaa1} (b), (c), we may also assume that $f_k\in Y_k$ is supported in the set $\{(\xi,\mu,\tau'):|\tau'-\omega(\xi,\mu)|\leq 2^{k-10}\}$. We write
\begin{equation*}
f_k(\xi,\mu,\tau')=\frac{\tau'-\omega(\xi,\mu)}{\tau'-\omega(\xi,\mu)+i}f_k(\xi,\mu,\tau')+\frac{i}{\tau'-\omega(\xi,\mu)+i}f_k(\xi,\mu,\tau').
\end{equation*}
Using Lemma \ref{Lemmaa1} (b), $||i(\tau'-\omega(\xi,\mu)+i)^{-1}f_k(\xi,\mu, \tau')||_{X_k}\leq C||f_k||_{Y_k}$. In view of \eqref{ni7}, it suffices to prove that
\begin{equation}\label{ni8}
\begin{split}
\Big|\Big|\int_\mathbb{R}f_k(\xi,\mu,\tau')&\psi(\tau-\tau')\,d\tau'\Big|\Big|_{X_k+Y_k}\\
&+\Big|\Big|\psi(\tau-\omega(\xi))\int_\mathbb{R}f_k(\xi,\mu,\tau')\,d\tau'\Big|\Big|_{X_k}\leq C||f_k||_{Y_k}.
\end{split}
\end{equation}
The bound for the second term in the left-hand side of \eqref{ni8} follows from
\eqref{sd4.1}. To bound the first term we write
\begin{equation*}
f_k(\xi,\mu,\tau')=f_k(\xi,\mu,\tau')\Big[\frac{\tau'-\omega(\xi,\mu)+i}{\tau-\omega(\xi,\mu)+i}+\frac{\tau-\tau'}{\tau-\omega(\xi,\mu)+i}\Big].
\end{equation*}
The first term in the left-hand side of \eqref{ni8} is dominated by
\begin{equation}\label{ni9}
\begin{split}
&C\Big|\Big|\eta_{\leq k-5}(\tau-\omega(\xi,\mu))\int_\mathbb{R}f_k(\xi,\mu,\tau')\frac{\tau'-\omega(\xi,\mu)+i}{\tau-\omega(\xi,\mu)+i}\cdot \psi(\tau-\tau')\,d\tau'\Big|\Big|_{Y_k}\\
&+C\Big|\Big|\eta_{\leq k-5}(\tau-\omega(\xi,\mu))\int_\mathbb{R}f_k(\xi,\mu,\tau')\frac{\psi(\tau-\tau')\cdot (\tau-\tau')}{\tau-\omega(\xi,\mu)+i}\,d\tau'\Big|\Big|_{X_k}\\
&+C\Big|\Big|\eta_{\geq k-5}(\tau-\omega(\xi,\mu))\int_\mathbb{R}f_k(\xi,\mu,\tau')\psi(\tau-\tau')\,d\tau'\Big|\Big|_{X_k}.
\end{split}
\end{equation}
For the first term in \eqref{ni9}, we use Lemma \ref{Lemmaa1} (c) to bound it by
\begin{equation*}
C2^{-k/2}||\mathcal{F}^{-1}(\psi)\cdot\mathcal{F}^{-1}[(\tau'-\omega(\xi)+i)f_k(\xi,\mu,\tau')]||_{L^1_yL^2_{x,t}}\leq C||f_k||_{Y_k},
\end{equation*}
as desired. To bound the second term, we observe that
\begin{equation*}
\Big|\frac{\psi(\tau-\tau')\cdot (\tau-\tau')}{\tau-\omega(\xi,\mu)+i}\Big|\leq C\frac{(1+|\tau-\tau'|)^{-4}}{1+|\tau'-\omega(\xi,\mu)|}.
\end{equation*} 
Thus the second term in \eqref{ni9} is bounded by
\begin{equation}\label{ni41}
\begin{split}
C&\Big|\Big|\int_\mathbb{R}\frac{|f_k(\xi,\mu,\tau')|}{1+|\tau'-\omega(\xi,\mu)|}\cdot (1+|\tau-\tau'|)^{-4}\,d\tau'\Big|\Big|_{X_k}\leq C\Big|\Big|\frac{f_k(\xi,\mu,\tau')}{1+|\tau'-\omega(\xi,\mu)|}\Big|\Big|_{X_k},
\end{split}
\end{equation}
which is dominated by $C||f_k||_{Y_k}$ in view of Lemma \ref{Lemmaa1} (b). To bound the third term in \eqref{ni9}, recall that $f_k$ is supported in the set $\{(\xi,\mu,\tau'):|\tau'-\omega(\xi,\mu)|\leq 2^{k-10}\}$, thus $f_k(\xi,\mu,\tau')=f_k(\xi,\mu,\tau')\cdot \eta_{\leq k-10}(\tau'-\omega(\xi,\mu))$. In addition, it is easy to see that
\begin{equation*}
|\eta_{\geq k-5}(\tau-\omega(\xi,\mu))\cdot \eta_{\leq k-10}(\tau'-\omega(\xi,\mu))\cdot \psi(\tau-\tau')|\leq C\frac{(1+|\tau-\tau'|)^{-4}}{1+|\tau'-\omega(\xi,\mu)|},
\end{equation*}
so the third term in \eqref{ni9} is also bounded as in \eqref{ni41}.

Finally, we prove that
\begin{equation}\label{ni50}
||T(f_k)||_{W_k}\leq C\|f_k\|_{V_k\cap W_k}\text{ uniformly in }k\in \Z.
\end{equation}
In view of the definition \eqref{sp4}, the left-hand side of \eqref{ni50} is dominated by
\begin{equation*}
\begin{split}
&C\|T[(\partial_\mu+I)f_k]\|_{X_k+Y_k}+C\Big|\Big|\psi'(\tau-\omega(\xi,\mu))\cdot (\mu/ \xi)\int_{\R}f_k(\xi,\mu,\tau')\,d\tau'\Big|\Big|_{X_k+Y_k}\\
&+C\Big|\Big|\int_{\R}f_k(\xi,\mu,\tau')\frac{d}{d\mu}\frac{\psi(\tau-\tau')-
\psi(\tau-\omega(\xi,\mu))}{\tau'-\omega(\xi,\mu)}\,d\tau'\Big|\Big|_{X_k+Y_k}.
\end{split}
\end{equation*}
The first term in the expression above is dominated by $C\|f_k\|_{W_k}$, in view of \eqref{ni7} and \eqref{ni30}. The second term is dominated by $C\|f_k\|_{V_k}$, in view of \eqref{sd4}. Thus, for \eqref{ni50}, it suffices to prove that
\begin{equation}\label{ni51}
\begin{split}
\Big|\Big|\int_{\R}f_k(\xi,\mu,\tau')\cdot \frac{d}{d\mu}\frac{\psi(\tau-\tau')-
\psi(\tau-\omega(\xi,\mu))}{\tau'-\omega(\xi,\mu)}\,d\tau'\Big|\Big|_{X_k}\leq C\|f_k\|_{V_k}.
\end{split}
\end{equation}
By analyzing the cases $|\tau'-\omega(\xi,\mu)|\leq 1$ and $|\tau'-\omega(\xi,\mu)|\geq 1$, it is easy to see that
\begin{equation*}
\begin{split}
\Big|\frac{d}{d\mu}&\frac{\psi(\tau-\tau')-
\psi(\tau-\omega(\xi,\mu))}{\tau'-\omega(\xi,\mu)}\Big|\\
&\leq \frac{C|\mu/ \xi|}{1+|\tau'-\omega(\xi,\mu)|}\cdot[(1+|\tau-\tau'| )^{-4}+(1+|\tau-\omega(\xi,\mu)| )^{-4}].
\end{split}
\end{equation*}
In addition, using Lemma \ref{Lemmaa1} (b), $\|f_k\cdot (1+|\tau-\omega(\xi,\mu)| )^{-1}\cdot |\mu/ \xi|\|_{X_k}\leq C\|f_k\|_{V_k}$. Thus, for \eqref{ni51}, it suffices to prove that
\begin{equation*}
\Big|\Big|\int_{\R}|f_k(\xi,\mu,\tau')|\cdot [(1+|\tau-\tau'| )^{-4}+(1+|\tau-\omega(\xi,\mu)| )^{-4}]\,d\tau'\Big|\Big|_{X_k}\leq C\|f_k\|_{X_k}.
\end{equation*}
This follows from \eqref{ni7.3}. 

The main bound \eqref{ni5} follows from \eqref{ni7}, \eqref{ni30}, and \eqref{ni50}.
\end{proof}

The proof of Proposition \ref{Lemmab3} above (in particular the bounds \eqref{ni7.3} and \eqref{ni8}) shows that if $\varphi\in\mathcal{S}(\R)$ then
\begin{equation}\label{io1}
\begin{cases}
&\|\mathcal{F}[\varphi(t)\cdot \mathcal{F}^{-1}(f)]\|_{X_k}\leq C\|f\|_{X_k};\\
&\|\mathcal{F}[\varphi(t)\cdot \mathcal{F}^{-1}(f)]\|_{X_k+Y_k}\leq C\|f\|_{X_k+Y_k};\\
&\|\mathcal{F}[\varphi(t)\cdot \mathcal{F}^{-1}(f)]\|_{V_k\cap W_k}\leq C\|f\|_{V_k\cap W_k},
\end{cases}
\end{equation}
for any $k\in \Z$ and $f\in L^2(\R^3)$ supported in $I_k\times\R\times\R$.

\section{Proof of Theorem \ref{Main1}}\label{proofthm}

In this section we reduce Theorem \ref{Main1} to proving the following two dyadic bilinear estimates: assume $k_i\in \Z$, $f_{k_i}\in V_{k_i}\cap W_{k_i}$,  and $\mathcal{F}(f_{k_i})$ are supported in $\R^2\times [-2,2]$, $i=1,2$.

$\bullet$ If $k\in \Z$, $k_1\leq k-20$ and $|k_2-k|\leq 2$ then
\begin{equation}\label{BIG1}
\begin{split}
2^k\big|\big|\chi_k(\xi)&\cdot(\tau-\omega(\xi,\mu)+i)^{-1}\cdot (f_{k_1}\ast
f_{k_2})\big|\big|_{V_k\cap W_k}\\
&\leq C(2^{-|k_1|/8}+2^{-|k-k_1|/8})||f_{k_1}||_{V_{k_1}\cap W_{k_1}}||f_{k_2}||_{V_{k_2}\cap W_{k_2}}.
\end{split}
\end{equation}

$\bullet$ If $|k_1-k_2|\leq 100$ then
\begin{equation}\label{BIG2}
\begin{split}
\Big[\sum_{k\in\Z}\big|\big|2^k\chi_k(\xi)\cdot&(\tau-\omega(\xi,\mu)+i)^{-1}\cdot (f_{k_1}\ast
f_{k_2})\big|\big|_{V_k\cap W_k}^2\Big]^{1/2}\\
&\leq C||f_{k_1}||_{V_{k_1}\cap W_{k_1}}\cdot ||f_{k_2}||_{V_{k_2}\cap W_{k_2}}.
\end{split}
\end{equation}

\newtheorem{Lemmao1}{Proposition}[section]
\begin{Lemmao1}\label{Lemmao1}
If $u,v\in F$ then
\begin{equation*}
||\partial_x(uv)||_{N}\leq C||u||_{F}\cdot \|v\|_{F}.
\end{equation*}
\end{Lemmao1}

\begin{proof}[Proof of Proposition \ref{Lemmao1}] Let $f_k=\chi_k(\xi)\cdot \mathcal{F}(u)$ and $g_k=\chi_k(\xi)\cdot \mathcal{F}(v)$. Then
\begin{equation}\label{su1}
\begin{cases}
&\|u\|_F=\big[\sum_{k_1\in\Z}\|f_{k_1}\|_{V_{k_1}\cap W_{k_1}}^2\big]^{1/2};\\
&\|v\|_F=\big[\sum_{k_2\in\Z}\|g_{k_2}\|_{V_{k_2}\cap W_{k_2}}^2\big]^{1/2}.
\end{cases}
\end{equation}
For $k\in\Z$ let
\begin{equation*}
\begin{cases}
&A=\{(k_1,k_2)\in\Z^2:\,|k_1-k_2|\leq 100\};\\
&A_1(k)=\{(k_1,k_2)\in\Z^2:\,|k_2-k|\leq 2\text{ and }k_1\leq k-20\};\\
&A_2(k)=\{(k_1,k_2)\in\Z^2:\,|k_1-k|\leq 2\text{ and }k_2\leq k-20\}.
\end{cases}
\end{equation*}
Clearly
\begin{equation*}
\chi_k(\xi)\cdot \mathcal{F}(\partial_x(uv))=C\chi_k(\xi)\xi\sum_{(k_1,k_2)\in A\cup A_1(k)\cup A_2(k)}(f_{k_1}\ast g_{k_2}).
\end{equation*}
Thus
\begin{equation}\label{su2}
\begin{split}
&||\partial_x(uv)||_{N}^2\leq C\sum_{k\in\Z}\Big(\sum_{(k_1,k_2)\in A}\big|\big|2^k\chi_k(\xi)\cdot(\tau-\omega(\xi,\mu)+i)^{-1}\cdot (f_{k_1}\ast
g_{k_2})\big|\big|_{V_k\cap W_k}\Big)^2\\
&+C\sum_{k\in\Z}\Big(\sum_{(k_1,k_2)\in A_1(k)}\big|\big|2^k\chi_k(\xi)\cdot(\tau-\omega(\xi,\mu)+i)^{-1}\cdot (f_{k_1}\ast
g_{k_2})\big|\big|_{V_k\cap W_k}\Big)^2\\
&+C\sum_{k\in\Z}\Big(\sum_{(k_1,k_2)\in A_2(k)}\big|\big|2^k\chi_k(\xi)\cdot(\tau-\omega(\xi,\mu)+i)^{-1}\cdot (f_{k_1}\ast
g_{k_2})\big|\big|_{V_k\cap W_k}\Big)^2.
\end{split}
\end{equation}

Using \eqref{BIG2}, the first  term in the right-hand side of \eqref{su2} is bounded by
\begin{equation*}
C\Big[\sum_{(k_1,k_2)\in A}||f_{k_1}||_{V_{k_1}\cap W_{k_1}}\cdot ||g_{k_2}||_{V_{k_2}\cap W_{k_2}}\Big]^2\leq C\|u\|_F^2\cdot ||v||_{F}^2.
\end{equation*}
Using \eqref{BIG1}, the second term in the right-hand side of \eqref{su2} is bounded by
\begin{equation*}
C\sum_{k\in\Z}\Big(||u||_F\cdot \sum_{|k_2-k|\leq 2}||g_{k_2}||_{V_{k_2}\cap W_{k_2}}\Big)^2\leq C\|u\|_F^2\cdot ||v||_{F}^2.
\end{equation*}
The third term in the right-hand side of \eqref{su2} is similar, and the proposition follows.
\end{proof}

If follows from Proposition  \ref{Lemmab3} and Proposition \ref{Lemmao1} that
\begin{equation}\label{su4}
\Big|\Big|\eta_0(t)\cdot \int_0^t[W(\partial_x (uv))(s)](t-s)\,ds\Big|\Big|_{F}\leq C||u||_{F}\cdot ||v||_{F},
\end{equation}
for any $u,v\in F$. It is easy to show that $F$ is a  Banach space, and Theorem \ref{Main1} follows from \eqref{su4} and Proposition \ref{Lemmab1} by a standard fixed-point argument.

The rest of the paper is concerned with the proofs of the dyadic bilinear estimates \eqref{BIG1} and \eqref{BIG2}.

\section{Preliminary estimates}\label{prop}

In this section we prove several localized $L^\infty_yL^2_{x,t}$ and
$L^2_yL^\infty_{x,t}$ estimates and an $L^4$ Strichartz
estimate. These bounds will be used in the bilinear estimates in
Sections  \ref{bilinear2}, \ref{bilinear1}, and \ref{bilinear3}. We start with a representation formula for functions in $Y_k$, $k\geq 100$. Let $\mathbf{1}_+$ and $\mathbf{1}_-$ denote the characteristic functions of the intervals $[0,\infty)$ and $(-\infty,0]$ respectively. 

\newtheorem{Lemmav0}{Lemma}[section]
\begin{Lemmav0}\label{Lemmav0}
If $k\geq 100$ and $f\in Y_k$, then $f$ can be written in the form
\begin{equation}\label{io10}
\begin{split}
f(\xi,\mu,\tau)&=2^{-k/2}\mathbf{1}_{I_k}(\xi)\cdot \chi_{[2k-30,2k+30]}(M)\cdot \mathbf{1}_+(M)\\
&\times \Big(\frac{\eta_0(M-\mu)}{M-\mu+i/2^k}+\frac{\eta_0(M+\mu)}{M+\mu+i/2^k}\Big)\cdot\int_{\R}e^{-iy\cdot \mu}g(y,\xi,\tau)\,dy+h,
\end{split}
\end{equation}
where $M=M(\xi,\tau)=\sqrt{\xi\cdot (\tau-\xi^3)}$, $h$ is supported in the set $\{(\xi,\mu,\tau): \xi\in I_k,\,|\mu|\in[2^{2k-100},2^{2k+100}]\}$, and
\begin{equation}\label{io11}
\|h\|_{X_k}+\|g\|_{L^1_yL^2_{\xi,\tau}}\leq C\|f\|_{Y_k}.
\end{equation}
\end{Lemmav0}

\begin{proof}[Proof of Lemma \ref{Lemmav0}] We start from the identity \eqref{me1}. Since $|\xi|\in[2^{k-2},2^{k+2}]$, $|\mu|\in [2^{2k-35},2^{2k+35}]$, $|\tau-\xi^3-\mu^2/ \xi|\leq 2^{k+2}$, we have $\xi\cdot(\tau-\xi^3)\in [2^{4k-80},2^{4k+80}]$. So $M=M(\xi,\tau)=\sqrt{\xi\cdot (\tau-\xi^3)}$ is well-defined and $M\in[2^{2k-40},2^{2k+40}]$. For $\xi\in I_k$, an elementary computation shows that we can approximate
\begin{equation*}
\begin{split}
&\chi_{[2k-30,2k+30]}(\mu)\cdot \mathbf{1}_+(\mu)\cdot \frac{\eta_{\leq k+1}(\tau-\omega(\xi,\mu))}{\tau-\omega(\xi,\mu)+i}\\
&=\chi_{[2k-30,2k+30]}(M)\cdot \mathbf{1}_+(M)\cdot \frac{\xi}{2M}\cdot \frac{\eta_0(M-\mu)}{M-\mu+i/2^k}+E_+(\xi,\mu,\tau)
\end{split}
\end{equation*}
where, with $\beta=|\tau-\omega(\xi,\mu)|+1$,
\begin{equation}\label{io17}
|E_+(\xi,\mu,\tau)|\leq C\chi_{[2k-40,2k+40]}(\mu)\cdot\frac{\eta_{\leq k+100}(\beta)}{\beta }\cdot \Big(\frac{\beta}{2^k}+\frac{1}{\beta}\Big). 
\end{equation}
Similarly, we approximate
\begin{equation*}
\begin{split}
&\chi_{[2k-30,2k+30]}(\mu)\cdot \mathbf{1}_-(\mu)\cdot \frac{\eta_{\leq k+1}(\tau-\omega(\xi,\mu))}{\tau-\omega(\xi,\mu)+i}\\
&=\chi_{[2k-30,2k+30]}(M)\cdot \mathbf{1}_+(M)\cdot \frac{\xi}{2M}\cdot \frac{\eta_0(M+\mu)}{M+\mu+i/2^k}+E_-(\xi,\mu,\tau),
\end{split}
\end{equation*}
with $E_-$ satisfying the same bound \eqref{io17}. We substitute these formulas into \eqref{me1} and notice that the terms corresponding to $E_+$ and $E_-$ can be estimated in $X_k$ (as in the proof of Lemma \ref{Lemmaa1} (b)). The bound \eqref{io11} follows from \eqref{me2}.
\end{proof}

We prove now a localized $L^\infty_yL^2_{x,t}$ estimate.

\newtheorem{Lemmav1}[Lemmav0]{Lemma}
\begin{Lemmav1}\label{Lemmav1}
Assume $k\geq 0$, $l\geq 2k-100$, and $f$ is supported in the set
\begin{equation*}
\{(\xi,\mu,\tau)\in\mathbb{R}^3:\xi\in I_k,\,|\mu|\in[2^{l-1},2^{l+1}]\}. 
\end{equation*}

(a) Then
\begin{equation}\label{am1}
||\mathcal{F}^{-1}(f)||_{L^\infty_yL^2_{x,t}}\leq C2^{-(l-k)/2}||f||_{X_k+Y_k}.
\end{equation}

(b) More generally, if $\varphi:\mathbb{R}\to[0,1]$ is a smooth function supported in the interval $[-2,2]$, $\epsilon\geq 2^{-k}$, and
\begin{equation*}
f^m_{\pm}(\xi,\mu,\tau)=f(\xi,\mu,\tau)\cdot \varphi((\mu/ \xi\pm \sqrt{3}\xi)/ \epsilon-m)\text{ for }m\in\mathbb{Z},
\end{equation*}
then
\begin{equation}\label{am100}
\Big[\sum_{m\in\mathbb{Z}}||\mathcal{F}^{-1}(f^m_{\pm})||_{L^\infty_yL^2_{x,t}}^2\Big]^{1/2}\leq C2^{-(l-k)/2}||f||_{X_k+Y_k}.
\end{equation}
\end{Lemmav1}

\begin{proof}[Proof of Lemma \ref{Lemmav1}] For part (a), assume first $f\in X_k$. Then (see \cite[p. 753]{CoKeSt})
\begin{equation}\label{am2}
||\mathcal{F}^{-1}(f)||_{L^\infty_yL^2_{x,t}}\leq C2^{-(l-k)/2}||f||_{X_k},
\end{equation}
as desired. Assume now that $f\in Y_k$, $k\geq 100$. We use the representation \eqref{io10} and the bound \eqref{io11}. In view of \eqref{am2}, and using Plancherel's theorem, it suffices to prove that
\begin{equation}\label{pr41}
\Big |\int_{\mathbb{R}}e^{iy_0\cdot \mu }\cdot \frac{\eta_0(M\pm \mu)}{M\pm \mu+i/2^k}\,d\mu \Big|\leq C,
\end{equation}
uniformly in $y_0$, $M\in [2^{2k-40},2^{2k+40}]$, and $\xi\in I_k$. This is a standard uniform estimate for the inverse Fourier transform of a Calder\'{o}n--Zygmund kernel.

For part (b), if $f\in X_k$, then \eqref{am100} follows from \eqref{am2} by orthogonality. Assume now $f\in Y_k$, $k\geq 100$. We use \eqref{io10}, so we may assume
\begin{equation*}
\begin{split}
f^m_{\pm}(\xi,\mu,\tau)&=2^{-k/2}\mathbf{1}_{I_k}(\xi)\cdot \chi_{[2k-30,2k+30]}(M)\cdot \mathbf{1}_+(M)\cdot \varphi((\mu/ \xi\pm \sqrt{3}\xi)/ \epsilon-m) \\
&\times \Big(\frac{\eta_0(M-\mu)}{M-\mu+i/2^k}+\frac{\eta_0(M+\mu)}{M+\mu+i/2^k}\Big)\cdot\int_{R}e^{-iy\cdot \mu}g(y,\xi,\tau)\,dy.
\end{split}
\end{equation*}
By comparing the supports in $\mu$ of the functions and using the fact that $2^k\epsilon\geq 1$, we conclude that $f^m_{\pm}(\xi,\mu,\tau)\equiv 0$ unless $(\tau-\xi^3)/ \xi\in[C^{-1}2^{2k},C2^{2k}]$ and
\begin{equation*}
\Big|\frac{\sqrt{(\tau-\xi^3)/ \xi}\pm\sqrt{3}\xi}{\epsilon}-m\Big|\leq C_0\text{ or }\Big|\frac{-\sqrt{(\tau-\xi^3)/ \xi}\pm\sqrt{3}\xi}{\epsilon}-m\Big|\leq C_0.
\end{equation*}
We define
\begin{equation*}
\begin{split}
g^m_{\pm}&(y,\xi,\tau)=g(y,\xi,\tau)\\
&\times\Big[\eta_0\Big(\frac{\sqrt{(\tau-\xi^3)/ \xi}\pm\sqrt{3}\xi}{C_0\epsilon}-\frac{m}{C_0}\Big)+\eta_0\Big(\frac{-\sqrt{(\tau-\xi^3)/ \xi}\pm\sqrt{3}\xi}{C_0\epsilon}-\frac{m}{C_0}\Big)\Big].
\end{split}
\end{equation*}
In view of the support property above, we have
\begin{equation*}
\begin{split}
f^m_{\pm}(\xi,\mu,\tau)&=2^{-k/2}\mathbf{1}_{I_k}(\xi)\cdot \chi_{[2k-30,2k+30]}(M)\cdot \mathbf{1}_+(M)\cdot \varphi((\mu/ \xi\pm \sqrt{3}\xi)/ \epsilon-m) \\
&\times \Big(\frac{\eta_0(M-\mu)}{M-\mu+i/2^k}+\frac{\eta_0(M+\mu)}{M+\mu+i/2^k}\Big)\cdot\int_{R}e^{-iy\cdot \mu}g^m_{\pm}(y,\xi,\tau)\,dy.
\end{split}
\end{equation*}
Using part (a) (in fact a slightly modified version of the bound \eqref{pr41}),
\begin{equation*}
||\mathcal{F}^{-1}(f^m_{\pm})||_{L^\infty_yL^2_{x,t}}\leq C2^{-k/2}||g^m_{\pm}||_{L^1_yL^2_{\xi,\tau}}.
\end{equation*}
Thus, the left-hand side of \eqref{am100} is dominated by
\begin{equation*}
C2^{-k/2}\Big[\sum_{m\in\mathbb{Z}}||g^m_{\pm}||_{L^1_yL^2_{\xi,\tau}}^2\Big]^{1/2}\leq C2^{-k/2}||g||_{L^1_yL^2_{\xi,\tau}},
\end{equation*}
which suffices in view of \eqref{io11}.
\end{proof}

We prove now several localized maximal function estimates: 

\newtheorem{Lemmav2}[Lemmav0]{Lemma}
\begin{Lemmav2}\label{Lemmav2}
Assume $k,l,j\in\mathbb{Z}$, $k\leq 0$, $l\geq 0$, $j\geq 0$.

(a) If $f$ is supported in the set
\begin{equation*}
\{(\xi,\mu,\tau)\in\mathbb{R}^3:\xi\in I_k,\,|\mu|\leq 2^l,\,|\tau-\omega(\xi,\mu)|\leq 2^j\},
\end{equation*}
then
\begin{equation}\label{max1}
||\mathcal{F}^{-1}(f)||_{L^2_yL^\infty_{x,t}}\leq C2^{j/2}\cdot 2^{(2l+k)/4}||(I-\partial_\tau^2)f||_{L^2}.
\end{equation}

(b) If $m\in\mathbb{R}$, $\epsilon\geq 2^{-l}$, and $f$ is supported in the set
\begin{equation*}
\{(\xi,\mu,\tau)\in\mathbb{R}^3:\xi\in I_k,\,|\mu|\leq 2^l,\,|\tau-\omega(\xi,\mu)|\leq 2^j,\,|\mu/
\xi\pm\sqrt{3}\xi-m|\leq \epsilon \},
\end{equation*}
then
\begin{equation}\label{max3}
||\mathcal{F}^{-1}(f)||_{L^2_yL^\infty_{x,t}}\leq C2^{j/2}\cdot (2^l\epsilon)^{1/2}\cdot 2^{k/2}||(I-\partial_\tau^2)f||_{L^2}.
\end{equation}
\end{Lemmav2}

\begin{proof}[Proof of Lemma \ref{Lemmav2}] For any $f:\mathbb{R}^3\to\mathbb{C}$ let $f^\#(\xi,\mu,\theta)=f(\xi,\mu,\theta+\omega(\xi,\mu))$. Then
\begin{equation}\label{timedecay}
\mathcal{F}^{-1}(f)(x,y,t)=C(t^2+1)^{-1}\int_{\mathbb{R}^3}[(I-\partial_\tau^2)f]^\#(\xi,\mu,\theta)e^{it\theta}e^{i(x\cdot \xi+y\cdot\mu+t\cdot\omega(\xi,\mu))}\,d\xi d\mu d\theta.
\end{equation}
Thus, for \eqref{max1}, after noticing the time decay in \eqref{timedecay}, it suffices to prove that if
\begin{equation*}
g\text{ is supported in the set }\{(\xi,\mu):\xi\in I_k,\,|\mu|\leq 2^l\},
\end{equation*}
then
\begin{equation}\label{v22}
\Big|\Big|\int_{\mathbb{R}^2}g(\xi,\mu)e^{i(x\cdot \xi+y\cdot\mu+t\cdot\omega(\xi,\mu))}\,d\xi d\mu\Big|\Big|_{L^2_yL^\infty_{x,|t|\leq 1/2}}\leq C2^{(2l+k)/4}||g||_{L^2}.
\end{equation}
A standard $TT^{\ast}$ argument (see, for example, \cite[p. 50]{KeZi}), shows that for \eqref{v22} it suffices to prove that
\begin{equation}\label{v23}
\Big|\Big|\int_{\mathbb{R}^2}\chi_{[k-1,k+1]}^2(\xi)\eta_0^2(\mu/2^l)e^{i(x\cdot \xi+y\cdot\mu+t\cdot\omega(\xi,\mu))}\,d\xi d\mu\Big|\Big|_{L^1_yL^\infty_{x,|t|\leq 1}}\leq C2^{(2l+k)/2}.
\end{equation}
To prove  \eqref{v23} we estimate the $\mu$-integral first. Simple integration by parts and van der Corput-type arguments show that if $y\in\mathbb{R}$, $|t|\leq 1$, $|\xi|\in [2^{k-2},2^{k+2}]$, and $k,l$ are as in the hypothesis then
\begin{equation*}
\Big|\int_{\mathbb{R}}\eta_0^2(\mu/2^l)e^{i(y\cdot\mu+t\cdot \mu^2/ \xi)}\,d\mu\Big|\leq C
\begin{cases}
2^{l-k}|y|^{-2}&\text{ if }|y|\geq 100\cdot 2^{l-k};\\
2^{l/2}|y|^{-1/2}&\text{ if }|y|\in[1,100\cdot 2^{l-k}];\\
2^l&\text{ if }|y|\leq 1.
\end{cases}
\end{equation*}
This leads to \eqref{v23}.

Similarly, for \eqref{max3}, it suffices to prove that if
\begin{equation*}
g\text{ is supported in the set }\{(\xi,\mu):\xi\in I_k,\,|\mu|\leq 2^l,\,|\mu/ \xi\pm\sqrt{3}\xi-m|\leq \epsilon\},
\end{equation*}
then
\begin{equation}\label{v26}
\Big|\Big|\int_{\mathbb{R}^2}g(\xi,\mu)e^{i(x\cdot \xi+y\cdot\mu+t\cdot\omega(\xi,\mu))}\,d\xi d\mu\Big|\Big|_{L^2_yL^\infty_{x,|t|\leq 1/2}}\leq C(2^{l}\epsilon)^{1/2}\cdot 2^{k/2}||g||_{L^2}.
\end{equation}
In proving \eqref{v26}, by orthogonality, we may assume $\epsilon=2^{-l}$. We may also assume $|m|\leq C2^{l-k}$. As before, for \eqref{v26}, it suffices to prove that
\begin{equation}\label{v25}
\Big|\Big|\int_{\mathbb{R}^2}\chi_{[k-1,k+1]}^2(\xi)\eta_0^2(2^l(\mu/ \xi \pm\sqrt{3}\xi-m))e^{i(x\cdot \xi+y\cdot\mu+t\cdot\omega(\xi,\mu))}\,d\xi d\mu\Big|\Big|_{L^1_yL^\infty_{x,|t|\leq 1}}\leq C2^{k}.
\end{equation}
The change of variables $\mu=\xi(\mp\sqrt{3}\xi+m+2^{-l}\beta)$, with $d\mu=2^{-l}\xi d\beta$, and integration by parts show that
\begin{equation*}
\Big|\int_{\mathbb{R}}\eta_0^2(2^l(\mu/ \xi \pm\sqrt{3}\xi-m))e^{i(y\cdot\mu+t\cdot\mu^2/ \xi)}d\mu\Big|\leq C2^{k-l}(1+2^{k-l}|y| )^{-2},
\end{equation*}
if $y\in\mathbb{R}$, $|m|\leq C2^{l-k}$, $|\xi|\in[2^{k-2},2^{k+2}]$, and $|t|\leq 1$. This leads to \eqref{v25}.
\end{proof}

\newtheorem{Lemmav3}[Lemmav0]{Lemma}
\begin{Lemmav3}\label{Lemmav3}
Assume $k,l,j\in\mathbb{Z}_+$.

(a) If $f$ is supported in the set
\begin{equation*}
\{(\xi,\mu,\tau)\in\mathbb{R}^3:\xi\in I_k,\,|\mu|\leq 2^l,\,|\tau-\omega(\xi,\mu)|\leq 2^j\},
\end{equation*}
then, for any $\delta>0$,
\begin{equation}\label{max1.1}
||\mathcal{F}^{-1}(f)||_{L^2_yL^\infty_{x,t}}\leq C_\delta 2^{j/2}\cdot (2^k+2^{l-k})^{1/2+\delta }||(I-\partial_\tau^2)f||_{L^2}.
\end{equation}

(b) If $m\in\mathbb{R}$, $l\geq 2k$, $\epsilon\geq 2^{-l}$, and $f$ is supported in the set
\begin{equation*}
\{(\xi,\mu,\tau)\in\mathbb{R}^3:\xi\in I_k,\,|\mu|\leq 2^l,\,|\tau-\omega(\xi,\mu)|\leq 2^j,\,|\mu/
\xi\pm\sqrt{3}\xi-m|\leq \epsilon \},
\end{equation*}
then
\begin{equation}\label{max3.1}
||\mathcal{F}^{-1}(f)||_{L^2_yL^\infty_{x,t}}\leq C2^{j/2}\cdot (2^l\epsilon)^{1/2}||(I-\partial_\tau^2)f||_{L^2}.
\end{equation}
\end{Lemmav3}

\begin{proof}[Proof of Lemma \ref{Lemmav2}] As in the proof of Lemma \ref{Lemmav2}, for \eqref{max1.1} it suffices to show that if
\begin{equation*}
g\text{ is supported in the set }\{(\xi,\mu):\xi\in I_k,\,|\mu|\leq 2^l\},
\end{equation*}
then
\begin{equation*}
\Big|\Big|\int_{\mathbb{R}^2}g(\xi,\mu)e^{i(x\cdot \xi+y\cdot\mu+t\cdot\omega(\xi,\mu))}\,d\xi d\mu\Big|\Big|_{L^2_yL^\infty_{x,|t|\leq 1/2}}\leq C_\delta  (2^k+2^{l-k})^{1/2+\delta }||g||_{L^2}.
\end{equation*}
This follows from \cite[Theorem 2.1 (b)]{KeZi}.

Similarly, for \eqref{max3.1}, it suffices to prove that if
\begin{equation*}
g\text{ is supported in the set }\{(\xi,\mu):\xi\in I_k,\,|\mu|\leq 2^l,\,|\mu/ \xi\pm\sqrt{3}\xi-m|\leq \epsilon\},
\end{equation*}
then
\begin{equation}\label{v26.1}
\Big|\Big|\int_{\mathbb{R}^2}g(\xi,\mu)e^{i(x\cdot \xi+y\cdot\mu+t\cdot\omega(\xi,\mu))}\,d\xi d\mu\Big|\Big|_{L^2_yL^\infty_{x,|t|\leq 1/2}}\leq C(2^{l}\epsilon)^{1/2}||g||_{L^2}.
\end{equation}
In proving \eqref{v26}, by orthogonality, we may assume $\epsilon=2^{-l}$. We may also assume $|m|\leq 2^{l-k+3}$. As before, for \eqref{v26.1}, it suffices to prove that
\begin{equation}\label{v25.1}
\Big|\Big|\int_{\mathbb{R}^2}\chi_{[k-1,k+1]}^2(\xi)\eta_0^2(2^l(\mu/ \xi \pm\sqrt{3}\xi-m))e^{i(x\cdot \xi+y\cdot\mu+t\cdot\omega(\xi,\mu))}\,d\xi d\mu\Big|\Big|_{L^1_yL^\infty_{x,|t|\leq 1}}\leq C.
\end{equation}
We make the change of variables $\mu=\xi(\mp\sqrt{3}\xi+m+2^{-l}\beta)$, with $d\mu=2^{-l}\xi d\beta$. The estimate \eqref{v25.1} becomes
\begin{equation}\label{v27.1}
2^{-l}\Big|\Big|\int_{\mathbb{R}^2}\xi\cdot \chi_{[k-1,k+1]}^2(\xi)\eta_0^2(\beta)e^{i\Phi(x,y,t,\xi,\beta)}\,d\xi d\beta\Big|\Big|_{L^1_yL^\infty_{x,|t|\leq 1}}\leq C,
\end{equation}
where
\begin{equation}\label{v28.1}
\Phi(x,y,t,\xi,\beta)=x\cdot \xi+y\cdot \xi(\mp\sqrt{3}\xi+m+2^{-l}\beta)+t\cdot \xi^3+t\cdot \xi(\mp\sqrt{3}\xi+m+2^{-l}\beta)^2.
\end{equation}

It remains to prove \eqref{v27.1}. For $|y|\leq 2^{l-k+10}$ we notice that $|\partial_\xi^3\Phi(x,y,t,\xi,\beta)|\geq |t|$ and $|\partial_\xi^2\Phi(x,y,t,\xi,\beta)|\geq 2\sqrt{3}|y|-C2^{l-k}|t|$, provided that $|\xi|\approx  2^k$ and $|m|\leq2^{l-k+3}$. Thus, using van der Corput's lemma for the integral in $\xi$,
\begin{equation}\label{v30.1}
2^{-l}\Big|\int_{\mathbb{R}^2}\xi\cdot \chi_{[k-1,k+1]}^2(\xi)\eta_0^2(\beta)e^{i\Phi(x,y,t,\xi,\beta)}\,d\xi d\beta\Big|\leq C2^{k-l}\cdot 2^{(l-k)/2}|y|^{-1/2}.
\end{equation}
For $|y|\geq 2^{l-k+10}$ we integrate first by parts in $\beta$ (notice that $|\partial_\beta\Phi|\geq 2^{k-l-4}|y|$ and $|\partial^2_\beta\Phi|\leq C2^{k-2l}|$ if $|t|\leq 1$). Then we use van der Corput's lemma for the integral in $\xi$ as before. The result is
\begin{equation}\label{v30.2}
2^{-l}\Big|\int_{\mathbb{R}^2}\xi\cdot \chi_{[k-1,k+1]}^2(\xi)\eta_0^2(\beta)e^{i\Phi(x,y,t,\xi,\beta)}\,d\xi d\beta\Big|\leq C2^{k-l}\cdot(2^{k-l}|y|)^{-1}\cdot |y|^{-1/2}.
\end{equation}
The bound \eqref{v27.1} follows from \eqref{v30.1} and \eqref{v30.2}.
\end{proof}

We conclude this section with an $L^4$ estimate.

\newtheorem{Lemmav5}[Lemmav0]{Lemma}
\begin{Lemmav5}\label{Lemmav5}
If $k\in\Z$ and $f\in X_k+Y_k$ then
\begin{equation}\label{io2}
\|\mathcal{F}^{-1}(f)\|_{L^4_{x,y,t}}\leq C\|f\|_{X_k+Y_k}.
\end{equation}
\end{Lemmav5}

\begin{proof}[Proof of Lemma \ref{Lemmav5}] We use the scale-invariant Strichartz estimate of \cite{ArSa}:
\begin{equation}\label{io30}
\Big|\Big|\int_{\R^2}\phi(\xi,\mu)e^{ix\cdot\xi}e^{iy\cdot\mu}e^{it\cdot\omega(\xi,\mu)}\,d\xi d\mu\Big|\Big|_{L^4_{x,y,t}}\leq C\|\phi\|_{L^2},
\end{equation}
for any $\phi\in L^2(\R^2)$. 

Assume first that $f\in X_k$. With $f^\#$ defined as in the proof of Lemma \ref{Lemmav2}, for $j\geq 0$
\begin{equation}\label{re1}
\begin{split}
\Big|\Big|&\int_{\R^3}f(\xi,\mu,\tau)\cdot \eta_j(\tau-\omega(\xi,\mu))\cdot e^{ix\cdot\xi}e^{iy\cdot\mu}e^{it\cdot\tau}\,d\xi d\mu d\tau\Big|\Big|_{L^4_{x,y,t}}\\
&=\Big|\Big|\int_{\R^3}f^\#(\xi,\mu,\theta)\cdot \eta_j(\theta)e^{it\cdot\theta}\cdot e^{ix\cdot\xi}e^{iy\cdot\mu}e^{it\cdot\omega(\xi,\mu)}\,d\xi d\mu d\theta\Big|\Big|_{L^4_{x,y,t}}\\
&\leq C2^{j/2}||f^\#(\xi,\mu,\theta)\cdot \eta_j(\theta)||_{L^2},
\end{split}
\end{equation}
which gives \eqref{io2}.

Assume now that $f\in Y_k$. We use the representation \eqref{io10}. With the notation in Lemma \ref{Lemmav0}, using  \eqref{io11} and  the bound \eqref{io2} for $f\in X_k$, it suffices to prove that
\begin{equation*}
\begin{split}
2^{-k/2}&\Big|\Big|\int_{\R^3}g(\xi,\tau)\mathbf{1}_{I_k}(\xi)\cdot \chi_{[2k-30,2k+30]}(M)\cdot \mathbf{1}_+(M)\\
&\times\frac{\eta_0(M\pm\mu)}{M\pm\mu+i/2^k}\cdot e^{ix\cdot\xi}e^{iy\cdot\mu}e^{it\cdot\tau}\,d\xi d\mu d\tau\Big|\Big|_{L^4_{x,y,t}}\leq C||g||_{L^2},
\end{split}
\end{equation*}
for any $g\in L^2(\R^2)$. We take the integral in $\mu$ first; it remains  to prove that
\begin{equation*}
\begin{split}
\Big|\Big|\int_{\R^2}g(\xi,\tau)&\mathbf{1}_{I_k}(\xi)\cdot \chi_{[2k-30,2k+30]}(M)\cdot \mathbf{1}_+(M)\\
&\times e^{ix\cdot\xi}e^{iy\cdot M}e^{it\cdot\tau}\,d\xi d\tau\Big|\Big|_{L^4_{x,y,t}}\leq C2^{k/2}||g||_{L^2}.
\end{split}
\end{equation*}
We make the change of variables $\tau=\xi^3+\nu^2/ \xi$, $\nu\in [C^{-1}2^{2k},C2^{2k}]$, $d\tau =2(\nu/ \xi)d\nu$. Clearly, $M(\xi,\tau)=\nu$. Thus, it suffices to prove that
  \begin{equation*}
\begin{split}
\Big|\Big|\int_{\R^2}&g(\xi,\xi^3+\nu^2/ \xi)\mathbf{1}_{I_k}(\xi)\cdot \chi_{[2k-30,2k+30]}(\nu)\\
&\times \mathbf{1}_+(\nu)e^{ix\cdot\xi}e^{iy\cdot \nu}e^{it\cdot\omega(\xi,\nu)}\,d\xi d\nu\Big|\Big|_{L^4_{x,y,t}}\leq C2^{-k/2}||g||_{L^2}.
\end{split}
\end{equation*}
This follows from \eqref{io30} with $\phi(\xi,\nu)=g(\xi,\xi^3+\nu^2/ \xi)\mathbf{1}_{I_k}(\xi)\cdot \chi_{[2k-30,2k+30]}(\nu)\mathbf{1}_+(\nu)$.
\end{proof}

\section{An $L^2$ bilinear estimate}\label{L2bi}

In this section we prove an $L^2$ bilinear estimate. For $k\in\Z$ and $j\in \Z_+$ let
\begin{equation*}
D_{k,j}=\{(\xi,\mu,\tau):\,\xi\in I_k,\,\mu\in\R,\,|\tau-\omega(\xi,\mu)|\leq 2^j\}.
\end{equation*}

\newtheorem{Main9l}{Lemma}[section]
\begin{Main9l}\label{Main9l}
Assume $k_1,k_2,k_3\in \Z$, $j_1,j_2,j_3\in\Z_+$, and $f_i:\R^3\to\R_+$ are $L^2$ functions supported in $D_{k_i,j_i}$, $i=1,2,3$. If
\begin{equation}\label{jj1.1}
\max(j_1,j_2,j_3)\leq k_1+k_2+k_3-20
\end{equation}
then
\begin{equation}\label{jj1}
\int_{\R^3}(f_1\ast f_2)\cdot f_3\leq C2^{(j_1+j_2+j_3)/2}\cdot 2^{-(k_1+k_2+k_3)/2}\cdot \|f_1\|_{L^2}\|f_2\|_{L^2}\|f_3\|_{L^2}.
\end{equation}
\end{Main9l}
Before we proceed to the proof of this lemma we state a simple corollary that follows by duality.
\newtheorem{Main9}[Main9l]{Corollary}
\begin{Main9}\label{Main9}
Assume $k_1,k_2,k_3\in \Z$, $j_1,j_2,j_3\in\Z_+$, and $f_i:\R^3\to\R_+$  are $L^2$ functions supported in $D_{k_i,j_i}$, $i=1,2$. If 
\begin{equation*}
\max(j_1,j_2,j_3)\leq k_1+k_2+k_3-20
\end{equation*}
then
\begin{equation*}
\|(f_1\ast f_2)\cdot \mathbf{1}_{D_{k_3,j_3}}\|_{L^2}\leq C2^{(j_1+j_2+j_3)/2}\cdot 2^{-(k_1+k_2+k_3)/2}\cdot \|f_1\|_{L^2}\|f_2\|_{L^2}.
\end{equation*}
\end{Main9}
\begin{proof}[Proof of Lemma \ref{Main9l}]
Clearly,
\begin{equation*}
\int_{\R^3}(f_1\ast f_2)\cdot f_3=\int_{\R^3}(\widetilde{f}_1\ast f_3)\cdot f_2=\int_{\R^3}(\widetilde{f}_2\ast f_3)\cdot f_1,
\end{equation*}
where $\widetilde{f}_i(\xi,\mu,\tau)=f_i(-\xi,-\mu,-\tau)$, $i=1,2$. In view of the symmetry of \eqref{jj1} we may assume
\begin{equation}\label{jj0}
j_3=\max(j_1,j_2,j_3).
\end{equation}
As in the proof of Lemma \ref{Lemmav2}, we define $f_i^\#(\xi,\mu,\theta)=f_i(\xi,\mu,\theta+\omega(\xi,\mu))$, $i=1,2,3$, $\|f_i^\#\|_{L^2}=\|f_i\|_{L^2}$. We rewrite the left-hand side of \eqref{jj1} in the form
\begin{equation}\label{jj3}
\begin{split}
&\int_{\R^6}f_1^\#(\xi_1,\mu_1,\theta_1)\cdot f_2^\#(\xi_2,\mu_2,\theta_2)\\
&\times f_3^\#(\xi_1+\xi_2,\mu_1+\mu_2,\theta_1+\theta_2+\Omega((\xi_1,\mu_1),(\xi_2,\mu_2)))\,d\xi_1d\xi_2d\mu_1d\mu_2d\theta_1d\theta_2,
\end{split}
\end{equation}
where
\begin{equation}\label{jj2}
\begin{split}
\Omega((\xi_1,\mu_1),(\xi_2,\mu_2))&=-\omega(\xi_1+\xi_2,\mu_1+\mu_2)+\omega(\xi_1,\mu_1)+\omega(\xi_2,\mu_2)\\
&=\frac{-\xi_1\xi_2}{\xi_1+\xi_2}\Big[(\sqrt{3}\xi_1+\sqrt{3}\xi_2)^2-\Big(\frac{\mu_1}{\xi_1}-\frac{\mu_2}{\xi_2}\Big)^2\Big].
\end{split}
\end{equation}
The functions $f_i^\#$ are supported in the sets $\{\xi,\mu,\theta):\,\xi\in I_{k_i},\,\mu\in\R,\,|\theta|\leq 2^{j_i}\}$. 

We will prove that if $g_i:\R^2\to\R_+$ are $L^2$ functions supported in $I_{k_i}\times\R$, $i=1,2$, and $g:\R^3\to\R_+$ is an $L^2$ function supported in $I_k\times\R\times [-2^j,2^j]$, $j\leq k_1+k_2+k-15$, then
\begin{equation}\label{jj5}
\begin{split}
\int_{\R^4}g_1(\xi_1,\mu_1)\cdot &g_2(\xi_2,\mu_2)\cdot g(\xi_1+\xi_2,\mu_1+\mu_2,\Omega((\xi_1,\mu_1),(\xi_2,\mu_2)))\,d\xi_1d\xi_2d\mu_1d\mu_2\\
&\leq C2^{j/2}\cdot 2^{-(k_1+k_2+k)/2}\cdot\|g_1\|_{L^2}\|g_2\|_{L^2}\|g\|_{L^2}. 
\end{split}
\end{equation}
This suffices for \eqref{jj1}, in view of \eqref{jj0} and \eqref{jj3}.

To prove \eqref{jj5}, we observe\footnote{There are four identical integrals of this type.} first that we may assume that the integral in the left-hand side of \eqref{jj5} is taken over the set
\begin{equation*}
\mathcal{R}_{++}=\{(\xi_1,\mu_1,\xi_2,\mu_2):\,\xi_1+\xi_2\geq 0\text{ and }\mu_1/ \xi_1-\mu_2/ \xi_2\geq 0\}.
\end{equation*}
Using the restriction $j\leq k_1+k_2+k-15$ and \eqref{jj2}, we may assume also that the integral in the left-hand side of \eqref{jj5} is taken over the set
\begin{equation*}
\widetilde{\mathcal{R}}_{++}=\{(\xi_1,\mu_1,\xi_2,\mu_2)\in\mathcal{R}_{++}:|\sqrt{3}(\xi_1+\xi_2)|-|\mu_1/ \xi_1-\mu_2/ \xi_2|\leq 2^{-10}|\xi_1+\xi_2|\}.
\end{equation*}
To summarize, it suffices to prove that
\begin{equation}\label{jj11}
\begin{split}
\int_{\widetilde{\mathcal{R}}_{++}}g_1(\xi_1,\mu_1)\cdot &g_2(\xi_2,\mu_2)\cdot g(\xi_1+\xi_2,\mu_1+\mu_2,\Omega((\xi_1,\mu_1),(\xi_2,\mu_2)))\,d\xi_1d\xi_2d\mu_1d\mu_2\\
&\leq C2^{j/2}\cdot 2^{-(k_1+k_2+k)/2}\cdot\|g_1\|_{L^2}\|g_2\|_{L^2}\|g\|_{L^2}. 
\end{split}
\end{equation}

We make the changes of variables
\begin{equation*}
\mu_1=\sqrt{3}\xi_1^2+\beta_1\xi_1\text{ and }\mu_2=-\sqrt{3}\xi_2^2+\beta_2\xi_2,
\end{equation*}
with $d\mu_1d\mu_2=\xi_1\xi_2\,d\beta_1d\beta_2$. The left-hand side of \eqref{jj11} is bounded by
\begin{equation}\label{jj12}
\begin{split}
&C2^{k_1+k_2}\int_{S}g_1(\xi_1,\sqrt{3}\xi_1^2+\beta_1\xi_1)\cdot g_2(\xi_2,-\sqrt{3}\xi_2^2+\beta_2\xi_2)\\
&\times g(\xi_1+\xi_2,\sqrt{3}\xi_1^2-\sqrt{3}\xi_2^2+\beta_1\xi_1+\beta_2\xi_2,\widetilde{\Omega}((\xi_1,\beta_1),(\xi_2,\beta_2)))\,d\xi_1d\xi_2d\beta_1d\beta_2,
\end{split}
\end{equation}
where
\begin{equation}\label{jj10}
S=\{(\xi_1,\beta_1,\xi_2,\beta_2):\xi_1+\xi_2\geq 0\text{ and }|\beta_1-\beta_2|\leq 2^{-10}(\xi_1+\xi_2)\},
\end{equation}
and
\begin{equation}\label{jj15}
\widetilde{\Omega}((\xi_1,\beta_1),(\xi_2,\beta_2))=\xi_1\xi_2(\beta_1-\beta_2)\Big(2\sqrt{3}+\frac{\beta_1-\beta_2}{\xi_1+\xi_2}\Big).
\end{equation}

We define the functions $h_i:\R^2\to\R_+$ supported in $I_{k_i}\times\R$, $i=1,2$,
\begin{equation*}
\begin{cases}
&h_1(\xi_1,\beta_1)=2^{k_1/2}\cdot g_1(\xi_1,\sqrt{3}\xi_1^2+\beta_1\xi_1);\\
&h_2(\xi_2,\beta_2)=2^{k_2/2}\cdot g_2(\xi_2,-\sqrt{3}\xi_2^2+\beta_2\xi_2),
\end{cases}
\end{equation*}
with $\|h_i\|_{L^2}\approx\|g_i\|_{L^2}$. Thus, for \eqref{jj5} it suffices to prove that
\begin{equation}\label{jj20}
\begin{split}
&2^{(k_1+k_2)/2}\int_{S}h_1(\xi_1,\beta_1)\cdot h_2(\xi_2,\beta_2)\\
&\times g(\xi_1+\xi_2,\sqrt{3}\xi_1^2-\sqrt{3}\xi_2^2+\beta_1\xi_1+\beta_2\xi_2,\widetilde{\Omega}((\xi_1,\beta_1),(\xi_2,\beta_2)))\,d\xi_1d\xi_2d\beta_1d\beta_2\\
&\leq C2^{j/2}\cdot 2^{-(k_1+k_2+k)/2}\cdot\|h_1\|_{L^2}\|h_2\|_{L^2}\|g\|_{L^2}.
\end{split}
\end{equation}

To prove \eqref{jj20}, we may assume without loss of generality that
\begin{equation}\label{jj11.1}
k_1\leq k_2.
\end{equation}
We make the change of variables $\beta_1=\beta_2+\beta$. In view of \eqref{jj10}, \eqref{jj15}, and the restriction on the support of $g$, we may assume $|\beta|\leq 2^{j-k_1-k_2+4}$. Thus, the integral in the left-hand side of \eqref{jj20} is equal to
\begin{equation}\label{jj21}
\begin{split}
&2^{(k_1+k_2)/2}\int_{\widetilde{S}}h_1(\xi_1,\beta+\beta_2)\cdot h_2(\xi_2,\beta_2)\cdot\mathbf{1}_{[-1,1]}( \beta/2^{j-k_1-k_2+4} )\\
&\times g(\xi_1+\xi_2,A(\xi_1,\xi_2,\beta)+\beta_2(\xi_1+\xi_2),B(\xi_1,\xi_2,\beta))\,d\xi_1d\xi_2d\beta d\beta_2,
\end{split}
\end{equation}
where $\widetilde{S}=\{(\xi_1,\xi_2,\beta,\beta_2)\in\R^4:\xi_1+\xi_2\geq 0\text{ and }|\beta|\leq 2^{-10}(\xi_1+\xi_2)\}$, and
\begin{equation}\label{jj25}
\begin{cases}
&A(\xi_1,\xi_2,\beta)=\sqrt{3}\xi_1^2-\sqrt{3}\xi_2^2+\beta \xi_1;\\
&B(\xi_1,\xi_2,\beta)=\xi_1\xi_2\beta\cdot (2\sqrt{3}+\beta/(\xi_1+\xi_2)).
\end{cases}
\end{equation}
Let $j'=j-k_1-k_2+4$ and decompose, for $i=1,2$, 
\begin{equation*}
h_i(\xi',\beta')=\sum_{m\in\Z}h_i(\xi',\beta')\cdot \mathbf{1}_{[0,1)}(\beta'/2^{j'}-m)=\sum_{m\in\Z}h_i^m(\xi',\beta').
\end{equation*}
The expression in \eqref{jj21} is dominated by to
\begin{equation}\label{jj26}
\begin{split}
&2^{(k_1+k_2)/2}\sum_{|m-m'|\leq 4}\int_{\widetilde{S}}h^m_1(\xi_1,\beta+\beta_2)\cdot h^{m'}_2(\xi_2,\beta_2)\\
&\times g(\xi_1+\xi_2,A(\xi_1,\xi_2,\beta)+\beta_2(\xi_1+\xi_2),B(\xi_1,\xi_2,\beta))\,d\xi_1d\xi_2d\beta d\beta_2.
\end{split}
\end{equation}
Also, for $i=1,2$, 
\begin{equation*}
\|h_i\|_{L^2}=\big[\sum_{m\in\Z}\|h_i^m\|_{L^2}^2\big].
\end{equation*}
Thus, to prove \eqref{jj20}, we may assume $h_1=h_1^m$ and $h_2=h_2^{m'}$ for some fixed $m,m'\in\Z$ with $|m-m'|\leq 4$. To summarize, it suffices to prove that if $F_i:\R^2\to[0,\infty)$ are $L^2$ functions supported in $I_{k_i}\times\R$, $g$ is as before, and $m\in\Z$ then
\begin{equation}\label{jj30}
\begin{split}
&2^{(k_1+k_2)/2}\int_{\widetilde{S}}F_1(\xi_1,\beta+\beta_2)\cdot F_2(\xi_2,\beta_2)\cdot \mathbf{1}_{[m-1,m+1]}(\beta_2/2^{j'})\\
&\times g(\xi_1+\xi_2,A(\xi_1,\xi_2,\beta)+\beta_2(\xi_1+\xi_2),B(\xi_1,\xi_2,\beta))\,d\xi_1d\xi_2d\beta d\beta_2\\
&\leq C2^{j/2}\cdot 2^{-(k_1+k_2+k)/2}\cdot\|F_1\|_{L^2}\|F_2\|_{L^2}\|g\|_{L^2}.
\end{split}
\end{equation}
To prove  \eqref{jj30} we use the Cauchy-Shwartz  inequality in the variables $(\xi_1,\xi_2,\beta)$: with 
\begin{equation*}
S'=\{(\xi_1,\xi_2,\beta)\in\R^3:\xi_i\in I_{k_i},\,\xi_1+\xi_2\geq 0,\,|\beta|\leq 2^{-10}(\xi_1+\xi_2)\},
\end{equation*}
the left-hand side of \eqref{jj30} is dominated by
\begin{equation}\label{jj35}
\begin{split}
C2^{(k_1+k_2)/2}\int_\R\mathbf{1}_{[m-1,m+1]}(\beta_2/2^{j'})\cdot \Big(\int_{S'}|F_1(\xi_1,\beta+\beta_2)\cdot F_2(\xi_2,\beta_2)|^2\,d\xi_1d\xi_2d\beta\Big)^{1/2}\\
\times\Big(\int_{S'}|g(\xi_1+\xi_2,A(\xi_1,\xi_2,\beta)+\beta_2(\xi_1+\xi_2),B(\xi_1,\xi_2,\beta))|^2\,d\xi_1d\xi_2d\beta \Big)^{1/2}\,d\beta_2.
\end{split}
\end{equation}
For \eqref{jj30}, it is easy to see that it suffices to prove that
\begin{equation}\label{jj40}
\begin{split}
\Big(\int_{S'}|g(\xi_1+\xi_2,A(\xi_1,\xi_2,\beta)+\beta_2(\xi_1+\xi_2),B(\xi_1,\xi_2,\beta))|^2\,d\xi_1d\xi_2d\beta \Big)^{1/2}\\
\leq C2^{-(k_1+k_2+k)/2}||g||_{L^2}.
\end{split}
\end{equation}
for any $\beta_2\in\R$. Indeed, assuming \eqref{jj40}, we can bound the expression in \eqref{jj35} by
\begin{equation*}
C2^{(k_1+k_2)/2}\int_\R\mathbf{1}_{[m-1,m+1]}(\beta_2/2^{j'})\cdot ||F_1||_{L^2}||F_2(.,\beta_2)||_{L^2_{\xi_2}}\cdot 2^{-(k_1+k_2+k)/2}||g||_{L^2}\,d\beta_2,
\end{equation*}
which suffices since $2^{j'/2}2^{(k_1+k_2)/2}\approx 2^{j/2}$.

Finally, to prove \eqref{jj40}, we may assume first that $\beta_2=0$. We examine \eqref{jj25} and make the change of variable $\beta=\sqrt{3}(\xi_1+\xi_2)\cdot \nu$. The left-hand side of \eqref{jj40} is dominated by
\begin{equation}\label{jj41}
C\Big(2^{k}\int_{S''}|g(\xi_1+\xi_2,\sqrt{3}(\xi_1+\xi_2)(\xi_1-\xi_2+\nu\xi_1),3\xi_1\xi_2(\xi_1+\xi_2)\nu(2+\nu))|^2\,d\xi_1d\xi_2d\nu \Big)^{1/2},
\end{equation}
where $S''=\{(\xi_1,\xi_2,\nu)\in\R^3:\xi_i\in I_{k_i},\,|\nu|\leq 2^{-10}\}$. We define the function 
\begin{equation*}
h(\xi,x,y)=2^{2k}\cdot |g(\xi,\sqrt{3}\xi\cdot x,3\xi\cdot y)|^2,
\end{equation*}
so $||h|||_{L^1}\approx ||g||_{L^2}^2$. The expression in \eqref{jj41} is dominated by
\begin{equation*}
C2^{-k/2}\Big(\int_{S''}|h(\xi_1+\xi_2,\xi_1-\xi_2+\nu\xi_1,\xi_1\xi_2\cdot \nu(2+\nu))|\,d\xi_1d\xi_2d\nu \Big)^{1/2}.
\end{equation*}
Therefore, it remains to prove that
\begin{equation*}
\int_{S''}|h(\xi_1+\xi_2,\xi_1-\xi_2+\nu\xi_1,\xi_1\xi_2\cdot \nu(2+\nu))|\,d\xi_1d\xi_2d\nu\leq C2^{-(k_1+k_2)}||h||_{L^1}
\end{equation*}
for any function $h\in L^1(\R^3)$. This is clear since the absolute value of the determinant of the change of variables $(\xi_1,\xi_2,\nu)\to [\xi_1+\xi_2,\xi_1-\xi_2+\nu\xi_1,\xi_1\xi_2\cdot \nu(2+\nu)]$ is equal to $(2+\nu)|\xi_1|\cdot |\xi_2(2+\nu)+\xi_1\nu|\approx 2^{k_1+k_2}$, see \eqref{jj11.1} and the definition of the set $S''$.
\end{proof}

\section{Dyadic bilinear estimates I}\label{bilinear2}

In this section we prove the bound \eqref{BIG1} for $k\geq 40$ and $k_1\in[0,k-20]$.

\newtheorem{Lemmak1}{Proposition}[section]
\begin{Lemmak1}\label{Lemmak1}
Assume $k\geq 40$, $k_2\in[k-2,k+2]$, $k_1\in[0,k-20]$, $f_{k_1}\in V_{k_1}\cap W_{k_1}$, $f_{k_2}\in  V_{k_2}\cap W_{k_2}$, and $\mathcal{F}^{-1}(f_{k_1})(x,y,t)$ is supported in $\R^2\times[-2,2]$.
Then
\begin{equation}\label{bt1}
\begin{split}
2^k\big|\big|\chi_k(\xi)&\cdot(\tau-\omega(\xi,\mu)+i)^{-1}\cdot (f_{k_1}\ast
f_{k_2})\big|\big|_{V_k\cap W_k}\\
&\leq
C(2^{-k_1/8}+2^{-(k-k_1)/8})||f_{k_1}||_{V_{k_1}\cap W_{k_1}}||f_{k_2}||_{V_{k_2}\cap W_{k_2}}.
\end{split}
\end{equation}
\end{Lemmak1}

Proposition \ref{Lemmak1} follows from Lemma \ref{Lemmak2}, Lemma \ref{Lemmak3}, and Lemma \ref{Lemmak4} below. We start by decomposing\footnote{In the  decomposition below we make an abuse of notation when we write that  $f_{k_i,2k_i}=
\sum_{l_i<  2k_i+1} f_{k_i,l_i}$. One can see in the rest of the paper that this notation avoids some unnecessary technicalities. One example of its efficiency is in the fact that  for any $l_i<  2k_i+1$
\[(1+|\xi|+|\mu/\xi|)| f_{k_i,l_i}|\sim (1+2^{k_i})| f_{k_i,l_i}|\] 
and hence  we can simply write 
\[(1+|\xi|+|\mu/\xi|)| f_{k_i,2k_i}|\sim (1+2^{k_i})| f_{k_i,2k_i}|.\] 
Our notation also explains why in the proof of the lemmas below we will always assume that $l_1\geq 2k_1$.}
\begin{equation*}
\begin{split}
f_{k_2}=f_{k_2,2k_2-10}+\sum_{l_2\geq 2k_2-9}f_{k_2,l_2}=f_{k_2}\cdot \eta_{\leq 2k_2-10}(\mu_2)+\sum_{l_2\geq 2k_2-9}f_{k_2}\cdot \eta_{l_2}(\mu_2).
\end{split}
\end{equation*}
and
\begin{equation*}
f_{k_1}=f_{k_1,2k_1}+\sum_{l_1\geq 2k_1+1} f_{k_1,l_1}=f_{k_1}\cdot \eta_{\leq 2k_1}(\mu_1)+\sum_{l_1\geq 2k_1+1}f_{k_1}\cdot\eta_{l_1}(\mu_1).
\end{equation*}

Finally for any $J\in\Z$ let $f_{k_i,l_i,J}=f_{k_i,l_i}\cdot \eta_{J}(\tau-\omega(\xi,\mu))$, $f_{k_i,l_i,\leq J}=f_{k_i,l_i}\cdot \eta_{\leq J}(\tau-\omega(\xi,\mu))$, and $f_{k_i,l_i,>J}=f_{k_i,l_i}\cdot \eta_{\geq J+1}(\tau-\omega(\xi,\mu))$, $i=1,2$. 

\newtheorem{Lemmak2}[Lemmak1]{Lemma}
\begin{Lemmak2}\label{Lemmak2}
With the notation in Proposition \ref{Lemmak1}, for any $l_2\in[2k_2-9,2k_2+9]$
\begin{equation*}
\begin{split}
2^k\big|\big|\chi_k(\xi)\cdot&(\tau-\omega(\xi,\mu)+i)^{-1}\cdot (f_{k_1}\ast
f_{k_2,l_2})\big|\big|_{V_k\cap W_k}\\
&\leq C(2^{-k_1/8}+2^{-(k-k_1)/8})||f_{k_1}||_{V_{k_1}\cap W_{k_1}}||f_{k_2,l_2}||_{V_{k_2}\cap W_{k_2}}.
\end{split}
\end{equation*}
\end{Lemmak2}

\begin{proof}[Proof of Lemma \ref{Lemmak2}] In view of the definitions and Lemma \ref{Lemmaa1} (b), it suffices to prove that
\begin{equation*}
\begin{split}
&\big|\big|\chi_k(\xi)\cdot(2^{2k}+i\mu)(\tau-\omega(\xi,\mu)+i)^{-1}\cdot (f_{k_1}\ast
f_{k_2,l_2})\big|\big|_{X_k+Y_k}\\
&+2^k\big|\big|\chi_k(\xi)\cdot(\tau-\omega(\xi,\mu)+i)^{-1}\cdot (f_{k_1}\ast (\partial_\mu+I)f_{k_2,l_2})\big|\big|_{X_k+Y_k}\\
&\leq C(2^{-k_1/8}+2^{-(k-k_1)/8})||f_{k_1}||_{V_{k_1}\cap W_{k_1}}\negmedspace(2^{k}||f_{k_2,l_2}||_{X_{k_2}+Y_{k_2}}+\negmedspace||(\partial_\mu+I)f_{k_2,l_2}||_{X_{k_2}+Y_{k_2}}).
\end{split}
\end{equation*}
For this, it suffices to prove that
\begin{equation}\label{bt2.1}
\begin{split}
\big|\big|\chi_k(\xi)\cdot(2^{k}+&i\mu/2^k)(\tau-\omega(\xi,\mu)+i)^{-1}\cdot (f_{k_1}\ast
f_{k_2,l_2})\big|\big|_{X_k+Y_k}\\
&\leq C(2^{-k_1/8}+2^{-(k-k_1)/8})||f_{k_1}||_{V_{k_1}\cap  W_{k_1}}\cdot ||f_{k_2,l_2}||_{X_{k_2}+Y_{k_2}}.
\end{split}
\end{equation}

In view of Lemma \ref{Lemmaa1} (a) and (b), Lemma \ref{Lemmav3} (a), \eqref{io1}, and the support assumption on $\mathcal{F}^{-1}(f_{k_1})$,
\begin{equation}\label{bj3.1}
\begin{split}
\|\mathcal{F}^{-1}(f_{k_1,l_1,>J})\|_{L^2_yL^\infty_{x,t}}&\leq C\sum_{j> J}2^{j/2}2^{(l_1-k_1)\cdot 3/5}||(I-\partial_{\tau_1}^2)f_{k_1,l_1,j}||_{L^2}\\
&\leq  C2^{(l_1-k_1)\cdot 3/5}(1+2^{(J-2k_1)/2})^{-1}\|(I-\partial_{\tau_1}^2)f_{k_1,l_1}\|_{X_{k_1}}\\
&\leq C(k_1+1)2^{-(l_1-k_1)\cdot 2/5}(1+2^{(J-2k_1)/2})^{-1}\cdot \|f_{k_1}\|_{V_{k_1}}
\end{split}
\end{equation}
for any $l_1\geq 2k_1$ and $J\in\Z\cap [-1,\infty)$. 

We estimate first the contribution of $f_{k_1,l_1}\ast f_{k_2,l_2}$, $2k_1\leq l_1\leq k+k_1-10$. In this range we will show that
\begin{equation}\label{bj125.1}
\begin{split}
2^k\big|\big|\chi_k(\xi)\cdot&(\tau-\omega(\xi,\mu)+i)^{-1}\cdot (f_{k_1,l_1}\ast
f_{k_2,l_2})\big|\big|_{X_k+Y_k}\\
&\leq C2^{-(l_1-k_1)/8}||f_{k_1}||_{V_{k_1}\cap W_{k_1}}\cdot \|f_{k_2,l_2}\|_{X_{k_2}+Y_{k_2}}.
\end{split}
\end{equation}
Let
\begin{equation}\label{j0}
J_0\text{ denote the smallest integer }\geq k-(l_1-k_1)/2-10.
\end{equation}
Using \eqref{sp7.7}, Lemma \ref{Lemmaa1} (a), Lemma \ref{Lemmav1} (a), and \eqref{bj3.1} with $J=-1$, we estimate
\begin{equation}\label{bj10.1}
\begin{split}
2^k\big|\big|&\chi_k(\xi)\cdot(\tau-\omega(\xi,\mu)+i)^{-1}\eta_{\geq J_0+1}(\tau-\omega(\xi,\mu))\cdot (f_{k_1,l_1}\ast
f_{k_2,l_2})\big|\big|_{X_k}\\
&\leq C2^{k}\cdot 2^{-J_0/2}||f_{k_1,l_1}\ast f_{k_2,l_2}||_{L^2_{\xi,\mu,\tau}}\\
&\leq C2^{k}2^{-J_0/2}||\mathcal{F}^{-1}(f_{k_1,l_1})||_{L^2_yL^\infty_{x,t}}\cdot ||\mathcal{F}^{-1}(f_{k_2,l_2})||_{L^\infty_yL^2_{x,t}}\\
&\leq C2^{-(l_1-k_1)/8}||f_{k_1}||_{V_{k_1}}\cdot \|f_{k_2,l_2}\|_{X_{k_2}+Y_{k_2}}.
\end{split}
\end{equation}
We decompose
\begin{equation}\label{hh1}
\begin{split}
f_{k_2,l_2}=f^+_{k_2,l_2,\leq J_0}+f^-_{k_2,l_2,\leq J_0}+f_{k_2,l_2,> J_0}&=f_{k_2,l_2}\cdot \eta_{\leq J_0}(\tau_2-\omega(\xi_2,\mu_2))\mathbf{1}_+(\mu_2)\\
&+f_{k_2,l_2}\cdot \eta_{\leq J_0}(\tau_2-\omega(\xi_2,\mu_2))\mathbf{1}_-(\mu_2)\\
&+f_{k_2,l_2}\cdot \eta_{\geq J_0+1}(\tau_2-\omega(\xi_2,\mu_2))
\end{split}
\end{equation} 
Using \eqref{sp7.8},
\begin{equation}\label{hh2}
\|f_{k_2,l_2,>J_0}\|_{L^2}\leq C2^{-J_0/2}\|f_{k_2,l_2}\|_{X_{k_2}+Y_{k_2}}.
\end{equation}
Thus, using the definitions, Lemma \ref{Lemmaa1} (a), (c), and \eqref{bj3.1} we estimate
\begin{equation}\label{bj11.1}
\begin{split}
2^k\big|\big|&\chi_k(\xi)\cdot(\tau-\omega(\xi,\mu)+i)^{-1}\eta_{\leq J_0}(\tau-\omega(\xi,\mu))\cdot (f_{k_1,l_1}\ast
f_{k_2,l_2,>J_0})\big|\big|_{Y_k}\\
&\leq C2^{k/2}\cdot ||\mathcal{F}^{-1}(f_{k_1,l_1}\ast f_{k_2,l_2,>J_0})||_{L^1_yL^2_{x,t}}\\
&\leq C2^{k/2}||\mathcal{F}^{-1}(f_{k_1,l_1})||_{L^2_yL^\infty_{x,t}}\cdot ||\mathcal{F}^{-1}(f_{k_2,l_2,>J_0})||_{L^2_yL^2_{x,t}}\\
&\leq C2^{-(l_1-k_1)/8}||f_{k_1}||_{V_{k_1}}\cdot \|f_{k_2,l_2}\|_{X_{k_2}+Y_{k_2}}.
\end{split}
\end{equation} 
An estimate similar to \eqref{bj11.1}, using \eqref{bj3.1} gives
\begin{equation}\label{bj12.1}
\begin{split}
2^k\big|\big|&\chi_k(\xi)\cdot(\tau-\omega(\xi,\mu)+i)^{-1}\eta_{\leq J_0}(\tau-\omega(\xi,\mu))\cdot (f_{k_1,l_1,>k+2k_1-10}\ast
f^\pm_{k_2,l_2,\leq J_0})\big|\big|_{Y_k}\\
&\leq C2^{-(l_1-k_1)/4}||f_{k_1}||_{V_{k_1}}\cdot \|f_{k_2,l_2}\|_{X_{k_2}+Y_{k_2}}.
\end{split}
\end{equation} 

It remains to estimate $$2^k\big|\big|\chi_k(\xi)\cdot(\tau-\omega(\xi,\mu)+i)^{-1}\eta_{\leq J_0}(\tau-\omega(\xi,\mu))\cdot (f_{k_1,l_1,\leq k+2k_1-10}\ast
f^\pm_{k_2,l_2,\leq J_0})\big|\big|_{X_k+Y_k}.$$
For $j_2\in\Z_+$ let $f^\pm_{k_2,l_2,j_2}=f_{k_2,l_2}\cdot \eta_{j_2}(\tau_2-\omega(\xi_2,\mu_2))\cdot \mathbf{1}_{\pm}(\mu_2)$. Using Corollary  \ref{Main9}, Lemma \ref{Lemmaa1} (b), and the definitions, we estimate
\begin{equation}\label{bj50.1}
\begin{split}
2^k&\sum_{j_1=J_0+1}^{k+2k_1-10}\big|\big|\chi_k(\xi)\cdot(\tau-\omega(\xi,\mu)+i)^{-1}\eta_{\leq J_0}(\tau-\omega(\xi,\mu))\cdot (f_{k_1,l_1,j_1}\ast
f^\pm_{k_2,l_2,\leq J_0})\big|\big|_{X_k}\\
&\leq C2^{k}\sum_{j_1=J_0+1}^{k+2k_1-10}\sum_{j,j_2=0}^{J_0}2^{-j/2}\big|\big|\eta_{j}(\tau-\omega(\xi,\mu))\cdot (f_{k_1,l_1,j_1}\ast
f^\pm_{k_2,l_2,j_2})\big|\big|_{L^2}\\
&\leq  C2^{k}\sum_{j_1=J_0+1}^{k+2k_1-10}\sum_{j,j_2=0}^{J_0}2^{-(2k+k_1)/2}\cdot 2^{j_1/2}\|f_{k_1,l_1,j_1}\|_{L^2}\cdot 2^{j_2/2}\|f^{\pm}_{k_2,l_2,j_2}\|_{L^2}\\
&\leq C2^{-k_1/2}\cdot k^3\cdot (2^{(J_0-2k_1)/2}+1)^{-1}\cdot 2^{-(l_1-k_1)}||f_{k_1}||_{V_{k_1}}\cdot \|f_{k_2,l_2}\|_{X_{k_2}+Y_{k_2}}\\
&\leq C2^{-(l_1-k_1)/4}||f_{k_1}||_{V_{k_1}}\cdot \|f_{k_2,l_2}\|_{X_{k_2}+Y_{k_2}}.
\end{split}
\end{equation}

Finally, we prove that
\begin{equation}\label{bj40.1}
\begin{split}
2^k\big|\big|&\chi_k(\xi)\cdot(\tau-\omega(\xi,\mu)+i)^{-1}\eta_{\leq J_0}(\tau-\omega(\xi,\mu))\cdot (f_{k_1,l_1,\leq J_0}\ast
f^+_{k_2,l_2,\leq J_0})\big|\big|_{Y_k}\\
&\leq C2^{-(l_1-k_1)/8}||f_{k_1}||_{V_{k_1}\cap W_{k_1}}\cdot \|f_{k_2,l_2}\|_{X_{k_2}+Y_{k_2}}.
\end{split}
\end{equation} 
Recall that (see \eqref{jj2})
\begin{equation}\label{bj13}
\begin{split}
&\Omega[(\xi_1,\mu_1),(\xi_2,\mu_2)]=-\omega(\xi_1+\xi_2,\mu_1+\mu_2)+\omega(\xi_1,\mu_1)+\omega(\xi_2,\mu_2)=-\frac{\xi_1\xi_2}{\xi_1+\xi_2}\\
&\times[(\sqrt{3}\xi_1-\mu_1/ \xi_1)+(\sqrt{3}\xi_2+\mu_2/ \xi_2)]\cdot [(\sqrt{3}\xi_1+\mu_1/ \xi_1)+(\sqrt{3}\xi_2-\mu_2/ \xi_2)].
\end{split}
\end{equation}
Thus, for $\xi_2\in I_{k_2}$, $\mu_2\in[2^{2k-11},2^{2k+11}]$, $\xi_1\in I_{k_1}$, and $|\mu_1|\leq 2^{k-k_1/2-9}$
\begin{equation}\label{bj20.1}
|\Omega[(\xi_1,\mu_1),(\xi_2,\mu_2)]|\geq 2^{k+k_1-4}|(\sqrt{3}\xi_1+\mu_1/ \xi_1)+(\sqrt{3}\xi_2-\mu_2/ \xi_2)|.
\end{equation}

Let $\varphi:\mathbb{R}\to[0,1]$ denote a smooth function supported in $[-1,1]$ with the property that
\begin{equation*}
\sum_{m\in\mathbb{Z}}\varphi(s-m)\equiv 1.
\end{equation*}
Let $\epsilon=2^{-(l_1+k_1)/2}$. For $m\in\mathbb{Z}$ we define
\begin{equation}\label{bj131.1}
\begin{cases}
&f_{k_1,l_1,\leq J_0}^{+,m}(\xi_1,\mu_1,\tau_1)=f_{k_1,l_1,\leq J_0}(\xi_1,\mu_1,\tau_1)\cdot \varphi((\sqrt{3}\xi_1+\mu_1/ \xi_1)/ \epsilon-m);\\
&f_{k_2,l_2,\leq J_0}^{+,m}(\xi_2,\mu_2,\tau_2)=f_{k_2,l_2,\leq J_0}^+(\xi_2,\mu_2,\tau_2)\cdot \varphi((\sqrt{3}\xi_2-\mu_2/ \xi_2)/ \epsilon+m).
\end{cases}
\end{equation}
The important observation is that, in view of \eqref{bj20.1} and the definition of $J_0$,
\begin{equation*}
\eta_{\leq J_0}(\tau-\omega(\xi,\mu))\cdot (f_{k_1,l_1,\leq J_0}^{+,m}\ast
f_{k_2,l_2,\leq J_0}^{+,m'})\equiv 0\text{ unless }|m-m'|\leq 4.
\end{equation*}
Thus, using the definitions and Lemma \ref{Lemmaa1} (c),
\begin{equation}\label{bj30.1}
\begin{split}
&2^k\big|\big|\chi_k(\xi)\cdot(\tau-\omega(\xi,\mu)+i)^{-1}\eta_{\leq J_0}(\tau-\omega(\xi,\mu))\cdot (f_{k_1,l_1,\leq J_0}\ast
f_{k_2,l_2,\leq J_0}^+)\big|\big|_{Y_k}\\
&\leq \negmedspace \sum_{|m-m'|\leq  4}\negmedspace 2^k\big|\big|\chi_k(\xi)(\tau-\omega(\xi,\mu)+i)^{-1}\eta_{\leq J_0}(\tau-\omega(\xi,\mu))\cdot (f_{k_1,l_1,\leq J_0}^{+,m}\ast
f_{k_2,l_2,\leq J_0}^{+,m'})\big|\big|_{Y_k}\\
&\leq C\sum_{|m-m'|\leq  4}2^{k/2}\|\mathcal{F}^{-1}(f_{k_1,l_1,\leq J_0}^{+,m})\|_{L^1_yL^\infty_{x,t}}\cdot\|\mathcal{F}^{-1}(f_{k_2,l_2,\leq J_0}^{+,m'})\|_{L^\infty_yL^2_{x,t}}.
\end{split}
\end{equation}
We use the elementary bound 
\begin{equation}\label{bj120}
\|g\|_{L^1(\mathbb{R})}^2\leq C\|g\|_{L^2(\mathbb{R})}\cdot \|(y+i)\cdot g\|_{L^2(\mathbb{R})}
\end{equation}
for any  $g\in L^2(\R)$, Lemma \ref{Lemmav3} (b), and the definitions to estimate
\begin{equation*}
\begin{split}
\|\mathcal{F}^{-1}&(f_{k_1,l_1,\leq J_0}^{+,m})\|_{L^1_yL^\infty_{x,t}}\leq C2^{(l_1-k_1)/4}\Big(\sum_{j\leq J_0}2^{j/2}\|(I-\partial_{\tau_1}^2)f_{k_1,l_1,j}^{+,m}\|_{L^2}\Big)^{1/2}\\
&\times \Big(\sum_{j\leq J_0}2^{j/2}\|(I-\partial_{\tau_1}^2)(\partial_{\mu_1}+I)f_{k_1,l_1,j}^{+,m}\|_{L^2}\Big)^{1/2}\\
&\leq C2^{(l_1-k_1)/4}\|(1+|\tau_1-\omega(\xi_1,\mu_1)| )^{1/2+1/40}(I-\partial_{\tau_1}^2)f_{k_1,l_1,\leq J_0}^{+,m}\|_{L^2}^{1/2}\\
&\times \|(1+|\tau_1-\omega(\xi_1,\mu_1)| )^{1/2+1/40}(I-\partial_{\tau_1}^2)(\partial_{\mu_1}+I)f_{k_1,l_1,\leq J_0}^{+,m}\|_{L^2}^{1/2},
\end{split}
\end{equation*}
where $f_{k_1,l_1,j}^{+,m}=f_{k_1,l_1,j}\cdot \varphi((\sqrt{3}\xi_1+\mu_1/ \xi_1)/ \epsilon-m)$. Thus, using Lemma \ref{Lemmaa1} (b), with $A=\|(1+|\tau_1-\omega(\xi_1,\mu_1)| )^{1/2+1/40}(I-\partial_{\tau_1}^2)f_{k_1,l_1}\|_{L^2}$
\begin{equation*}
\begin{split}
&\Big[\sum_{m\in\mathbb{Z}}\|\mathcal{F}^{-1}(f_{k_1,l_1,\leq J_0}^{+,m})\|_{L^1_yL^\infty_{x,t}}^2\Big]^{1/2}\\
&\leq C2^{(l_1-k_1)/4}\cdot A^{1/2}\\
&\times \big[\|(1+|\tau_1-\omega(\xi_1,\mu_1)| )^{1/2+1/40}(I-\partial_{\tau_1}^2)(\partial_{\mu_1}+I)f_{k_1,l_1}\|_{L^2}+2^{l_1-k_1}\cdot A\big]^{1/2}\\
&\leq C2^{(l_1-k_1)/4}\cdot 2^{-(l_1-k_1)/2}\cdot (2^{k_1/20}(k_1+1))\|(I-\partial_{\tau_1}^2)f_{k_1,l_1}\|_{V_{k_1}}^{1/2}\\
&\times (2^{k_1/20}(k_1+1))\|(I-\partial_{\tau}^2)f_{k_1,l_1}\|_{V_{k_1}\cap W_{k_1}}^{1/2}.
\end{split}
\end{equation*}
We substitute this last bound into \eqref{bj30.1} and, using Lemma \ref{Lemmav1} (b), \eqref{io1} and $2k_1\leq l_1$, we conclude that the right-hand side of \eqref{bj30.1} is dominated by
\begin{equation*}
\begin{split}
C2^{k/2}\Big[\sum_{m\in\mathbb{Z}}\|\mathcal{F}^{-1}(f_{k_1,l_1,\leq J_0}^{+,m})\|_{L^1_yL^\infty_{x,t}}^2\Big]^{1/2}\cdot \Big[\sum_{m\in\mathbb{Z}}\|\mathcal{F}^{-1}(f_{k_2,l_2,\leq J_0}^{+,m})\|_{L^\infty_yL^2_{x,t}}^2\Big]^{1/2}\\
\leq C2^{-(l_1-k_1)/8}\|f_{k_1}\|_{V_{k_1}\cap W_{k_1}}\cdot \|f_{k_2,l_2}\|_{X_{k_2}+Y_{k_2}}.
\end{split}
\end{equation*}
This gives the bound \eqref{bj40.1}. The bound \eqref{bj125.1} follows from the bounds \eqref{bj10.1}, \eqref{bj11.1}, \eqref{bj12.1}, \eqref{bj50.1}, and \eqref{bj40.1}.

We estimate now the contribution of $f_{k_1,l_1}\ast f_{k_2,l_2}$, $k+k_1-10\leq l_1\leq 2k_2+12$. In this range we will show that
\begin{equation}\label{bj125.2}
\begin{split}
2^k\big|\big|\chi_k(\xi)\cdot&(\tau-\omega(\xi,\mu)+i)^{-1}\cdot (f_{k_1,l_1}\ast
f_{k_2,l_2})\big|\big|_{X_k+Y_k}\\
&\leq C2^{-(l_1-k_1)/4}||f_{k_1}||_{V_{k_1}}\cdot \|f_{k_2,l_2}\|_{X_{k_2}+Y_{k_2}}.
\end{split}
\end{equation}
Using \eqref{sp7.7}, Lemma \ref{Lemmav1} (a), and \eqref{bj3.1} with $J=-1$, we estimate
\begin{equation}\label{bj10.2}
\begin{split}
2^k\big|\big|&\chi_k(\xi)\cdot(\tau-\omega(\xi,\mu)+i)^{-1}\eta_{\geq k-4}(\tau-\omega(\xi,\mu))\cdot (f_{k_1,l_1}\ast
f_{k_2,l_2})\big|\big|_{X_k}\\
&\leq C2^{k}\cdot 2^{-k/2}||f_{k_1,l_1}\ast f_{k_2,l_2}||_{L^2_{\xi,\mu,\tau}}\\
&\leq C2^{k}2^{-k/2}||\mathcal{F}^{-1}(f_{k_1,l_1})||_{L^2_yL^\infty_{x,t}}\cdot ||\mathcal{F}^{-1}(f_{k_2,l_2})||_{L^\infty_yL^2_{x,t}}\\
&\leq C2^{-(l_1-k_1)/4}||f_{k_1}||_{V_{k_1}}\cdot \|f_{k_2,l_2}\|_{X_{k_2}+Y_{k_2}}.
\end{split}
\end{equation}
Using \eqref{hh2}, Lemma \ref{Lemmaa1} (a), (c), and \eqref{bj3.1} we estimate
\begin{equation}\label{bj11.2}
\begin{split}
2^k\big|\big|&\chi_k(\xi)\cdot(\tau-\omega(\xi,\mu)+i)^{-1}\eta_{\leq k-5}(\tau-\omega(\xi,\mu))\cdot (f_{k_1,l_1}\ast
f_{k_2,l_2,>k-5})\big|\big|_{Y_k}\\
&\leq C2^{k/2}\cdot ||\mathcal{F}^{-1}(f_{k_1,l_1}\ast f_{k_2,l_2,>k-5})||_{L^1_yL^2_{x,t}}\\
&\leq C2^{k/2}||\mathcal{F}^{-1}(f_{k_1,l_1})||_{L^2_yL^\infty_{x,t}}\cdot ||\mathcal{F}^{-1}(f_{k_2,l_2,>k-5})||_{L^2_yL^2_{x,t}}\\
&\leq C2^{-(l_1-k_1)/4}||f_{k_1}||_{V_{k_1}}\cdot \|f_{k_2,l_2}\|_{X_{k_2}+Y_{k_2}}.
\end{split}
\end{equation} 
An estimate similar to \eqref{bj12.1} gives
\begin{equation}\label{bj12.2}
\begin{split}
2^k\big|\big|&\chi_k(\xi)\cdot(\tau-\omega(\xi,\mu)+i)^{-1}\eta_{\leq k-5}(\tau-\omega(\xi,\mu))\cdot (f_{k_1,l_1,>k+2k_1-10}\ast
f_{k_2,l_2,\leq k-5})\big|\big|_{Y_k}\\
&\leq C2^{-(l_1-k_1)/4}||f_{k_1}||_{V_{k_1}}\cdot \|f_{k_2,l_2}\|_{X_{k_2}+Y_{k_2}}.
\end{split}
\end{equation} 
Finally,  we  use Corollary \ref{Main9} and Lemma \ref{Lemmaa1} (b) to estimate
\begin{equation}\label{bj50.2}
\begin{split}
2^k&\sum_{j_1=0}^{k+2k_1-10}\big|\big|\chi_k(\xi)\cdot(\tau-\omega(\xi,\mu)+i)^{-1}\eta_{\leq k-5}(\tau-\omega(\xi,\mu))(f_{k_1,l_1,j_1}\ast
f_{k_2,l_2,\leq k-5})\big|\big|_{X_k}\\
&\leq C2^{k}\sum_{j_1=0}^{k+2k_1-10}\sum_{j,j_2=0}^{k-5}2^{-j/2}\big|\big|\eta_{j}(\tau-\omega(\xi,\mu))\cdot (f_{k_1,l_1,j_1}\ast
f_{k_2,l_2,j_2})\big|\big|_{L^2}\\
&\leq  C2^{k}\sum_{j_1=0}^{k+2k_1-10}\sum_{j,j_2=0}^{k-5}2^{-(2k+k_1)/2}\cdot 2^{j_1/2}\|f_{k_1,l_1,j_1}\|_{L^2}\cdot 2^{j_2/2}\|f_{k_2,l_2,j_2}\|_{L^2}\\
&\leq C2^{-k_1/2}\cdot k^3\cdot 2^{-(l_1-k_1)}||f_{k_1}||_{V_{k_1}}\cdot \|f_{k_2,l_2}\|_{X_{k_2}+Y_{k_2}}\\
&\leq C2^{-(l_1-k_1)/2}||f_{k_1}||_{V_{k_1}}\cdot \|f_{k_2,l_2}\|_{X_{k_2}+Y_{k_2}}.
\end{split}
\end{equation}
The bound \eqref{bj125.2} follows from \eqref{bj10.2}, \eqref{bj11.2}, \eqref{bj12.2}, and \eqref{bj50.2}.

We estimate now the contribution of  $\sum_{l_1\geq 2k_2+13}f_{k_1,l_1}\ast f_{k_2,l_2}$: using \eqref{sp7.7} and Lemma \ref{Lemmav5}
\begin{equation}\label{bj6.9}
\begin{split}
\big|\big|&\chi_k(\xi)\cdot (2^k+i\mu/2^k)(\tau-\omega(\xi,\mu)+i)^{-1}\cdot (\sum_{l_1\geq 2k_2+13}f_{k_1,l_1}\ast
f_{k_2,l_2})\big|\big|_{X_k}\\
&\leq C2^{-k}\big|\big|\chi_k(\xi)\cdot\mu \cdot (\sum_{l_1\geq 2k_2+13}f_{k_1,l_1}\ast f_{k_2,l_2})\big|\big|_{L^2}\\
&\leq C2^{-k}\big[\sum_{l_1\geq 2k_2+13}||2^{l_1}f_{k_1,l_1}\ast f_{k_2,l_2}||_{L^2}^2\big]^{1/2}\\
&\leq C2^{-k+k_1}\big[\sum_{ l_1\geq 2k_2+13}||2^{l_1-k_1}f_{k_1,l_1}||_{X_{k_1}+Y_{k_1}}^2\cdot ||f_{k_2,l_2}||_{X_k+Y_k}^2\big]^{1/2}\\
&\leq C2^{k_1-k}||f_{k_1}||_{V_{k_1}}\cdot \|f_{k_2,l_2}\|_{X_{k_2}+Y_{k_2}}.
\end{split}
\end{equation}

The main bound \eqref{bt2.1} follows from \eqref{bj125.1}, \eqref{bj125.2}, and \eqref{bj6.9}.
\end{proof}

\newtheorem{Lemmak3}[Lemmak1]{Lemma}
\begin{Lemmak3}\label{Lemmak3}
With the notation in Proposition \ref{Lemmak1}, 
\begin{equation*}
\begin{split}
2^k\big|\big|\chi_k(\xi)\cdot&(\tau-\omega(\xi,\mu)+i)^{-1}\cdot (f_{k_1}\ast
f_{k_2,2k_2-10})\big|\big|_{V_k\cap W_k}\\
&\leq C(2^{-k_1/8}+2^{-(k-k_1)/8})||f_{k_1}||_{V_{k_1}}\cdot ||f_{k_2,2k_2-10}||_{V_{k_2}\cap W_{k_2}}.
\end{split}
\end{equation*}
\end{Lemmak3}

\begin{proof}[Proof of Lemma \ref{Lemmak3}] As in Lemma \ref{Lemmak2}, it suffices to prove that
\begin{equation}\label{bo1.1}
\begin{split}
\big|\big|\chi_k(\xi)\cdot(2^{k}+&i\mu/2^k)(\tau-\omega(\xi,\mu)+i)^{-1}\cdot (f_{k_1}\ast
f_{k_2,2k_2-10})\big|\big|_{X_k}\\
&\leq C(2^{-k_1/8}+2^{-(k-k_1)/8})||f_{k_1}||_{V_{k_1}}\cdot ||f_{k_2,2k_2-10}||_{X_{k_2}+Y_{k_2}}.
\end{split}
\end{equation}

We estimate first the contribution of $f_{k_1,l_1}\ast f_{k_2,2k_2-10}$, $l_1\in[2k_1,2k+10]$. Let
\begin{equation}\label{bo3.1}
J_0=2k+k_1-40.
\end{equation}
Using \eqref{sp7.7}, \eqref{sp7.8}, and Lemma \ref{Lemmav5}, we estimate
\begin{equation}\label{rr1.1}
\begin{split}
2^{k}\big|\big|\chi_k(\xi)&\cdot(\tau-\omega(\xi,\mu)+i)^{-1}\eta_{\geq 2k-39}(\tau-\omega(\xi,\mu))\cdot (f_{k_1,l_1}\ast f_{k_2,2k_2-10})\big|\big|_{X_k}\\
&\leq C2^k2^{-k}\big|\big|f_{k_1,l_1}\ast f_{k_2,2k_2-10}\big|\big|_{L^2}\\
&\leq C||\mathcal{F}^{-1}(f_{k_1,l_1})||_{L^4}\cdot ||\mathcal{F}^{-1}(f_{k_2,2k_2-10})||_{L^4}\\
&\leq C2^{k_1-l_1}||f_{k_1}||_{V_{k_1}}\cdot ||f_{k_2,2k_2-10}||_{X_{k_2}+Y_{k_2}}.
\end{split}
\end{equation}
We have the $L^\infty$ bound
\begin{equation}\label{nj1}
\begin{split}
||\mathcal{F}^{-1}(f_{k_1,l_1})||_{L^\infty}&\leq C\sum_{j\geq 0}2^{j/2}2^{k_1/2}2^{l_1/2}||f_{k_1,l_1,j}||_{L^2}\\
&\leq C(k_1+1)2^{3k_1/2}2^{-l_1/2}||f_{k_1}||_{V_{k_1}}.
\end{split}
\end{equation}
Thus, using \eqref{sp7.7} and \eqref{sp7.8}, we estimate
\begin{equation}\label{rr2.1}
\begin{split}
2^k\big|\big|\chi_k(\xi)\cdot&(\tau-\omega(\xi,\mu)+i)^{-1}\eta_{\leq 2k-40}(\tau-\omega(\xi,\mu))\cdot (f_{k_1,l_1}\ast f_{k_2,2k_2-10,>J_0})\big|\big|_{X_k}\\
&\leq C2^k\big|\big|f_{k_1,l_1}\ast f_{k_2,2k_2-10,>J_0}\big|\big|_{L^2}\\
&\leq C2^k\|\mathcal{F}^{-1}(f_{k_1,l_1})\|_{L^\infty}\cdot \|f_{k_2,2k_2-10,>J_0}\|_{L^2}\\
&\leq C(k_1+1)2^{(k_1-l_1)/2}||f_{k_1}||_{V_{k_1}}\cdot ||f_{k_2,2k_2-10}||_{X_{k_2}+Y_{k_2}}.
\end{split}
\end{equation}
As in the proof of Lemma \ref{Lemmav5} (see  \eqref{re1}) we have
\begin{equation*}
\|\mathcal{F}^{-1}(f_{k_1,l_1,>J_0})\|_{L^4}\leq C\sum_{j\geq J_0+1}2^{j/2}\|f_{k_1,l_1,j}\|_{L^2}\leq C2^{-(2k-k_1)/2}2^{k_1-l_1}\|f_{k_1}\|_{V_{k_1}}.
\end{equation*}
Thus, using \eqref{sp7.7} and \eqref{sp7.8} and Lemma \ref{Lemmav5},
\begin{equation}\label{rr1.2}
\begin{split}
2^k\big|\big|\chi_k(\xi)&\cdot(\tau-\omega(\xi,\mu)+i)^{-1}\eta_{\leq  2k-40}(\tau-\omega(\xi,\mu))\cdot (f_{k_1,l_1,>J_0}\ast f_{k_2,2k_2-10,\leq  J_0})\big|\big|_{X_k}\\
&\leq C2^k\big|\big|f_{k_1,l_1,>J_0}\ast f_{k_2,2k_2-10,\leq J_0}\big|\big|_{L^2}\\
&\leq C2^k||\mathcal{F}^{-1}(f_{k_1,l_1,>J_0})||_{L^4}\cdot ||\mathcal{F}^{-1}(f_{k_2,2k_2-10,\leq J_0})||_{L^4}\\
&\leq C2^{3 k_1/2-l_1}||f_{k_1}||_{V_{k_1}}\cdot ||f_{k_2,2k_2-10}||_{X_{k_2}+Y_{k_2}}.
\end{split}
\end{equation}
Finally, we observe that
\begin{equation*}
\eta_{\leq 2k-40}(\tau-\omega(\xi,\mu))\cdot (f_{k_1,l_1,\leq J_0}\ast
f_{k_2,2k_2-10,\leq J_0})\equiv 0,
\end{equation*}
unless $l_1\in[k+k_1-10,k+k_1+10]$, which is a consequence of the identity \eqref{bj13}. Using Corollary \ref{Main9}, Lemma \ref{Lemmaa1} (b), and the definitions, we estimate
\begin{equation}\label{rr3.1}
\begin{split}
2^k&\big|\big|\chi_k(\xi)\cdot(\tau-\omega(\xi,\mu)+i)^{-1}\eta_{\leq 2k-40}(\tau-\omega(\xi,\mu))\cdot (f_{k_1,l_1,\leq J_0}\ast
f_{k_2,2k_2-10,\leq J_0})\big|\big|_{X_k}\\
&\leq C2^{k}\sum_{j_1,j_2=0}^{J_0}\sum_{j=0}^{2k-40}2^{-j/2}\big|\big|\eta_{j}(\tau-\omega(\xi,\mu))\cdot (f_{k_1,l_1,j_1}\ast
f_{k_2,2k_2-10,j_2})\big|\big|_{L^2}\\
&\leq  C2^{k}\sum_{j,j_1,j_2=0}^{J_0}2^{-(2k+k_1)/2}\cdot 2^{j_1/2}\|f_{k_1,l_1,j_1}\|_{L^2}\cdot 2^{j_2/2}\|f_{k_2,2k_2-10,j_2}\|_{L^2}\\
&\leq C2^{-k_1/2}\cdot k^3\cdot 2^{-(l_1-k_1)}||f_{k_1}||_{V_{k_1}}\cdot \|f_{k_2,2k_2-10}\|_{X_{k_2}+Y_{k_2}}.
\end{split}
\end{equation}
Thus, using \eqref{rr1.1}, \eqref{rr2.1}, \eqref{rr1.2}, and \eqref{rr3.1} with $l_1\in[k+k_1-10,k+k_1+10]$, we have
\begin{equation}\label{bo4.1}
\begin{split}
\sum_{l_1=2k_1}^{2k+10}2^k\big|\big|&\chi_k(\xi)\cdot(\tau-\omega(\xi,\mu)+i)^{-1}\cdot (f_{k_1,l_1}\ast
f_{k_2,2k_2-10})\big|\big|_{X_k}\\
&\leq C2^{-k_1/4}||f_{k_1}||_{V_{k_1}}\cdot ||f_{k_2,2k_2-10}||_{X_{k_2}+Y_{k_2}}.
\end{split}
\end{equation}

We estimate now the contribution of  $\sum_{l_1\geq 2k+11}f_{k_1,l_1}\ast f_{k_2,2k_2-10}$: using \eqref{sp7.7} and Lemma \ref{Lemmav5}, we estimate as in \eqref{bj6.9}
\begin{equation}\label{rr9.1}
\begin{split}
\big|\big|\chi_k(\xi)\cdot &(2^k+i\mu/2^k)(\tau-\omega(\xi,\mu)+i)^{-1}\cdot (\sum_{l_1\geq 2k+11}f_{k_1,l_1}\ast
f_{k_2,2k_2-10})\big|\big|_{X_k}\\
&\leq C2^{k_1-k}||f_{k_1}||_{V_{k_1}}\cdot \|f_{k_2,2k_2-10}\|_{X_{k_2}+Y_{k_2}}.
\end{split}
\end{equation}

The main bound \eqref{bo1.1} follows from \eqref{bo4.1} and \eqref{rr9.1}.
\end{proof}

\newtheorem{Lemmak4}[Lemmak1]{Lemma}
\begin{Lemmak4}\label{Lemmak4}
With the notation in Proposition \ref{Lemmak1}, for any $l_2\geq 2k_2+10$
\begin{equation*}
\begin{split}
2^k\big|&\big|\chi_k(\xi)\cdot(\tau-\omega(\xi,\mu)+i)^{-1}\cdot (f_{k_1}\ast
f_{k_2,l_2})\big|\big|_{V_k\cap W_k}\\
&\leq C2^{-(l_2-2k_2)/4}(2^{-k_1/8}+2^{-(k-k_1)/8})||f_{k_1}||_{V_{k_1}}\cdot ||f_{k_2,l_2}||_{V_{k_2}\cap W_{k_2}}.
\end{split}
\end{equation*}
\end{Lemmak4}

\begin{proof}[Proof of Lemma \ref{Lemmak4}] As in Lemma \ref{Lemmak2}, it suffices to prove that
\begin{equation}\label{go1.1}
\begin{split}
\big|\big|&\chi_k(\xi)\cdot(2^{l_2-k}+i\mu/2^k)(\tau-\omega(\xi,\mu)+i)^{-1}\cdot (f_{k_1}\ast
f_{k_2,l_2})\big|\big|_{X_k}\\
&\leq C2^{(l_2-2k_2)\cdot (3/4)}(2^{-k_1/8}+2^{-(k-k_1)/8})||f_{k_1}||_{V_{k_1}}\cdot ||f_{k_2,l_2}||_{X_{k_2}+Y_{k_2}}.
\end{split}
\end{equation}

We estimate first the contribution of $f_{k_1,l_1}\ast f_{k_2,l_2}$, for
\begin{equation}\label{l1}
l_1\in[2k_1,l_2+10]\setminus [l_2-k_2+k_1-10,l_2-k_2+k_1+10].
\end{equation}
Let
\begin{equation}\label{go3.1}
J_0=l_2+k_1-40.
\end{equation}
Using \eqref{sp7.7}, \eqref{bj3.1}, and Lemma \ref{Lemmav1}, we estimate
\begin{equation}\label{gr1.1}
\begin{split}
2^{l_2-k}\big|\big|\chi_k(\xi)&\cdot(\tau-\omega(\xi,\mu)+i)^{-1}\eta_{\geq J_0+1}(\tau-\omega(\xi,\mu))\cdot (f_{k_1,l_1}\ast f_{k_2,l_2})\big|\big|_{X_k}\\
&\leq C2^{l_2-k}2^{-k}\big|\big|f_{k_1,l_1}\ast f_{k_2,l_2}\big|\big|_{L^2}\\
&\leq C2^{l_2-2k_2}||\mathcal{F}^{-1}(f_{k_1,l_1})||_{L^2_yL^\infty_{x,t}}\cdot ||\mathcal{F}^{-1}(f_{k_2,l_2})||_{L^\infty_yL^2_{x,t}}\\
&\leq C2^{(l_2-2k_2)/2}2^{-(l_1-k_1)/4}||f_{k_1}||_{V_{k_1}}\cdot ||f_{k_2,l_2}||_{X_{k_2}+Y_{k_2}}.
\end{split}
\end{equation}
Using the $L^\infty$ bound \eqref{nj1}, \eqref{sp7.7}, and \eqref{sp7.8}, we estimate
\begin{equation}\label{gr2.1}
\begin{split}
2^{l_2-k}\big|\big|\chi_k(\xi)\cdot&(\tau-\omega(\xi,\mu)+i)^{-1}\eta_{\leq J_0}(\tau-\omega(\xi,\mu))\cdot (f_{k_1,l_1}\ast f_{k_2,l_2,>J_0})\big|\big|_{X_k}\\
&\leq C2^{l_2-k}\big|\big|f_{k_1,l_1}\ast f_{k_2,l_2,>J_0}\big|\big|_{L^2}\\
&\leq C2^{l_2-k}\|\mathcal{F}^{-1}(f_{k_1,l_1})\|_{L^\infty}\cdot \|f_{k_2,l_2,>J_0}\|_{L^2}\\
&\leq C(k_1+1)2^{(k_1-l_1)/2}||f_{k_1}||_{V_{k_1}}\cdot ||f_{k_2,l_2}||_{X_{k_2}+Y_{k_2}}.
\end{split}
\end{equation}
Using \eqref{sp7.7}, \eqref{bj3.1}, and Lemma \ref{Lemmav1},
\begin{equation}\label{gr1.2}
\begin{split}
2^{l_2-k}\big|\big|\chi_k(\xi)&\cdot(\tau-\omega(\xi,\mu)+i)^{-1}\eta_{\leq  J_0}(\tau-\omega(\xi,\mu))\cdot (f_{k_1,l_1,>J_0}\ast f_{k_2,l_2,\leq  J_0})\big|\big|_{X_k}\\
&\leq C2^{l_2-k}\big|\big|f_{k_1,l_1,>J_0}\ast f_{k_2,l_2,\leq J_0}\big|\big|_{L^2}\\
&\leq C2^{l_2-k}||\mathcal{F}^{-1}(f_{k_1,l_1,>J_0})||_{L^2_yL^\infty_{x,t}}\cdot ||\mathcal{F}^{-1}(f_{k_2,l_2,\leq J_0})||_{L^\infty_yL^2_{x,t}}\\
&\leq C2^{-(l_1-k_1)/4}||f_{k_1}||_{V_{k_1}}\cdot ||f_{k_2,l_2}||_{X_{k_2}+Y_{k_2}}.
\end{split}
\end{equation}
Finally, we observe that for $l_1$ as in \eqref{l1}
\begin{equation*}
\eta_{\leq J_0}(\tau-\omega(\xi,\mu))\cdot (f_{k_1,l_1,\leq J_0}\ast
f_{k_2,l_2,\leq J_0})\equiv 0,
\end{equation*}
which is a consequence of the identity \eqref{bj13}. Thus, for $l_1$ as in \eqref{l1},
\begin{equation}\label{gr3.1}
\begin{split}
2^{l_2-k}\big|\big|\chi_k(\xi)&\cdot(\tau-\omega(\xi,\mu)+i)^{-1}\cdot (f_{k_1,l_1}\ast f_{k_2,l_2})\big|\big|_{X_k}\\
&\leq C2^{(l_2-2k_2)/2}2^{-(l_1-k_1)/4}||f_{k_1}||_{V_{k_1}}\cdot ||f_{k_2,l_2}||_{X_{k_2}+Y_{k_2}}.
\end{split}
\end{equation}

We estimate now the contribution of $f_{k_1,l_1}\ast f_{k_2,l_2}$, for
\begin{equation}\label{l2}
l_1\in[l_2-k_2+k_1-10,l_2-k_2+k_1+10].
\end{equation}
Let
\begin{equation*}
J_1=2k+k_1-40.
\end{equation*}
As in \eqref{gr1.1}, \eqref{gr2.1}, and \eqref{gr1.2}, using also $2^{k_1-l_1}\approx 2^{k_2-l_2}$, we estimate
\begin{equation}\label{gr4.1}
\begin{split}
&2^{l_2-k}\big|\big|\chi_k(\xi)\cdot(\tau-\omega(\xi,\mu)+i)^{-1}\eta_{\geq 2k-39}(\tau-\omega(\xi,\mu))\cdot (f_{k_1,l_1}\ast f_{k_2,l_2})\big|\big|_{X_k}\\
&+2^{l_2-k}\big|\big|\chi_k(\xi)\cdot(\tau-\omega(\xi,\mu)+i)^{-1}\eta_{\leq 2k-40}(\tau-\omega(\xi,\mu))\cdot (f_{k_1,l_1}\ast f_{k_2,l_2,>J_1})\big|\big|_{X_k}\\
&+2^{l_2-k}\big|\big|\chi_k(\xi)\cdot (\tau-\omega(\xi,\mu)+i)^{-1}\eta_{\leq 2k-40}(\tau-\omega(\xi,\mu))\cdot (f_{k_1,l_1,>J_1}\ast f_{k_2,l_2,\leq  J_1})\big|\big|_{X_k}\\
&\leq C 2^{(l_2-2k_2)/2}2^{-k/4}||f_{k_1}||_{V_{k_1}}\cdot ||f_{k_2,l_2}||_{X_{k_2}+Y_{k_2}}.
\end{split}
\end{equation}
In addition, using Corollary  \ref{Main9}, Lemma \ref{Lemmaa1} (b), and the definitions, we estimate
\begin{equation}\label{gr4.2}
\begin{split}
2^{l_2-k}&\big|\big|\chi_k(\xi)\cdot(\tau-\omega(\xi,\mu)+i)^{-1}\eta_{\leq 2k-40}(\tau-\omega(\xi,\mu))\cdot (f_{k_1,l_1,\leq J_1}\ast
f_{k_2,l_2,\leq J_1})\big|\big|_{X_k}\\
&\leq C2^{l_2-k}\sum_{j_1,j_2=0}^{J_1}\sum_{j=0}^{2k-40}2^{-j/2}\big|\big|\eta_{j}(\tau-\omega(\xi,\mu))\cdot (f_{k_1,l_1,j_1}\ast
f_{k_2,l_2,j_2})\big|\big|_{L^2}\\
&\leq  C2^{l_2-k}\sum_{j,j_1,j_2=0}^{J_1}2^{-(2k+k_1)/2}\cdot 2^{j_1/2}\|f_{k_1,l_1,j_1}\|_{L^2}\cdot 2^{j_2/2}\|f_{k_2,l_2,j_2}\|_{L^2}\\
&\leq C2^{-k_1/2}\cdot k^3\cdot 2^{-k}||f_{k_1}||_{V_{k_1}}\cdot \|f_{k_2,l_2}\|_{X_{k_2}+Y_{k_2}}.
\end{split}
\end{equation}

We estimate now the contribution of  $\sum_{l_1\geq l_2+11}f_{k_1,l_1}\ast f_{k_2,l_2}$: using \eqref{sp7.7} and Lemma \ref{Lemmav5}, we estimate as in \eqref{bj6.9}
\begin{equation}\label{gr9.1}
\begin{split}
\big|\big|\chi_k(\xi)\cdot &(2^{l_2-k}+i\mu/2^k)(\tau-\omega(\xi,\mu)+i)^{-1}\cdot (\sum_{l_1\geq l_2+11}f_{k_1,l_1}\ast
f_{k_2,l_2})\big|\big|_{X_k}\\
&\leq C2^{-k}\big|\big|\mu\cdot (\sum_{l_1\geq l_2+11}f_{k_1,l_1}\ast
f_{k_2,l_2})\big|\big|_{L^2}\\
&\leq C2^{k_1-k}\big[\sum_{l_1\geq l_2+11}\|2^{l_1-k_1}f_{k_1,l_1}\ast f_{k_2,l_2}\|_{L^2}^2\big]^{1/2}\\
&\leq C2^{k_1-k}\big[\sum_{l_1\geq l_2+11}\|2^{l_1-k_1}f_{k_1,l_1}\|_{X_{k_1}+Y_{k_1}}^2\cdot \|f_{k_2,l_2}\|_{X_{k_2}+Y_{k_2}}^2\big]^{1/2}\\
&\leq C2^{k_1-k}||f_{k_1}||_{V_{k_1}}\cdot \|f_{k_2,l_2}\|_{X_{k_2}+Y_{k_2}}.
\end{split}
\end{equation}

The main bound \eqref{go1.1} follows from \eqref{gr3.1}, \eqref{gr4.1}, \eqref{gr4.2}, and \eqref{gr9.1}.
\end{proof}

\section{Dyadic bilinear estimates II}\label{bilinear1}

In this section we prove the bound \eqref{BIG1} for $k\geq 40$ and $k_1\leq 0$.

\newtheorem{Lemmac1}{Proposition}[section]
\begin{Lemmac1}\label{Lemmac1}
Assume $k\geq 40$, $k_2\in[k-2,k+2]$, $k_1\leq0$, $f_{k_1}\in  V_{k_1}\cap W_{k_1}$, $f_{k_2}\in  V_{k_2}\cap W_{k_2}$, and $\mathcal{F}^{-1}(f_{k_1})$ is supported in $\R^2\times[-2,2]$.
Then
\begin{equation}\label{bj1}
\begin{split}
2^k\big|\big|\chi_k&(\xi)\cdot(\tau-\omega(\xi,\mu)+i)^{-1}\cdot (f_{k_1}\ast
f_{k_2})\big|\big|_{V_k\cap W_k}\\
&\leq C2^{k_1/4}||f_{k_1}||_{V_{k_1}\cap W_{k_1}}||f_{k_2}||_{V_{k_2}\cap W_{k_2}}.
\end{split}
\end{equation}
\end{Lemmac1}

Proposition \ref{Lemmac1} follows from Lemma \ref{Lemmac2}, Lemma \ref{Lemmac3}, and Lemma \ref{Lemmac4} below. We start by decomposing\footnote{In the  decomposition below we make an abuse of notation when we write that  $f_{k_2,2k_2}=
\sum_{l_i<  2k_2+1} f_{k_2,l_2}$ and $f_{k_1,k_1}=
\sum_{l_1<  k_1+1} f_{k_1,l_1}$ . One can see in the rest of the paper that this notation avoids some unnecessary technicalities. One example of its efficiency is in the fact that  for any $k_1\leq 0, \, l_1<  k_1+1$
\[(1+|\xi|+|\mu/\xi|)| f_{k_1,l_1}|\sim | f_{k_2,l_2}|\] 
and hence  we can simply write 
\[(1+|\xi|+|\mu/\xi|)| f_{k_1,k_1}|\sim | f_{k_1,k_1}|.\] 
Our notation also explains why in the proof of the lemmas below we will always assume that $l_1\geq k_1$.}

\begin{equation*}
\begin{split}
f_{k_2}=f_{k_2,2k_2-10}+\sum_{l_2\geq 2k_2-9}f_{k_2,l_2}=f_{k_2}\cdot \eta_{\leq 2k_2-10}(\mu_2)+\sum_{l_2\geq 2k_2-9}f_{k_2}\cdot \eta_{l_2}(\mu_2).
\end{split}
\end{equation*}
and
\begin{equation*}
f_{k_1}=f_{k_1,k_1}+\sum_{l_1\geq k_1+1} f_{k_1,l_1}=f_{k_1}\cdot \eta_{0}(\mu_1/2^{k_1})+\sum_{l_1\geq k_1+1}f_{k_1}\cdot\eta_{l_1}(\mu_1).
\end{equation*}
For any $J\in\Z$ let $f_{k_i,l_i,J}=f_{k_i,l_i}\cdot \eta_{J}(\tau-\omega(\xi,\mu))$, $f_{k_i,l_i,\leq J}=f_{k_i,l_i}\cdot \eta_{\leq J}(\tau-\omega(\xi,\mu))$, and $f_{k_i,l_i,>J}=f_{k_i,l_i}\cdot \eta_{\geq J+1}(\tau-\omega(\xi,\mu))$, $i=1,2$.

\newtheorem{Lemmac2}[Lemmac1]{Lemma}
\begin{Lemmac2}\label{Lemmac2}
With the notation in Proposition \ref{Lemmac1}, for any $l_2\in[2k_2-9,2k_2+9]$
\begin{equation*}
\begin{split}
2^k\big|\big|\chi_k(\xi)\cdot&(\tau-\omega(\xi,\mu)+i)^{-1}\cdot (f_{k_1}\ast
f_{k_2,l_2})\big|\big|_{V_k\cap W_k}\\
&\leq C2^{k_1/4}||f_{k_1}||_{V_{k_1}\cap W_{k_1}}||f_{k_2,l_2}||_{V_{k_2}\cap W_{k_2}}.
\end{split}
\end{equation*}
\end{Lemmac2}

\begin{proof}[Proof of Lemma \ref{Lemmac2}] In view of the definitions and Lemma \ref{Lemmaa1} (b), it suffices to prove that
\begin{equation}\label{bj2}
\begin{split}
&\big|\big|\chi_k(\xi)\cdot(2^{2k}+i\mu)(\tau-\omega(\xi,\mu)+i)^{-1}\cdot (f_{k_1}\ast
f_{k_2,l_2})\big|\big|_{X_k+Y_k}\\
&+2^k\big|\big|\chi_k(\xi)\cdot(\tau-\omega(\xi,\mu)+i)^{-1}\cdot (f_{k_1}\ast (\partial_\mu+I)f_{k_2,l_2})\big|\big|_{X_k+Y_k}\\
&\leq C2^{k_1/4}||f_{k_1}||_{V_{k_1}\cap W_{k_1}}\cdot (2^{k}||f_{k_2,l_2}||_{X_{k_2}+Y_{k_2}}+||(\partial_\mu+I)f_{k_2,l_2}||_{X_{k_2}+Y_{k_2}}).
\end{split}
\end{equation}
For this, it suffices to prove that
\begin{equation}\label{bj2.1}
\begin{split}
\big|\big|\chi_k(\xi)\cdot(2^{k}+&i\mu/2^k)(\tau-\omega(\xi,\mu)+i)^{-1}\cdot (f_{k_1}\ast
f_{k_2,l_2})\big|\big|_{X_k+Y_k}\\
&\leq C2^{k_1/4}||f_{k_1}||_{V_{k_1}\cap W_{k_1}}\cdot ||f_{k_2,l_2}||_{X_{k_2}+Y_{k_2}}.
\end{split}
\end{equation}

In view of Lemma \ref{Lemmaa1} (a), Lemma \ref{Lemmav2} (a), \eqref{io1}, and the support assumption on $\mathcal{F}^{-1}(f_{k_1})$,
\begin{equation}\label{bj3}
\begin{cases}
&\|\mathcal{F}^{-1}(\sum_{l_1=k_1}^0f_{k_1,l_1})\|_{L^2_yL^\infty_{x,t}}\leq C2^{k_1/4}||f_{k_1}||_{V_{k_1}};\\
&\|\mathcal{F}^{-1}(f_{k_1,l_1})\|_{L^2_yL^\infty_{x,t}}\leq C2^{(2l_1+k_1)/4}\cdot 2^{k_1-l_1}||f_{k_1}||_{V_{k_1}}\text{ for }l_1\geq 1.
\end{cases}
\end{equation}
Also, using the elementary inequality \eqref{bj120},
\begin{equation}\label{bj4}
\|\mathcal{F}^{-1}(\sum_{l_1=k_1}^0f_{k_1,l_1})\|_{L^1_yL^\infty_{x,t}}\leq C2^{k_1/4}||f_{k_1}||_{V_{k_1}\cap W_{k_1}}.
\end{equation}

We start by estimating the contribution of $\sum_{l_1=k_1}^0f_{k_1,l_1}\ast f_{k_2,l_2}$. Using the definitions, Lemma \ref{Lemmaa1} (a), (c), Lemma \ref{Lemmav1} (a), \eqref{bj3}, and \eqref{bj4}, we estimate
\begin{equation*}
\begin{split}
2^k\big|\big|&\chi_k(\xi)\cdot(\tau-\omega(\xi,\mu)+i)^{-1}\eta_{\leq k-10}(\tau-\omega(\xi,\mu))\cdot (\sum_{l_1=k_1}^0f_{k_1,l_1}\ast
f_{k_2,l_2})\big|\big|_{Y_k}\\
&\leq C2^{k/2}||\mathcal{F}^{-1}(\sum_{l_1=k_1}^0f_{k_1,l_1}\ast f_{k_2,l_2})||_{L^1_yL^2_{x,t}}\\
&\leq C2^{k/2}||\mathcal{F}^{-1}(\sum_{l_1=k_1}^0f_{k_1,l_1})||_{L^1_yL^\infty_{x,t}}\cdot ||\mathcal{F}^{-1}(f_{k_2,l_2})||_{L^\infty_yL^2_{x,t}}\\
&\leq C2^{k_1/4}||f_{k_1}||_{V_{k_1}\cap W_{k_1}}\cdot \|f_{k_2,l_2}\|_{X_{k_2}+Y_{k_2}}.
\end{split}
\end{equation*}
In addition, using \eqref{sp7.7},
\begin{equation*}
\begin{split}
2^k\big|\big|&\chi_k(\xi)\cdot(\tau-\omega(\xi,\mu)+i)^{-1}\eta_{\geq k-9}(\tau-\omega(\xi,\mu))\cdot (\sum_{l_1=k_1}^0f_{k_1,l_1}\ast
f_{k_2,l_2})\big|\big|_{X_k}\\
&\leq C2^{k}\cdot 2^{-k/2}||\sum_{l_1=k_1}^0f_{k_1,l_1}\ast f_{k_2,l_2}||_{L^2_{\xi,\mu,\tau}}\\
&\leq C2^{k/2}||\mathcal{F}^{-1}(\sum_{l_1=k_1}^0f_{k_1,l_1})||_{L^2_yL^\infty_{x,t}}\cdot ||\mathcal{F}^{-1}(f_{k_2,l_2})||_{L^\infty_yL^2_{x,t}}\\
&\leq C2^{k_1/4}||f_{k_1}||_{V_{k_1}}\cdot \|f_{k_2,l_2}\|_{X_{k_2}+Y_{k_2}}.
\end{split}
\end{equation*}
Thus
\begin{equation}\label{bj5}
\begin{split}
\big|\big|\chi_k(\xi)\cdot&(2^k+i\mu/ 2^k)(\tau-\omega(\xi,\mu)+i)^{-1}\cdot (\sum_{l_1=k_1}^0f_{k_1,l_1}\ast
f_{k_2,l_2})\big|\big|_{X_k+Y_k}\\
&\leq C2^{k_1/4}||f_{k_1}||_{V_{k_1}\cap W_{k_1}}\cdot ||f_{k_2,l_2}||_{X_{k_2}+Y_{k_2}}.
\end{split}
\end{equation}

We estimate now the contribution of $f_{k_1,l_1}\ast f_{k_2,l_2}$, $1\leq l_1\leq k+2k_1-10$. In this range we will show that
\begin{equation}\label{bj125}
\begin{split}
2^k\big|\big|\chi_k(\xi)\cdot&(\tau-\omega(\xi,\mu)+i)^{-1}\cdot (f_{k_1,l_1}\ast
f_{k_2,l_2})\big|\big|_{X_k+Y_k}\\
&\leq C2^{3k_1/4}2^{-l_1/4}||f_{k_1}||_{V_{k_1}\cap W_{k_1}}\cdot \|f_{k_2,l_2}\|_{X_{k_2}+Y_{k_2}}.
\end{split}
\end{equation}
Let
\begin{equation*}
J_0\text{ denote the smallest integer }\geq k+k_1-l_1/2-10.
\end{equation*}
Using \eqref{sp7.7}, Lemma \ref{Lemmaa1} (a), Lemma \ref{Lemmav1} (a), and \eqref{bj3} we estimate
\begin{equation}\label{bj10}
\begin{split}
2^k\big|\big|&\chi_k(\xi)\cdot(\tau-\omega(\xi,\mu)+i)^{-1}\eta_{\geq J_0+1}(\tau-\omega(\xi,\mu))\cdot (f_{k_1,l_1}\ast
f_{k_2,l_2})\big|\big|_{X_k}\\
&\leq C2^{k}\cdot 2^{-J_0/2}||f_{k_1,l_1}\ast f_{k_2,l_2}||_{L^2_{\xi,\mu,\tau}}\\
&\leq C2^{k}2^{-J_0/2}||\mathcal{F}^{-1}(f_{k_1,l_1})||_{L^2_yL^\infty_{x,t}}\cdot ||\mathcal{F}^{-1}(f_{k_2,l_2})||_{L^\infty_yL^2_{x,t}}\\
&\leq C2^{3k_1/4}2^{-l_1/4}||f_{k_1}||_{V_{k_1}}\cdot \|f_{k_2,l_2}\|_{X_{k_2}+Y_{k_2}}.
\end{split}
\end{equation}
As in \eqref{hh1}, we decompose
\begin{equation*}
\begin{split}
f_{k_2,l_2}=f^+_{k_2,l_2,\leq J_0}+f^-_{k_2,l_2,\leq J_0}+f_{k_2,l_2,> J_0}&=f_{k_2,l_2}\cdot \eta_{\leq J_0}(\tau_2-\omega(\xi_2,\mu_2))\mathbf{1}_+(\mu_2)\\
&+f_{k_2,l_2}\cdot \eta_{\leq J_0}(\tau_2-\omega(\xi_2,\mu_2))\mathbf{1}_-(\mu_2)\\
&+f_{k_2,l_2}\cdot \eta_{\geq J_0+1}(\tau_2-\omega(\xi_2,\mu_2)).
\end{split}
\end{equation*} 
Using \eqref{sp7.8},
\begin{equation*}
\|f_{k_2,l_2,>J_0}\|_{L^2}\leq C2^{-J_0/2}\|f_{k_2,l_2}\|_{X_{k_2}+Y_{k_2}}.
\end{equation*}
Thus, using the definitions, Lemma \ref{Lemmaa1} (a), (c), and \eqref{bj3} we estimate
\begin{equation}\label{bj11}
\begin{split}
2^k\big|\big|&\chi_k(\xi)\cdot(\tau-\omega(\xi,\mu)+i)^{-1}\eta_{\leq J_0}(\tau-\omega(\xi,\mu))\cdot (f_{k_1,l_1}\ast
f_{k_2,l_2,>J_0})\big|\big|_{Y_k}\\
&\leq C2^{k/2}\cdot ||\mathcal{F}^{-1}(f_{k_1,l_1}\ast f_{k_2,l_2,>J_0})||_{L^1_yL^2_{x,t}}\\
&\leq C2^{k/2}||\mathcal{F}^{-1}(f_{k_1,l_1})||_{L^2_yL^\infty_{x,t}}\cdot ||\mathcal{F}^{-1}(f_{k_2,l_2,>J_0})||_{L^2_yL^2_{x,t}}\\
&\leq C2^{3k_1/4}2^{-l_1/4}||f_{k_1}||_{V_{k_1}}\cdot \|f_{k_2,l_2}\|_{X_{k_2}+Y_{k_2}}.
\end{split}
\end{equation} 
Also, using Lemma \ref{Lemmaa1}(a), Lemma \ref{Lemmav2} (a), and the definitions
\begin{equation*}
\|\mathcal{F}^{-1}(f_{k_1,l_1,>J_0})\|_{L^2_yL^\infty_{x,t}}\leq C2^{-J_0/2}2^{(2l_1+k_1)/4}\cdot 2^{k_1-l_1}||f_{k_1}||_{V_{k_1}}.
\end{equation*}
An estimate similar to \eqref{bj11} then gives
\begin{equation}\label{bj12}
\begin{split}
2^k\big|\big|&\chi_k(\xi)\cdot(\tau-\omega(\xi,\mu)+i)^{-1}\eta_{\leq J_0}(\tau-\omega(\xi,\mu))\cdot (f_{k_1,l_1,>J_0}\ast
f^\pm_{k_2,l_2,\leq J_0})\big|\big|_{Y_k}\\
&\leq C2^{3k_1/4}2^{-l_1/4}||f_{k_1}||_{V_{k_1}}\cdot \|f_{k_2,l_2}\|_{X_{k_2}+Y_{k_2}}.
\end{split}
\end{equation} 

It remains to estimate $2^k\big|\big|\chi_k(\xi)\cdot(\tau-\omega(\xi,\mu)+i)^{-1}\eta_{\leq J_0}(\tau-\omega(\xi,\mu))\cdot (f_{k_1,l_1,\leq J_0}\ast
f^\pm_{k_2,l_2,\leq J_0})\big|\big|_{Y_k}$. We will use Lemma \ref{Lemmav1} (b) and Lemma \ref{Lemmav2} (b) to exploit some additional orthogonality. By symmetry, it suffices to prove that
\begin{equation}\label{bj40}
\begin{split}
2^k\big|\big|&\chi_k(\xi)\cdot(\tau-\omega(\xi,\mu)+i)^{-1}\eta_{\leq J_0}(\tau-\omega(\xi,\mu))\cdot (f_{k_1,l_1,\leq J_0}\ast
f^+_{k_2,l_2,\leq J_0})\big|\big|_{Y_k}\\
&\leq C2^{k_1}2^{-l_1/4}||f_{k_1}||_{V_{k_1}\cap W_{k_1}}\cdot \|f_{k_2,l_2}\|_{X_{k_2}+Y_{k_2}}.
\end{split}
\end{equation} 
For $(\xi_1,\mu_1),(\xi_2,\mu_2)\in\mathbb{R}^2$ recall that (see \eqref{bj13})
\begin{equation*}
\begin{split}
&\Omega[(\xi_1,\mu_1),(\xi_2,\mu_2)]=-\omega(\xi_1+\xi_2,\mu_1+\mu_2)+\omega(\xi_1,\mu_1)+\omega(\xi_2,\mu_2)=-\frac{\xi_1\xi_2}{\xi_1+\xi_2}\\
&\times[(\sqrt{3}\xi_1-\mu_1/ \xi_1)+(\sqrt{3}\xi_2+\mu_2/ \xi_2)]\cdot [(\sqrt{3}\xi_1+\mu_1/ \xi_1)+(\sqrt{3}\xi_2-\mu_2/ \xi_2)].
\end{split}
\end{equation*}
Thus, for $\xi_2\in I_{k_2}$, $\mu_2\in[2^{2k-11},2^{2k+11}]$, $\xi_1\in I_{k_1}$, and $|\mu_1|\leq 2^{k+2k_1-9}$
\begin{equation}\label{bj20}
|\Omega[(\xi_1,\mu_1),(\xi_2,\mu_2)]|\geq 2^{k+k_1-4}|(\sqrt{3}\xi_1+\mu_1/ \xi_1)+(\sqrt{3}\xi_2-\mu_2/ \xi_2)|.
\end{equation}

Let $\varphi:\mathbb{R}\to[0,1]$ denote a smooth function supported in $[-1,1]$ with the property that
\begin{equation*}
\sum_{m\in\mathbb{Z}}\varphi(s-m)\equiv 1.
\end{equation*}
Let $\epsilon=2^{-l_1/2}$. For $m\in\mathbb{Z}$ we define
\begin{equation}\label{bj131}
\begin{cases}
&f_{k_1,l_1,\leq J_0}^{+,m}(\xi_1,\mu_1,\tau_1)=f_{k_1,l_1,\leq J_0}(\xi_1,\mu_1,\tau_1)\cdot \varphi((\sqrt{3}\xi_1+\mu_1/ \xi_1)/ \epsilon-m);\\
&f_{k_2,l_2,\leq J_0}^{+,m}(\xi_2,\mu_2,\tau_2)=f_{k_2,l_2,\leq J_0}^+(\xi_2,\mu_2,\tau_2)\cdot \varphi((\sqrt{3}\xi_2-\mu_2/ \xi_2)/ \epsilon+m).
\end{cases}
\end{equation}
The important observation is that, in view of \eqref{bj20} and the definition of $J_0$,
\begin{equation*}
\eta_{\leq J_0}(\tau-\omega(\xi,\mu))\cdot (f_{k_1,l_1,\leq J_0}^{+,m}\ast
f_{k_2,l_2,\leq J_0}^{+,m'})\equiv 0\text{ unless }|m-m'|\leq 4.
\end{equation*}
Thus, using the definitions and Lemma \ref{Lemmaa1} (c),
\begin{equation}\label{bj30}
\begin{split}
&2^k\big|\big|\chi_k(\xi)\cdot(\tau-\omega(\xi,\mu)+i)^{-1}\eta_{\leq J_0}(\tau-\omega(\xi,\mu))\cdot (f_{k_1,l_1,\leq J_0}\ast
f_{k_2,l_2,\leq J_0}^+)\big|\big|_{Y_k}\\
&\leq \negmedspace \sum_{|m-m'|\leq  4}\negmedspace 2^k\big|\big|\chi_k(\xi)(\tau-\omega(\xi,\mu)+i)^{-1}\eta_{\leq J_0}(\tau-\omega(\xi,\mu))\cdot (f_{k_1,l_1,\leq J_0}^{+,m}\ast
f_{k_2,l_2,\leq J_0}^{+,m'})\big|\big|_{Y_k}\\
&\leq C\sum_{|m-m'|\leq  4}2^{k/2}\|\mathcal{F}^{-1}(f_{k_1,l_1,\leq J_0}^{+,m})\|_{L^1_yL^\infty_{x,t}}\cdot\|\mathcal{F}^{-1}(f_{k_2,l_2,\leq J_0}^{+,m'})\|_{L^\infty_yL^2_{x,t}}.
\end{split}
\end{equation}
Using the bound \eqref{bj120}, Lemma \ref{Lemmav2} (b), and the definitions
\begin{equation*}
\begin{split}
\|\mathcal{F}^{-1}(f_{k_1,l_1,\leq J_0}^{+,m})&\|_{L^1_yL^\infty_{x,t}}\leq C2^{l_1/4}2^{k_1/2}\| (\tau_1-\omega(\xi_1,\mu_1)+i)\cdot (I-\partial_{\tau_1}^2)f_{k_1,l_1,\leq J_0}^{+,m}\|_{L^2}^{1/2}\\
&\times \| (\tau_1-\omega(\xi_1,\mu_1)+i)\cdot (I-\partial_{\tau_1}^2)(\partial_{\mu_1}+I)f_{k_1,l_1,\leq J_0}^{+,m}\|_{L^2}^{1/2}.
\end{split}
\end{equation*}
Thus
\begin{equation*}
\begin{split}
\Big[\sum_{m\in\mathbb{Z}}&\|\mathcal{F}^{-1}(f_{k_1,l_1,\leq J_0}^{+,m})\|_{L^1_yL^\infty_{x,t}}^2\Big]^{1/2}\\
&\leq C2^{l_1/4}2^{k_1/2}2^{-(l_1-k_1)/2}\|(I-\partial_{\tau}^2)f_{k_1,l_1}\|_{V_{k_1}}^{1/2}\cdot \|(I-\partial_{\tau}^2)f_{k_1,l_1}\|_{V_{k_1}\cap W_{k_1}}^{1/2}.
\end{split}
\end{equation*}
We substitute this last bound into \eqref{bj30} and, using Lemma \ref{Lemmav1} (b) and \eqref{io1}, we conclude that the right-hand side of \eqref{bj30} is dominated by
\begin{equation*}
\begin{split}
C2^{k/2}\Big[\sum_{m\in\mathbb{Z}}\|\mathcal{F}^{-1}(f_{k_1,l_1,\leq J_0}^{+,m})\|_{L^1_yL^\infty_{x,t}}^2\Big]^{1/2}\cdot \Big[\sum_{m\in\mathbb{Z}}\|\mathcal{F}^{-1}(f_{k_2,l_2,\leq J_0}^{+,m})\|_{L^\infty_yL^2_{x,t}}^2\Big]^{1/2}\\
\leq C2^{-l_1/4}2^{k_1}\|f_{k_1}\|_{V_{k_1}\cap W_{k_1}}\cdot \|f_{k_2,l_2}\|_{X_{k_2}+Y_{k_2}}.
\end{split}
\end{equation*}
This gives the bound \eqref{bj40}. The bound \eqref{bj125} follows from the bounds \eqref{bj10}, \eqref{bj11}, \eqref{bj12}, and \eqref{bj40}.

We estimate now the contribution of $f_{k_1,l_1}\ast f_{k_2,l_2}$, $l_1\in[k+2k_1-10,3k]\cap [1,\infty)$: using \eqref{sp7.7}, Lemma \ref{Lemmaa1} (a), Lemma \ref{Lemmav1} (a), and \eqref{bj3},
\begin{equation}\label{bj6}
\begin{split}
\big|\big|&\chi_k(\xi)\cdot(2^k+i\mu/ 2^k)(\tau-\omega(\xi,\mu)+i)^{-1}\cdot (f_{k_1,l_1}\ast
f_{k_2,l_2})\big|\big|_{X_k}\\
&\leq C(2^k+2^{l_1-k})\cdot ||f_{k_1,l_1}\ast f_{k_2,l_2}||_{L^2_{\xi,\mu,\tau}}\\
&\leq C(2^{k}+2^{l_1-k})\cdot ||\mathcal{F}^{-1}(f_{k_1,l_1})||_{L^2_yL^\infty_{x,t}}\cdot ||\mathcal{F}^{-1}(f_{k_2,l_2})||_{L^\infty_yL^2_{x,t}}\\
&\leq C2^{k_1/4}(2^{-(l_1-k-2k_1)/2}+2^{-(3k-l_1)/2})||f_{k_1}||_{V_{k_1}}\cdot \|f_{k_2,l_2}\|_{X_{k_2}+Y_{k_2}}.
\end{split}
\end{equation}

Finally, we estimate the contribution of  $\sum_{l_1\geq 3k}f_{k_1,l_1}\ast f_{k_2,l_2}$: using \eqref{sp7.7} and Lemma \ref{Lemmav5}
\begin{equation}\label{bj6.6}
\begin{split}
\big|\big|&\chi_k(\xi)\cdot (2^k+i\mu/2^k)(\tau-\omega(\xi,\mu)+i)^{-1}\cdot (\sum_{l_1\geq 3k}f_{k_1,l_1}\ast
f_{k_2,l_2})\big|\big|_{X_k}\\
&\leq C2^{-k}\big|\big|\chi_k(\xi)\cdot\mu \cdot (\sum_{l_1\geq 3k}f_{k_1,l_1}\ast f_{k_2,l_2})\big|\big|_{L^2}\\
&\leq C2^{-k}\big[\sum_{l_1\geq 3k}||2^{l_1}f_{k_1,l_1}\ast f_{k_2,l_2}||_{L^2}^2\big]^{1/2}\\
&\leq C\big[\sum_{l_1\geq 3k}||2^{l_1-k}f_{k_1,l_1}||_{X_{k_1}}^2\cdot ||f_{k_2,l_2}||_{X_k+Y_k}^2\big]^{1/2}\\
&\leq C2^{k_1-k}||f_{k_1}||_{V_{k_1}}\cdot \|f_{k_2,l_2}\|_{X_{k_2}+Y_{k_2}}.
\end{split}
\end{equation}

The main bound \eqref{bj2.1} follows from \eqref{bj5}, \eqref{bj125}, \eqref{bj6}, and \eqref{bj6.6}.
\end{proof}

\newtheorem{Lemmac3}[Lemmac1]{Lemma}
\begin{Lemmac3}\label{Lemmac3}
With the notation in Proposition \ref{Lemmac1}, 
\begin{equation*}
\begin{split}
2^k\big|\big|\chi_k(\xi)\cdot&(\tau-\omega(\xi,\mu)+i)^{-1}\cdot (f_{k_1}\ast
f_{k_2,2k_2-10})\big|\big|_{V_k\cap W_k}\\
&\leq C2^{k_1/4}||f_{k_1}||_{V_{k_1}}\cdot ||f_{k_2,2k_2-10}||_{V_{k_2}\cap W_{k_2}}.
\end{split}
\end{equation*}
\end{Lemmac3}

\begin{proof}[Proof of Lemma \ref{Lemmac3}] As in Lemma \ref{Lemmac2}, it suffices to prove that
\begin{equation}\label{bo1}
\begin{split}
\big|\big|\chi_k(\xi)\cdot(2^{k}+&i\mu/2^k)(\tau-\omega(\xi,\mu)+i)^{-1}\cdot (f_{k_1}\ast
f_{k_2,2k_2-10})\big|\big|_{X_k}\\
&\leq C2^{k_1/4}||f_{k_1}||_{V_{k_1}}\cdot ||f_{k_2,2k_2-10}||_{X_{k_2}+Y_{k_2}}.
\end{split}
\end{equation}

We estimate first the contribution of $f_{k_1,l_1}\ast f_{k_2,2k_2-10}$, $l_1\in[k_1,2k+10]\cap\Z$. Using \eqref{sp7.8}, for any $J\in\Z\cap[-1,\infty)$,
\begin{equation}\label{bo2}
\begin{split}
||\mathcal{F}^{-1}(f_{k_1,l_1,>J})]||_{L^\infty}&\leq C\sum_{j> J}2^{j/2}2^{k_1/2}2^{l_1/2}||f_{k_1,l_1,j}||_{L^2}\\
&\leq C2^{-J/2}2^{k_1/2}2^{l_1/2}||f_{k_1,l_1}||_{X_{k_1}}\\
&\leq C2^{-J/2}2^{3k_1/2}2^{-l_1/2}||f_{k_1}||_{V_{k_1}}.
\end{split}
\end{equation}
Let
\begin{equation}\label{bo3}
J_0=2k+k_1-40.
\end{equation}
Using \eqref{sp7.7}, \eqref{sp7.8}, and \eqref{bo2} (with $J=-1$), we estimate
\begin{equation}\label{rr1}
\begin{split}
2^k\big|\big|\chi_k(\xi)&\cdot(\tau-\omega(\xi,\mu)+i)^{-1}\eta_{\geq J_0+1}(\tau-\omega(\xi,\mu))\cdot (f_{k_1,l_1}\ast f_{k_2,2k_2-10})\big|\big|_{X_k}\\
&\leq C2^k2^{-J_0/2}\big|\big|f_{k_1,l_1}\ast f_{k_2,2k_2-10}\big|\big|_{L^2}\\
&\leq C2^k2^{-J_0/2}\cdot ||\mathcal{F}^{-1}(f_{k_1,l_1})||_{L^\infty}\cdot ||f_{k_2,2k_2-10}||_{L^2}\\
&\leq C2^{k_1-l_1/2}||f_{k_1}||_{V_{k_1}}\cdot ||f_{k_2,2k_2-10}||_{X_{k_2}+Y_{k_2}}.
\end{split}
\end{equation}
Similarly, using \eqref{sp7.7}, \eqref{sp7.8}, and \eqref{bo2}, we estimate
\begin{equation}\label{rr2}
\begin{split}
&2^k\big|\big|\chi_k(\xi)\cdot(\tau-\omega(\xi,\mu)+i)^{-1}\eta_{\leq J_0}(\tau-\omega(\xi,\mu))\cdot (f_{k_1,l_1,\leq J_0}\ast f_{k_2,2k_2-10,>J_0})\big|\big|_{X_k}\\
&+2^k\big|\big|\chi_k(\xi)\cdot(\tau-\omega(\xi,\mu)+i)^{-1}\eta_{\leq J_0}(\tau-\omega(\xi,\mu))\cdot (f_{k_1,l_1,>J_0}\ast f_{k_2,2k_2-10})\big|\big|_{X_k}\\
&\leq C2^{k_1-l_1/2}||f_{k_1}||_{V_{k_1}}\cdot ||f_{k_2,2k_2-10}||_{X_{k_2}+Y_{k_2}}.
\end{split}
\end{equation}

We observe now that
\begin{equation*}
\eta_{\leq J_0}(\tau-\omega(\xi,\mu))\cdot (f_{k_1,l_1,\leq J_0}\ast
f_{k_2,2k_2-10,\leq J_0})\equiv 0,
\end{equation*}
unless $l_1\in[k+k_1-10,k+k_1+10]\cap\Z$, which is a consequence of the identity \eqref{bj13}. As in \eqref{rr1} and \eqref{rr2}, we estimate
\begin{equation}\label{rr5}
\begin{split}
2^k\big|\big|\chi_k(\xi)&\cdot(\tau-\omega(\xi,\mu)+i)^{-1}\eta_{\leq  J_0}(\tau-\omega(\xi,\mu))\cdot (f_{k_1,l_1,\leq J_0}\ast f_{k_2,2k_2-10,\leq J_0})\big|\big|_{X_k}\\
&\leq C2^k\big|\big|f_{k_1,l_1,\leq J_0}\ast f_{k_2,\leq 2k_2-10,\leq J_0}\big|\big|_{L^2}\\
&\leq C2^k2^{3k_1/2}2^{-l_1/2}||f_{k_1}||_{V_{k_1}}\cdot ||f_{k_2,2k_2-10}||_{X_{k_2}+Y_{k_2}}.
\end{split}
\end{equation}
Using Corollary  \ref{Main9}, Lemma \ref{Lemmaa1} (b), and the definitions, we estimate
\begin{equation}\label{rr3}
\begin{split}
2^k&\big|\big|\chi_k(\xi)\cdot(\tau-\omega(\xi,\mu)+i)^{-1}\eta_{\leq J_0}(\tau-\omega(\xi,\mu))\cdot (f_{k_1,l_1,\leq J_0}\ast
f_{k_2,2k_2-10,\leq J_0})\big|\big|_{X_k}\\
&\leq C2^{k}\sum_{j,j_1,j_2=0}^{J_0}2^{-j/2}\big|\big|\eta_{j}(\tau-\omega(\xi,\mu))\cdot (f_{k_1,l_1,j_1}\ast
f_{k_2,2k_2-10,j_2})\big|\big|_{L^2}\\
&\leq  C2^{k}\sum_{j,j_1,j_2=0}^{J_0}2^{-(2k+k_1)/2}\cdot 2^{j_1/2}\|f_{k_1,l_1,j_1}\|_{L^2}\cdot 2^{j_2/2}\|f_{k_2,2k_2-10,j_2}\|_{L^2}\\
&\leq C2^{-k_1/2}\cdot k^3\cdot 2^{-(l_1-k_1)}||f_{k_1}||_{V_{k_1}}\cdot \|f_{k_2,2k_2-10}\|_{X_{k_2}+Y_{k_2}}.
\end{split}
\end{equation}
It follows from \eqref{rr5} and \eqref{rr3} that 
\begin{equation*}
\begin{split}
2^k&\big|\big|\chi_k(\xi)\cdot(\tau-\omega(\xi,\mu)+i)^{-1}\eta_{\leq J_0}(\tau-\omega(\xi,\mu))\cdot (f_{k_1,l_1,\leq J_0}\ast
f_{k_2,2k_2-10,\leq J_0})\big|\big|_{X_k}\\
&\leq C2^{k_1/4}||f_{k_1}||_{V_{k_1}}\cdot \|f_{k_2,2k_2-10}\|_{X_{k_2}+Y_{k_2}},
\end{split}
\end{equation*}
for $l_1\in[k+k_1-10,k+k_1+10]\cap\Z$.
Thus, using also \eqref{rr1} and \eqref{rr2},
\begin{equation}\label{bo4}
\begin{split}
\sum_{l_1=k_1}^{2k+10}2^k\big|\big|&\chi_k(\xi)\cdot(\tau-\omega(\xi,\mu)+i)^{-1}\cdot (f_{k_1,l_1}\ast
f_{k_2,2k_2-10})\big|\big|_{X_k}\\
&\leq C2^{k_1/4}||f_{k_1}||_{V_{k_1}}\cdot ||f_{k_2,2k_2-10}||_{X_{k_2}+Y_{k_2}}.
\end{split}
\end{equation}

We estimate now the contribution of  $\sum_{l_1\geq 2k+11}f_{k_1,l_1}\ast f_{k_2,2k_2-10}$: using \eqref{sp7.7} and Lemma \ref{Lemmav5}, we estimate as in \eqref{bj6.6}
\begin{equation}\label{rr9}
\begin{split}
\big|\big|\chi_k(\xi)\cdot &(2^k+i\mu/2^k)(\tau-\omega(\xi,\mu)+i)^{-1}\cdot (\sum_{l_1\geq 2k+11}f_{k_1,l_1}\ast
f_{k_2,2k_2-10})\big|\big|_{X_k}\\
&\leq C2^{k_1-k}||f_{k_1}||_{V_{k_1}}\cdot \|f_{k_2,\leq 2k_2-10}\|_{X_{k_2}+Y_{k_2}}.
\end{split}
\end{equation}

The main bound \eqref{bo1} follows from \eqref{bo4} and \eqref{rr9}.
\end{proof}

\newtheorem{Lemmac4}[Lemmac1]{Lemma}
\begin{Lemmac4}\label{Lemmac4}
With the notation in Proposition \ref{Lemmac1}, for any $l_2\geq 2k_2+10$
\begin{equation*}
\begin{split}
2^k\big|&\big|\chi_k(\xi)\cdot(\tau-\omega(\xi,\mu)+i)^{-1}\cdot (f_{k_1}\ast
f_{k_2,l_2})\big|\big|_{V_k\cap W_k}\\
&\leq C2^{-(l_2-2k_2)/4}2^{k_1/4}||f_{k_1}||_{V_{k_1}}\cdot ||f_{k_2,l_2}||_{V_{k_2}\cap W_{k_2}}.
\end{split}
\end{equation*}
\end{Lemmac4}

\begin{proof}[Proof of Lemma \ref{Lemmac4}] As in Lemma \ref{Lemmak4}, it suffices to prove that
\begin{equation}\label{fo1.1}
\begin{split}
\big|\big|&\chi_k(\xi)\cdot(2^{l_2-k}+i\mu/2^k)(\tau-\omega(\xi,\mu)+i)^{-1}\cdot (f_{k_1}\ast
f_{k_2,l_2})\big|\big|_{X_k}\\
&\leq C2^{(l_2-2k_2)\cdot (3/4)}2^{k_1/4}||f_{k_1}||_{V_{k_1}}\cdot ||f_{k_2,l_2}||_{X_{k_2}+Y_{k_2}}.
\end{split}
\end{equation}
Using Lemma \ref{Lemmav2} and the definitions, for any $J\in[-1,\infty)\cap\Z$, $k_1\leq 0$, and $l_1\geq k_1$, 
\begin{equation}\label{fff1}
\begin{split}
\|\mathcal{F}^{-1}(f_{k_1,l_1,>J})\|_{L^2_yL^\infty_{x,t}}&\leq C2^{-J/2}2^{k_1/4}(2^{l_1/2}+1)\cdot 2^{k_1-l_1}||f_{k_1}||_{V_{k_1}}\\
&\leq C2^{-J/2}2^{k_1/4}2^{-(l_1-k_1)/2}||f_{k_1}||_{V_{k_1}}.
\end{split}
\end{equation}
Recall also the $L^\infty$ estimate \eqref{bo2}
\begin{equation}\label{fff2}
\begin{split}
||\mathcal{F}^{-1}(f_{k_1,l_1,>J})]||_{L^\infty}\leq C2^{-J/2}2^{k_1}2^{-(l_1-k_1)/2}||f_{k_1}||_{V_{k_1}}.
\end{split}
\end{equation}

We estimate first the contribution of $f_{k_1,l_1}\ast f_{k_2,l_2}$ for
\begin{equation}\label{f1}
l_1\in[k_1,l_2+10]\setminus [l_2-k_2+k_1-10,l_2-k_2+k_1+10].
\end{equation}
Let
\begin{equation}\label{fo3.1}
J_0=l_2+k_1-40.
\end{equation}
If $l_2-2k_2+k_1\geq 0$ then we use \eqref{sp7.7}, \eqref{fff1}, and Lemma \ref{Lemmav1} to estimate
\begin{equation*}
\begin{split}
2^{l_2-k}&\big|\big|\chi_k(\xi)\cdot(\tau-\omega(\xi,\mu)+i)^{-1}\eta_{\geq J_0+1}(\tau-\omega(\xi,\mu))\cdot (f_{k_1,l_1}\ast f_{k_2,l_2})\big|\big|_{X_k}\\
&\leq C2^{l_2-k}2^{-k}\big|\big|f_{k_1,l_1}\ast f_{k_2,l_2}\big|\big|_{L^2}\\
&\leq C2^{l_2-2k_2}||\mathcal{F}^{-1}(f_{k_1,l_1})||_{L^2_yL^\infty_{x,t}}\cdot ||\mathcal{F}^{-1}(f_{k_2,l_2})||_{L^\infty_yL^2_{x,t}}\\
&\leq C2^{(l_2-2k_2)/2}\cdot 2^{k_1/4}2^{-(l_1-k_1)/2}||f_{k_1}||_{V_{k_1}}\cdot ||f_{k_2,l_2}||_{X_{k_2}+Y_{k_2}}.
\end{split}
\end{equation*}
If $l_2-2k_2+k_1\leq 0$ then we use \eqref{sp7.7}, \eqref{fff2}, and \eqref{sp7.8} to estimate
\begin{equation}\label{fr1.10}
\begin{split}
2^{l_2-k}&\big|\big|\chi_k(\xi)\cdot(\tau-\omega(\xi,\mu)+i)^{-1}\eta_{\geq J_0+1}(\tau-\omega(\xi,\mu))\cdot (f_{k_1,l_1}\ast f_{k_2,l_2})\big|\big|_{X_k}\\
&\leq C2^{l_2-k}2^{-(l_2+k_1)/2}\big|\big|f_{k_1,l_1}\ast f_{k_2,l_2}\big|\big|_{L^2}\\
&\leq C2^{l_2-k_2}2^{-(l_2+k_1)/2}||\mathcal{F}^{-1}(f_{k_1,l_1})||_{L^\infty}\cdot ||f_{k_2,l_2}||_{L^2}\\
&\leq C2^{(l_2-2k_2)/2}\cdot 2^{k_1/2}2^{-(l_1-k_1)/2}||f_{k_1}||_{V_{k_1}}\cdot ||f_{k_2,l_2}||_{X_{k_2}+Y_{k_2}}.
\end{split}
\end{equation}
Thus
\begin{equation}\label{fr1.1}
\begin{split}
2^{l_2-k}&\big|\big|\chi_k(\xi)\cdot(\tau-\omega(\xi,\mu)+i)^{-1}\eta_{\geq J_0+1}(\tau-\omega(\xi,\mu))\cdot (f_{k_1,l_1}\ast f_{k_2,l_2})\big|\big|_{X_k}\\
&\leq C2^{(l_2-2k_2)/2}\cdot 2^{k_1/4}2^{-(l_1-k_1)/2}||f_{k_1}||_{V_{k_1}}\cdot ||f_{k_2,l_2}||_{X_{k_2}+Y_{k_2}}.
\end{split}
\end{equation}

Using the $L^\infty$ bound \eqref{fff2}, \eqref{sp7.7}, and \eqref{sp7.8}, we estimate
\begin{equation}\label{fr2.1}
\begin{split}
2^{l_2-k}\big|\big|\chi_k(\xi)\cdot&(\tau-\omega(\xi,\mu)+i)^{-1}\eta_{\leq J_0}(\tau-\omega(\xi,\mu))\cdot (f_{k_1,l_1}\ast f_{k_2,l_2,>J_0})\big|\big|_{X_k}\\
&\leq C2^{l_2-k}\big|\big|f_{k_1,l_1}\ast f_{k_2,l_2,>J_0}\big|\big|_{L^2}\\
&\leq C2^{l_2-k}\|\mathcal{F}^{-1}(f_{k_1,l_1})\|_{L^\infty}\cdot \|f_{k_2,l_2,>J_0}\|_{L^2}\\
&\leq C2^{(l_2-2k_2)/2}\cdot 2^{k_1/2}2^{-(l_1-k_1)/2}||f_{k_1}||_{V_{k_1}}\cdot ||f_{k_2,l_2}||_{X_{k_2}+Y_{k_2}}.
\end{split}
\end{equation}
Using \eqref{fff2}, \eqref{sp7.7}, and \eqref{sp7.8},
\begin{equation}\label{fr1.2}
\begin{split}
2^{l_2-k}\big|\big|\chi_k(\xi)&\cdot(\tau-\omega(\xi,\mu)+i)^{-1}\eta_{\leq  J_0}(\tau-\omega(\xi,\mu))\cdot (f_{k_1,l_1,>J_0}\ast f_{k_2,l_2,\leq  J_0})\big|\big|_{X_k}\\
&\leq C2^{l_2-k}\big|\big|f_{k_1,l_1,>J_0}\ast f_{k_2,l_2,\leq J_0}\big|\big|_{L^2}\\
&\leq C2^{l_2-k}||\mathcal{F}^{-1}(f_{k_1,l_1,>J_0})||_{L^\infty}\cdot ||f_{k_2,l_2,\leq J_0}||_{L^2}\\
&\leq C2^{(l_2-2k_2)/2}\cdot 2^{k_1/2}2^{-(l_1-k_1)/2}||f_{k_1}||_{V_{k_1}}\cdot ||f_{k_2,l_2}||_{X_{k_2}+Y_{k_2}}.
\end{split}
\end{equation}
Finally, we observe that for $l_1$ as in \eqref{f1}
\begin{equation*}
\eta_{\leq J_0}(\tau-\omega(\xi,\mu))\cdot (f_{k_1,l_1,\leq J_0}\ast
f_{k_2,l_2,\leq J_0})\equiv 0,
\end{equation*}
which is a consequence of the identity \eqref{bj13}. Thus, for $l_1$ as in \eqref{f1},
\begin{equation}\label{fr3.1}
\begin{split}
2^{l_2-k}\big|\big|\chi_k(\xi)&\cdot(\tau-\omega(\xi,\mu)+i)^{-1}\cdot (f_{k_1,l_1}\ast f_{k_2,l_2})\big|\big|_{X_k}\\
&\leq C2^{(l_2-2k_2)/2}2^{k_1/4}2^{-(l_1-k_1)/2}||f_{k_1}||_{V_{k_1}}\cdot ||f_{k_2,l_2}||_{X_{k_2}+Y_{k_2}}.
\end{split}
\end{equation}

We estimate now the contribution of $f_{k_1,l_1}\ast f_{k_2,l_2}$, for
\begin{equation*}
l_1\in[l_2-k_2+k_1-10,l_2-k_2+k_1+10].
\end{equation*}
Let
\begin{equation*}
J_1=2k+2k_1-40.
\end{equation*}
As in \eqref{fr1.10}, \eqref{fr2.1}, and \eqref{fr1.2}, using also $2^{k_1-l_1}\approx 2^{k_2-l_2}$, we estimate
\begin{equation}\label{fr4.1}
\begin{split}
&2^{l_2-k}\big|\big|\chi_k(\xi)\cdot(\tau-\omega(\xi,\mu)+i)^{-1}\eta_{\geq J_1+1}(\tau-\omega(\xi,\mu))\cdot (f_{k_1,l_1}\ast f_{k_2,l_2})\big|\big|_{X_k}\\
&+2^{l_2-k}\big|\big|\chi_k(\xi)\cdot(\tau-\omega(\xi,\mu)+i)^{-1}\eta_{\leq J_1}(\tau-\omega(\xi,\mu))\cdot (f_{k_1,l_1}\ast f_{k_2,l_2,>J_1})\big|\big|_{X_k}\\
&+2^{l_2-k}\big|\big|\chi_k(\xi)\cdot (\tau-\omega(\xi,\mu)+i)^{-1}\eta_{\leq J_1}(\tau-\omega(\xi,\mu))\cdot (f_{k_1,l_1,>J_1}\ast f_{k_2,l_2,\leq  J_1})\big|\big|_{X_k}\\
&\leq C 2^{l_2-k}(2^{J_1}+1)^{-1/2}2^{k_1}2^{-(l_1-k_1)/2}||f_{k_1}||_{V_{k_1}}\cdot ||f_{k_2,l_2}||_{X_{k_2}+Y_{k_2}}\\
&\leq C 2^{(l_2-2k_2)/2}2^{k_1/2}||f_{k_1}||_{V_{k_1}}\cdot ||f_{k_2,l_2}||_{X_{k_2}+Y_{k_2}}.
\end{split}
\end{equation}
In addition, using Corollary \ref{Main9}, Lemma \ref{Lemmaa1} (b), and the definitions, we estimate
\begin{equation}\label{fr4.2}
\begin{split}
2^{l_2-k}&\big|\big|\chi_k(\xi)\cdot(\tau-\omega(\xi,\mu)+i)^{-1}\eta_{\leq J_1}(\tau-\omega(\xi,\mu))\cdot (f_{k_1,l_1,\leq J_1}\ast
f_{k_2,l_2,\leq J_1})\big|\big|_{X_k}\\
&\leq C2^{l_2-k}\sum_{j,j_1,j_2=0}^{J_1}2^{-j/2}\big|\big|\eta_{j}(\tau-\omega(\xi,\mu))\cdot (f_{k_1,l_1,j_1}\ast
f_{k_2,l_2,j_2})\big|\big|_{L^2}\\
&\leq  C2^{l_2-k}\sum_{j,j_1,j_2=0}^{J_1}2^{-(2k+k_1)/2}\cdot 2^{j_1/2}\|f_{k_1,l_1,j_1}\|_{L^2}\cdot 2^{j_2/2}\|f_{k_2,l_2,j_2}\|_{L^2}\\
&\leq C2^{k_1/4}||f_{k_1}||_{V_{k_1}}\cdot \|f_{k_2,l_2}\|_{X_{k_2}+Y_{k_2}},
\end{split}
\end{equation}
since we may assume that $k+k_1\geq 0$ (compare with the definition of $J_1$).

We estimate now the contribution of  $\sum_{l_1\geq l_2+11}f_{k_1,l_1}\ast f_{k_2,l_2}$: using \eqref{sp7.7} and Lemma \ref{Lemmav5}, we estimate as in \eqref{gr9.1}
\begin{equation}\label{fr9.1}
\begin{split}
\big|\big|\chi_k(\xi)\cdot &(2^{l_2-k}+i\mu/2^k)(\tau-\omega(\xi,\mu)+i)^{-1}\cdot (\sum_{l_1\geq l_2+11}f_{k_1,l_1}\ast
f_{k_2,l_2})\big|\big|_{X_k}\\
&\leq C2^{k_1-k}||f_{k_1}||_{V_{k_1}}\cdot \|f_{k_2,l_2}\|_{X_{k_2}+Y_{k_2}}.
\end{split}
\end{equation}

The main bound \eqref{fo1.1} follows from \eqref{fr3.1}, \eqref{fr4.1}, \eqref{fr4.2}, and \eqref{fr9.1}.
\end{proof}

\section{Dyadic bilinear estimates III}\label{bilinear3}

In this section we prove the bound \eqref{BIG1} for $k\leq 40$.

\newtheorem{Lemmaj1}{Proposition}[section]
\begin{Lemmaj1}\label{Lemmaj1}
Assume $k\leq 40$, $k_2\in[k-2,k+2]$, $k_1\leq k-20$, $f_{k_1}\in  V_{k_1}\cap W_{k_1}$, and $f_{k_2}\in  V_{k_2}\cap W_{k_2}$.
Then
\begin{equation}\label{ss1}
\begin{split}
2^k\big|\big|\chi_k(\xi)&\cdot(\tau-\omega(\xi,\mu)+i)^{-1}\cdot (f_{k_1}\ast
f_{k_2})\big|\big|_{V_k\cap W_k}\\
&\leq
C2^{k_1/2}||f_{k_1}||_{V_{k_1}\cap W_{k_1}}\cdot ||f_{k_2}||_{V_{k_2}\cap W_{k_2}}.
\end{split}
\end{equation}
\end{Lemmaj1}

\begin{proof}[Proof of Proposition \ref{Lemmaj1}] We show first that
\begin{equation}\label{ss2}
\begin{split}
2^k\big|\big|\chi_k(\xi)&\cdot(\tau-\omega(\xi,\mu)+i)^{-1}\cdot (f_{k_1}\ast
f_{k_2})\big|\big|_{V_k}\leq
C2^{k_1/2}||f_{k_1}||_{V_{k_1}}\cdot ||f_{k_2}||_{V_{k_2}}.
\end{split}
\end{equation}
Using \eqref{sp7.7} and the definition \eqref{sp1}, the left-hand side of \eqref{ss2} is dominated by
\begin{equation}\label{ss3}
\begin{split}
&C||(1+|\mu|)\cdot \chi_k(\xi)\cdot (f_{k_1}\ast f_{k_2})\big|\big|_{L^2}\\
&\leq C||(|\mu_1f_{k_1}|)\ast |f_{k_2}|\,||_{L^2}+C||\,|f_{k_1}|\ast [(1+|\mu_2|)|f_{k_2}|]||_{L^2}.
\end{split}
\end{equation}
We observe now that, for $i=1,2$
\begin{equation}\label{ss4}
\begin{split}
||\mathcal{F}^{-1}(|f_{k_i}|)||_{L^\infty}&\leq C\sum_{l,j\geq 0}||f_{k_i}\cdot \eta_{l}(\mu)\cdot \eta_j(\tau-\omega(\xi,\mu))||_{L^1}\\
&\leq C\sum_{l,j\geq 0}2^{(k_i+l+j)/2}||f_{k_1}\cdot \eta_{l}(\mu)\cdot \eta_j(\tau-\omega(\xi,\mu))||_{L^2}\\
&\leq C2^{k_i/2}\sum_{j\geq 0}2^{j/2}||f_{k_1}\cdot (1+|\mu|)\cdot \eta_j(\tau-\omega(\xi,\mu))||_{L^2}\\
&\leq C2^{k_i/2}||f_{k_i}||_{V_{k_i}}.
\end{split}
\end{equation}
Thus, using also \eqref{sp7.8}, the right hand side of \eqref{ss3} is bounded by
\begin{equation*}
\begin{split}
&C||\mu_1\cdot f_{k_1}||_{L^2}\cdot ||\mathcal{F}^{-1}(|f_{k_2}|)||_{L^\infty}+C||\mathcal{F}^{-1}(|f_{k_1}|)||_{L^\infty}||(1+|\mu_2|)f_{k_2}||_{L^2}\\
&\leq C2^{k_1}||(\mu_1/2^{k_1})\cdot f_{k_1}||_{X_{k_1}}||f_{k_2}||_{V_{k_2}}+C2^{k_1/2}||f_{k_1}||_{V_{k_1}}||(1+|\mu_2|)f_{k_2}||_{X_{k_2}}\\
&\leq C2^{k_1/2}||f_{k_1}||_{V_{k_1}}\cdot ||f_{k_2}||_{V_{k_2}}.
\end{split}
\end{equation*}
This completes the proof of \eqref{ss2}.

We show now that
\begin{equation}\label{ss8}
\begin{split}
2^k\big|\big|\chi_k(\xi)&\cdot(\tau-\omega(\xi,\mu)+i)^{-1}\cdot (f_{k_1}\ast
f_{k_2})\big|\big|_{W_k}\leq
C2^{k_1/2}||f_{k_1}||_{V_{k_1}}||f_{k_2}||_{V_{k_2}\cap W_{k_2}}.
\end{split}
\end{equation}
In view of the definition \eqref{sp4}, the left-hand side of \eqref{ss8} is bounded by
\begin{equation}\label{ss9}
\begin{split}
&C2^k\big|\big|\chi_k(\xi)\cdot(\mu/\xi)(\tau-\omega(\xi,\mu)+i)^{-2}\cdot (f_{k_1}\ast
f_{k_2})\big|\big|_{X_k}\\
&+C2^k\big|\big|\chi_k(\xi)\cdot(\tau-\omega(\xi,\mu)+i)^{-1}\cdot (f_{k_1}\ast
(\partial_\mu+I)f_{k_2})\big|\big|_{X_k}.
\end{split}
\end{equation}
The first term in \eqref{ss9} is dominated by the left-hand side of \eqref{ss3}. We use \eqref{sp7.7}, \eqref{ss4}, and  \eqref{sp7.8}  to estimate the second term in \eqref{ss9} by
\begin{equation*}
\begin{split}
C2^k\big|\big|f_{k_1}\ast(\partial_\mu+I)f_{k_2}\big|\big|_{L^2}&\leq C||\mathcal{F}^{-1}(f_{k_1})||_{L^\infty}\cdot ||(\partial_\mu+I)f_{k_2}||_{L^2}\\
&\leq C2^{k_1/2}||f_{k_1}||_{X_{k_1}}\cdot ||(\partial_\mu+I)f_{k_2}||_{X_{k_2}},
\end{split}
\end{equation*}
which suffices for \eqref{ss8}.
\end{proof}

\section{Dyadic bilinear estimates IV}\label{bilinear4}

In this section we prove the bound \eqref{BIG2}.

\newtheorem{Lemmal1}{Proposition}[section]
\begin{Lemmal1}\label{Lemmal1}
Assume $k_1,k_2\in\Z$, $|k_1-k_2|\leq 100$, $f_{k_1}\in  V_{k_1}\cap W_{k_1}$, and $f_{k_2}\in  V_{k_2}\cap W_{k_2}$. Then
\begin{equation}\label{ll1}
\begin{split}
\Big[\sum_{k\in\Z}\big|\big|2^k\chi_k(\xi)&(\tau-\omega(\xi,\mu)+i)^{-1}\cdot (f_{k_1}\ast
f_{k_2})\big|\big|_{V_k\cap W_k}^2\Big]^{1/2}\\
&\leq C||f_{k_1}||_{V_{k_1}\cap W_{k_1}}\cdot ||f_{k_2}||_{V_{k_2}\cap W_{k_2}}.
\end{split}
\end{equation}
\end{Lemmal1}

\begin{proof}[Proof of Proposition \ref{Lemmal1}] We show first that
\begin{equation}\label{ll2}
\begin{split}
\Big[\sum_{k\in\Z}\big|\big|2^k\chi_k(\xi)&(\tau-\omega(\xi,\mu)+i)^{-1}\cdot (f_{k_1}\ast
f_{k_2})\big|\big|_{V_k}^2\Big]^{1/2}\leq C||f_{k_1}||_{V_{k_1}}\cdot ||f_{k_2}||_{V_{k_2}}.
\end{split}
\end{equation}
Using \eqref{sp7.7} and the definition \eqref{sp1}, the left-hand side of \eqref{ll2} is dominated by
\begin{equation}\label{ll3}
\begin{split}
C\Big[\sum_{k\in\Z}&\big|\big|(1+2^{2k}+|\mu|)\chi_k(\xi)\cdot (f_{k_1}\ast
f_{k_2})\big|\big|_{L^2}^2\Big]^{1/2}\\
&\leq C(2^{k_1+k_2}+1)||f_{k_1}\ast
f_{k_2}||_{L^2}+C||\mu\cdot (f_{k_1}\ast f_{k_2})||_{L^2}.
\end{split}
\end{equation}
Using Lemma \ref{Lemmav5}, the first term in the right-hand side of \eqref{ll3} is dominated by
\begin{equation*}
C(2^{k_1}+1)||\mathcal{F}^{-1}(f_{k_1})||_{L^4}\cdot (2^{k_2}+1)||\mathcal{F}^{-1}(f_{k_2})||_{L^4}\leq C||f_{k_1}||_{V_{k_1}}\cdot ||f_{k_2}||_{V_{k_2}}.
\end{equation*}
The second term is bounded by
\begin{equation*}
\begin{split}
&C||\mathcal{F}^{-1}(\mu_1\cdot f_{k_1})||_{L^4}\cdot ||\mathcal{F}^{-1}(f_{k_2})||_{L^4}+C||\mathcal{F}^{-1}(f_{k_1})||_{L^4}\cdot ||\mathcal{F}^{-1}(\mu_2\cdot f_{k_2})||_{L^4}\\
&\leq C2^{k_1}||f_{k_1}||_{V_{k_1}}\cdot 2^{-k_2}||f_{k_2}||_{V_{k_2}}+C2^{-k_1}||f_{k_1}||_{V_{k_1}}\cdot 2^{k_2}||f_{k_2}||_{V_{k_2}}.
\end{split}
\end{equation*}
This completes the proof of \eqref{ll2}.

We show now that
\begin{equation}\label{ll6}
\begin{split}
\Big[\sum_{k\in\Z}\big|\big|2^k\chi_k(\xi)&(\tau-\omega(\xi,\mu)+i)^{-1}\cdot (f_{k_1}\ast
f_{k_2})\big|\big|_{W_k}^2\Big]^{1/2}\leq C||f_{k_1}||_{V_{k_1}}\cdot ||f_{k_2}||_{V_{k_2}\cap W_{k_2}}.
\end{split}
\end{equation}
In view of the definition \eqref{sp4}, the left-hand side of \eqref{ll6} is bounded by
\begin{equation}\label{ll7}
\begin{split}
&C\Big[\sum_{k\in\Z}\big|\big|2^k\chi_k(\xi)\cdot (\mu/ \xi)(\tau-\omega(\xi,\mu)+i)^{-2}\cdot (f_{k_1}\ast
f_{k_2})\big|\big|_{X_k+Y_k}^2\Big]^{1/2}\\
&+C\Big[\sum_{k\in\Z}\big|\big|2^k\chi_k(\xi)(\tau-\omega(\xi,\mu)+i)^{-1}\cdot (f_{k_1}\ast
(\partial_\mu+I)f_{k_2})\big|\big|_{X_k+Y_k}^2\Big]^{1/2}.
\end{split}
\end{equation}
The first term in \eqref{ll7} is dominated by the left-hand side of \eqref{ll2}, which suffices. Using \eqref{sp7.7} and Lemma \ref{Lemmav5}, the second term in \eqref{ll7} is bounded by
\begin{equation*}
\begin{split}
C2^{k_2}||f_{k_1}&\ast (\partial_\mu+I)f_{k_2}||_{L^2}\leq C||\mathcal{F}^{-1}(f_{k_1})||_{L^4}\cdot 2^{k_2}||\mathcal{F}^{-1}((\partial_\mu+I)f_{k_2})||_{L^4}\\
&\leq C||f_{k_1}||_{V_{k_1}}\cdot ||(\partial_\mu+I)f_{k_2}||_{X_{k_2}+Y_{k_2}}.
\end{split}
\end{equation*}
This completes the proof of \eqref{ll6}. 

The proposition follows from the estimates \eqref{ll2} and \eqref{ll6}.
\end{proof}

\end{document}